\documentclass[11pt,a4paper]{article}
\usepackage[margin=20mm,headheight=14pt]{geometry}
\usepackage[T1]{fontenc}
\usepackage{lmodern}
\usepackage{amsmath,amssymb,amsthm,mathtools,mathrsfs,enumitem,booktabs,array,longtable,xcolor}
\usepackage[unicode,colorlinks=true,linkcolor=blue!45!black,citecolor=blue!45!black,urlcolor=blue!45!black]{hyperref}
\usepackage{fancyhdr}
\newtheorem{theorem}{Theorem}[section]
\newtheorem{proposition}[theorem]{Proposition}
\newtheorem{lemma}[theorem]{Lemma}
\newtheorem{corollary}[theorem]{Corollary}
\theoremstyle{definition}
\newtheorem{definition}[theorem]{Definition}

\newcommand{\Z}{\mathbb Z}\newcommand{\Q}{\mathbb Q}\newcommand{\R}{\mathbb R}\newcommand{\C}{\mathbb C}
\newcommand{\A}{\mathbb A}\newcommand{\Sph}{\mathbb S}\newcommand{\F}{\mathbb F}
\newcommand{\cS}{\mathcal S}\newcommand{\cD}{\mathcal D}\newcommand{\cB}{\mathcal B}\newcommand{\cE}{\mathcal E}\newcommand{\cL}{\mathcal L}
\DeclareMathOperator{\Perf}{Perf}\DeclareMathOperator{\Mod}{Mod}\DeclareMathOperator{\Fun}{Fun}
\DeclareMathOperator{\End}{End}\DeclareMathOperator{\Tr}{Tr}\DeclareMathOperator{\cofib}{cofib}
\DeclareMathOperator{\fib}{fib}\DeclareMathOperator{\ord}{ord}
\DeclareMathOperator{\Loc}{Loc}\DeclareMathOperator{\Tot}{Tot}
\DeclareMathOperator{\Res}{Res}
\DeclareMathOperator{\Shv}{Shv}\DeclareMathOperator{\Sp}{Sp}
\DeclareMathOperator{\Cospan}{Cospan}
\DeclareMathOperator{\coker}{coker}

\DeclareMathOperator{\covol}{covol}\DeclareMathOperator{\rank}{rank}
\DeclareMathOperator{\vol}{vol}

\setlist{itemsep=2pt,topsep=4pt}
\hypersetup{pdfauthor={},pdftitle={},pdfsubject={},pdfkeywords={},pdfcreator={},pdfproducer={}}
\title{Arithmetic Winding Curves in $\mathbb F_1$ Geometry}
\author{SeungJu Lee}
\date{\today}
\begin{document}
\maketitle
\begin{abstract}
We construct arithmetic winding curves from cyclic covers and the completions
of a number field $K$. Transfers between cyclic periods give the primitive
cyclotomic relations. Selected factors of the equivariant endomorphism rings
of their derived intersections are finite fields in odd residue characteristic
and length-two Witt rings in residue characteristic two. Valuation covers
carry these periods on the curve $Y_K$. Their normal coordinates realize
deck shifts as Schwartz operators in Meyer's arithmetic trace formula, and
their puncture differentials give the finite Euler factors of $\zeta_K$.
We construct a sheaf on $Y_K$ for which the middle Fourier--Poisson
representation embeds as a $C_K$-module in ordinary first sheaf cohomology. The
Dedekind-zeta component is an equivariant direct summand. The injection
uses valuation monodromy at the finite places and adelic Poisson
summation on the archimedean charts. The cup product of two image classes vanishes with tensor-product
coefficients, also after conjugating one factor. On the section side, tensor
degree and integral jets give multiplicative lattices, explicit duality
quotients, and relative theta identities. A principal cyclotomic twist
with determinant one on truncated jets yields a quadratic lower bound
for the logarithmic index in a family of increasing rank. The
fixed-weight Poisson quotients preserve finite collections of Mellin
derivatives at Hecke zeros of orders below their respective multiplicities.
A normalized pairing defined by relative section determinants recovers
the complete arithmetic trace, including its identity term. For
tests on the norm quotient it is a signed sum of contact theta differences
and regularized Green intersections on the marked coefficient line.
\end{abstract}
\section{Introduction}
\label{sec:introduction}

A closed point of degree $d$ on a curve over $\F_q$ determines both
a residue field of order $q^d$ and a Frobenius orbit of length $d$.
The search for a geometry over $\F_1$ asks, in particular, how residue
fields, prime orbits, and the zeta function could be related on a
curve for a number field \cite[Introduction and \S1.2]{SouleFone}.
Our construction starts with a cyclic cover. Its transfer maps impose
cyclotomic relations, and its derived intersections have endomorphisms
that supply the additive arithmetic of the residue coefficients.
At each finite place $v$ the valuation cover carries these coefficients.
The boundary has period $\log q_v$ after passage to the norm quotient.
The same valuation shift acts on the test functions entering
the explicit formula. Definition~\ref{def:assembly:curve} assembles the
space $Y_K$, its covering towers, its section realizations, and the
action on the Fourier--Poisson complex. Their compatibilities are the
content of Theorem~\ref{thm:assembly:compatibility}.

There is also a comparison on $Y_K$ itself. With
$Q=H^1(C_K^{\mathrm{FP}})$ and $Q_{\mathrm D}=Q^{C_K^1}$, we construct
a sheaf $\mathscr F_K$ of complex vector spaces and equivariant maps
\[
 Q\xrightarrow{\ \iota\ }H^1(Y_K,\mathscr F_K)
   \xrightarrow{\ r_{v,w}\ }Q_{\mathrm D},\qquad r_{v,w}\iota=\Pi_{\mathrm D} ,
\]
where $\Pi_{\mathrm D}$ is the normalized average over the norm-one idele
class group. The finite charts force every extendible puncture
section to vanish. On the archimedean charts the extendible
constant sections are the image of adelic Poisson summation.
The resulting connecting map embeds the whole middle representation.
A winding loop at a finite place and an archimedean arc give a left inverse
on its Dedekind component
(Propositions~\ref{cmp:mv} and~\ref{cmp:prop-retraction}).
The cup product on pairs of these classes is computed in
Proposition~\ref{cmp:cup}.

The normal determinant calculation recovers the factors of the completed
Dedekind zeta function. We use \emph{extended Riemann hypothesis} (ERH)
for the assertion that all its nontrivial zeros have real part $1/2$
for every number field $K$. The Fourier--Poisson representation also
carries the unitary Hecke components, including characters of infinite
order. Section~\ref{sec:zero-comparison} gives fixed-weight completions
preserving evaluations at their zeros with multiplicity. Integral jets
on marked projective coefficient lines provide a second arithmetic
measurement. They detect a principal twist even when its action on
the truncated jet has determinant one.

\subsection{From cyclic periods to arithmetic coefficients}

Write $C_m=\Z/m\Z$ and let $R=\Sph[\Z]$ be the sphere group algebra,
with generator $t$.
The primitive period-$m$ object has the presentation
\[
 U_m\simeq\cofib\bigl(\Phi_m(t):R\longrightarrow R\bigr),
\]
where $\Phi_m$ is the $m$th cyclotomic polynomial. Its integral
degree-zero homology is $\Z[t]/(\Phi_m)$. For $m<n$, the derived
intersection $U_m\otimes_R U_n$ vanishes unless $n/m$ is a prime
power; in that case its underlying spectrum is a wedge of Moore spectra
(Proposition~\ref{prop:src:contact-spectrum}). Retaining the $C_m$-action
on this intersection also retains its cyclotomic coefficient information.
When $p\nmid m$, selected factors of
\[
 \pi_0\End_{C_m}(U_m\otimes_R U_{mp})
\]
are finite fields in odd characteristic and $W_2(k)$ in residue
characteristic two (Theorem~\ref{thm:src:endomorphisms}). The additive
coefficients thus arise from stable maps between geometric period objects.
The distinction between the underlying Moore spectrum and its indexed
endomorphisms is essential to this recovery.

Graph cospans supply the composition law. Sewing is defined by homotopy
pushout and gives a cofiber sequence on the associated coefficients
(Proposition~\ref{prop:src:sewing}). The correspondence coefficient algebra
constructed later preserves the source modules and their mapping spectra
through Amitsur descent (Theorem~\ref{thm:cor:descent}). Section~\ref{sec:arithmetic} also gives an arithmetic
presentation of the input field: a cyclic irreducible marked graph return
determines a number field, and its invariant lattices recover the ring of
integers and the completion data used in the curve.

\subsection{The curve and its arithmetic trace}

Write $C_K=\A_K^\times/K^\times$. For a place $v$, let
$j_v:K_v^\times\to C_K$ be the map induced by the $v$th idele coordinate.
The chart at $v$ is
\[
 N_v=(C_K\times K_v)/K_v^\times,\qquad
 (c,z)\cdot h=(c\,j_v(h),h^{-1}z).
\]
Its puncture is identified with $C_K$ by $(c,z)\mapsto c\,j_v(z)$.
Gluing all the charts along this common puncture defines $Y_K$.
The zero orbit in $N_v$ is $C_K/j_v(K_v^\times)$, so the chart records
both the idele-class action and the multiplicative stabilizer at $v$.

These properties explain the choice of the curve. At a finite place,
quotienting by $\mathcal O_v^\times$ leaves a valuation coordinate and
produces a principal $\Z$-cover $Q_v\to Y_K$. Its finite cyclic
quotients carry the winding periods, and refinements by open subgroups of
$\mathcal O_v^\times$ supply the angular data used to select residue
coefficients. Quotienting the zero orbit $Z_v$ by the norm-one subgroup $C_K^1$ gives the
boundary circle $\R_{>0}/q_v^\Z$, where $q_v$ is the residue-field size.
Thus the norm of a prime determines a period of the
same geometric action that carries the covers. Real and complex
completions enter through the same chart formula, while finite extensions
act through the norm maps of Theorem~\ref{thm:geo:norm}.
For $a=j_v(a_v)\ne1$, the fixed stratum is $Z_v$ and the transverse
coordinate is multiplied by $a_v$. The inverse-pullback contribution is
$|1-a_v^{-1}|_v^{-1}$, the factor occurring in the trace away from
the identity (Theorem~\ref{thm:ana:fixed}). The same charts therefore
carry the covering data and the transverse action used in the explicit
formula.

The connection to the explicit formula uses a Schwartz module realization.
The normal ball basis is sent to Schwartz functions by
\[
 t_v^n\longmapsto\mathbf1_{\varpi_v^n\mathcal O_v}.
\]
This map intertwines multiplication by $t_v$ with dilation by
$\varpi_v$ (Proposition~\ref{prop:ana:ball-module}).
The correspondence algebra acts on the
Fourier--Poisson complex, and Meyer's trace theorem evaluates the action,
including its prime and archimedean contributions
\cite[Theorems~5.8 and~5.11]{Meyer}.
Theorem~\ref{thm:ana:complete-trace} and
Proposition~\ref{prop:ana:source-trace} express this trace with the chosen
measures and source operators. The role of the atlas is therefore concrete:
it supplies coordinates and covering operations that have specified images
in the arithmetic trace calculation. The same puncture differential
$1-t_v$, specialized by $t_v\mapsto q_v^{-s}$, gives the Euler factor
$(1-q_v^{-s})^{-1}$. Its product over finite places, together with the
archimedean Gaussian integrals and the self-dual volume factor, is
\[
 |D_K|^{s/2}\Gamma_\R(s)^{r_1}\Gamma_\C(s)^{r_2}\zeta_K(s)
 \qquad (\Re(s)>1)
\]
(Proposition~\ref{prop:ana:completed-product}), where $D_K$ is the field
discriminant and $r_1,r_2$ count the real and complex places.

\subsection{Sections, integral jets, and metrics}

Section multiplication must respect a change of winding frame. The cover
$Q_v$ defines an invertible $\Z[t,t^{-1}]$-coefficient sheaf
$\mathscr L_v$. A transition $t^n$ acts on its $m$th tensor power by
$t^{mn}$, so tensor degree gives the descended algebra
$\bigoplus_{m\ge0}\mathscr L_v^{\otimes m}$
(Proposition~\ref{prop:prod-descent}). At finite period this algebra has
a frame-invariant coefficient cone and a positive Hermitian metric with
contractive multiplication. Quotient maps are contractions, and the
normalized transfers are isometric multiplicative embeddings with an
idempotent unit image.

To obtain finite section lattices, we realize the coefficient variable $t$
on a marked projective line over $\Z$. This is the projective realization
of the coefficient algebra. A marked divisor $D$ gives a finite free
lattice $S(D)$, and integrality of its jet at $t=1$ through order $M-1$
defines $S_M(D)$. This condition retains arithmetic information at the
marked point that rank and degree alone do not record. These lattices satisfy
\[
 S_M(D)\,S_M(E)\subseteq S_M(D+E),
\]
and multiplication by a nonzero section is injective
(Theorem~\ref{thm:prod-integral-jets}). When $\rank S(D)=M$, their
index is computed by an expansion of fractional lattices
(Theorem~\ref{thm:prod-index-expansion}).

The index detects information that the divisor class and determinant can miss.
Write $Z(P)$ for the divisor cut out by a polynomial $P$.
For the degree-zero divisor $D_0$ and the associated divisors $D_m^\pm$ given by
\[
 D_0=2Z(\Phi_2)-3Z(\Phi_4)+Z(\Phi_8),\qquad
 D_m^\pm=(m-1)[\infty]\pm mD_0,
\]
Corollary~\ref{cor:prod-signed-index} proves
\[
 \log[S(D_m^\pm):S_m(D_m^\pm)]
 \ge\left(\frac18m^2-O(m\log m)\right)\log2.
\]
The ranks here grow with $m$. The corresponding cyclotomic rational function
\[
 R(t)=\frac{\Phi_2(t)^2\Phi_8(t)}{\Phi_4(t)^3}
\]
satisfies
$D_0=\operatorname{div}(R)$ and $R(1)=1$, so
multiplication by $R(1+x)^m$ on $\Q_2[x]/(x^m)$ has determinant one.
The ordinary divisor class and this determinant therefore do not record
the displayed logarithmic index. The retained integral jet does.
Propositions~\ref{prop:prod-dual} and
\ref{prop:prod-jet-duality-quotient} compute the dual lattice and the
finite quotient relating it to the opposite jet condition. Poisson
summation then gives the relative theta identity of
Theorem~\ref{thm:prod-jet-theta}.

The metrics are chosen to work with these operations. The Bombieri metric
gives contractive multiplication across tensor degrees. A weighted Laurent
completion gives bounded division on cyclotomic vanishing ideals and
positive finite jet quotient metrics
(Theorem~\ref{thm:pos-jet-quotient}). Averaging and normalizing its
finite-period quotient metrics recovers the winding-cover metric
(Proposition~\ref{prop:pos-finite-cover-metric}). Thus the positive norms
used for sections and those used for cyclic coefficients are linked by a
proved metric comparison.

The final analytic construction concerns completion of the Poisson
quotient. For each unitary character component and each fixed exponential
weight $\alpha>1/2$, Theorem~\ref{thm:zero:jets} constructs a nuclear
Fr\'echet quotient $Q_\alpha$. If $\rho$ is an actual nontrivial zero
of multiplicity $m_\rho$, its Mellin derivative evaluations of orders
$0\le j<m_\rho$ remain continuous on $Q_\alpha$; any finite set of
distinct such evaluations is jointly surjective. The completion therefore
preserves these finite tests of zeros and their multiplicities.

Section~\ref{sec:section-trace-comparison} constructs a pairing from
relative determinants of two multiplication maps on linear sections.
Its normalization at the identity is fixed by the reference additive
characters and their self-dual Gaussians; degree-zero lattice duality
supplies the discriminant correction. Theorem~\ref{thm:st:comparison}
identifies this pairing with the complete arithmetic trace for all
compactly supported tests, including angular components. On the norm
quotient, Proposition~\ref{prop:st:finite-theta} and
Theorem~\ref{thm:st:green} express the same form through the existing
cyclotomic contact theta differences and a regularized signed sum of
Green intersections on the marked coefficient line.
\subsection{Cyclic coefficients and spherical arithmetic}
\label{sec:relations}

Several established approaches to $\F_1$ already encode cyclic arithmetic.
Deitmar's extension of a monoid $M$ to $\Z$ is the monoid algebra
$\Z[M]$ \cite[\S1, Theorem~1.1]{DeitmarFone}. Soul\'e uses
$\Z[T]/(T^L-1)$ for cyclic scalar extensions and includes Hermitian
lattices among his examples \cite[\S2.4, equation~(4), and \S5.3]{SouleFone}.
The integral realization of our level-$L$ coefficient is
\[
 H\Z\otimes_{\Sph}\Sph[C_L]\simeq H(\Z[C_L]).
\]
The level $\mathbb F_{1^L}$ in Definition~\ref{def:src:f1level}
also retains the covering torsor, primitive spectrum, module category,
quotient maps, and stable transfers. Thus $L$ has the geometric meaning
of a covering degree and a cyclic period.

Lorscheid's cyclotomic blueprint imposes subgroup-sum relations on
the roots of unity. Its extension to rings is $\Z[\zeta_m]$. For
$m=q-1$, the residue field at a prime above $p$ is $\F_q$, where
$q$ is a power of $p$ \cite[\S1.10]{LorscheidBlueprints}. Our transfer
cofiber realizes these relations in $H_0(U_m;\Z)$
(Proposition~\ref{prop:src:primitive-lattice}).
For example, at $m=3$ the transfer imposes $1+t+t^2=0$. Since
$\Phi_6(\zeta_3)=-2\zeta_3$, the underlying spectrum of $T_{3,6}=U_3\otimes_R U_6$ is a
wedge of two copies of $\Sph/2$. Theorems~\ref{thm:src:endomorphisms} and~\ref{thm:res:residue}
compute the ring of $C_3$-equivariant endomorphisms in degree zero,
\[
 \pi_0\End_{C_3}(T_{3,6})
     \simeq (\Z/4)[t]/(t^2+t+1)
     \simeq W_2(\F_4).
\]
Its reduction is $\F_4$ because $t^2+t+1$ is irreducible modulo $2$.
This endomorphism ring has characteristic four because the stable
identity of $\Sph/2$ has additive order four
\cite[Proposition~4.1]{StricklandMoore}. Blueprints describe
relations between formal sums, and bands record null sums in a pointed
multiplicative monoid \cite[Definitions~1.1--1.2]{LorscheidBlueprints}
\cite[Definition~1.1]{BakerJinLorscheidBands}.

The primitive spectra have a direct spectral predecessor.
Construction~3.3.2 of Zhang, following Rezk, takes the fibre of the
quotient map between cyclic permutation spectra
\cite[Construction~3.3.2, pp.~25--26, author v4]{ZhangDirichlet}.
Our transfer cofiber is its equivariant Spanier--Whitehead dual,
as explained after Definition~\ref{def:src:primitive}.
To\"en--Vaqui\'e include
schemes over the sphere spectrum in relative algebraic geometry
\cite[Definitions~5.8--5.9 and \S5.2]{ToenVaquieRelative}.
Carmeli--Nikolaus--Yuan prove that the functor
$G\mapsto\Sph[G]$ on finitely generated abelian groups is fully faithful
as a functor to commutative ring spectra
\cite[Theorem~1.1]{CarmeliNikolausYuan}.
Our primitive cofibers and their contacts are formed in module categories.

For a scheme flat over $\Z$, Borger describes a $\Lambda$-structure
by commuting Frobenius lifts, one for every prime. On a monoid algebra
the toric lifts act by $\psi_\ell([g])=[g^\ell]$
\cite[Introduction and \S2.2]{BorgerLambda}. Thus the integral
level ring $\Z[C_L]$ carries all the substitutions $t\mapsto t^\ell$.
On the primitive quotient $\Z[t]/(\Phi_m)$ this substitution descends
when $\ell\nmid m$, since it permutes the primitive $m$th roots.
For $\ell\mid m$ it sends a primitive root to one of order $m/\ell$,
as illustrated by $\Phi_3(\zeta_3^3)=3$. At the contact above, the map
$t\mapsto t^2$ fixes $\Z/4$ and lifts Frobenius on $\F_4$.
The compatible change of marking on $C_6\to C_3$ has exponent $5$
on $C_6$, as $5\equiv2\bmod 3$ and $5\equiv1\bmod 2$
(Proposition~\ref{prop:res:frobenius}).

Deck translation on the valuation torsor acts by increasing
its integer label by one. After choosing a uniformizer, its generator
$t_v$ maps to the idele-class translation $\delta_{j_v(\varpi_v)}$
in the correspondence algebra. For a finite set of places $S$,
Theorem~\ref{thm:cor:descent} constructs a module retraction of
this coefficient map. Mathew's descent theorem then recovers the
$\Sph[\Z^S]$-modules and their mapping spectra by Amitsur totalization,
also after finite cyclic base change
\cite[Definition~3.18 and Propositions~3.22, 3.28]{Mathew}.

\subsection{Boolean coefficients and Connes--Consani geometry}

The Boolean semifield $\mathbb B$ has idempotent addition $1+1=1$.
A unital semiring map from $\mathbb B$ to a ring therefore forces $1=0$
in the target. Connes and Consani also develop a spherical coefficient
formalism that includes rings \cite[\S3]{ConnesConsaniMetaphysics}.
In our construction the additive law on the residue coefficients comes
from addition of stable maps, and multiplication comes from composition.

Connes and Consani identify the points of their arithmetic site over the
tropical semifield with $\Q^\times\backslash\A_\Q/\widehat{\Z}^{\times}$
\cite[Theorem~1.1 and Definition~3.1]{ConnesConsaniArithmetic}.
Their scaling site gives the orbit $C_p=\R_{>0}/p^\Z$ a sheaf of
continuous convex piecewise affine functions with slopes in $\Z[1/p]$.
They compute its divisor class group and prove a Riemann--Roch formula
using a real-valued dimension
\cite[\S2.1, Lemma~5.1(ii), Theorems~5.6 and~5.17]{ConnesConsaniScaling}.
For $K=\Q$, Proposition~\ref{prop:geo:norm-flow} gives the same circle
in $C_K^1\backslash Y_K$. Before taking that quotient, our valuation
tower records the powers of a uniformizer and its angular refinements
retain finite quotients of the units
(Theorem~\ref{thm:geo:valuation-cover} and Proposition~\ref{prop:res:angular-covers}).

Connes and Consani's spherical approach uses the categorical sphere in
$\Gamma$-sets and $\Gamma$-rings. It supplies divisor cohomologies,
integer dimensions and Serre duality, leading to Riemann--Roch formulas
for $\overline{\operatorname{Spec}\Z}$
\cite[Theorems~1.2--1.3]{CCRiemannRoch2023}
\cite[Theorem~1]{CCRiemannRoch2024}.
Our section theory computes integral-jet lattices on the marked
coefficient realizations. Proposition~\ref{prop:prod-jet-duality-quotient}
identifies the finite quotient between the dual lattice and the opposite
jet condition. Corollary~\ref{cor:prod-signed-index} gives a quadratic
lower bound for the logarithm of the section index along powers of a
principal twist, with the rank increasing.

In their description of the Picard and Jacobian monoids, Connes and
Consani represent adelic classes by rank-one torsion-free groups with
rigidifying data and realize multiplication by tensor product. Their
description of the transverse field relates the stabilizer
$\Q_v^\times$ to the explicit formula
\cite[\S\S2--4, 8, 9.1--9.2]{CCJacobian2026}. On
$\operatorname{Spec}\Z$, the Schwartz algebras for finite sets of places
containing $\infty$ give the sections of one sheaf, and their crossed
products give the sections of a second sheaf \cite[\S9.4, (21) and Theorem~9.1]{CCJacobian2026}.
The absolute curve of
\cite[Proposition~2.2, Corollary~3.15 and Theorem~3.17]{CCAbsoluteCurve2026}
is defined on $\operatorname{Spec}\Z$ by pulling back a spherical
structure sheaf, with stalk $\mathbb F_1[T^{\Z[1/p]_+}]$ at $p$.
For an algebraically closed perfectoid field $F$ of residue characteristic
$p$, the quotient of the continuous geometric points at $p$ by the
$\Q_p^\times$-action parametrizes untilts of $F^\flat$ up to Frobenius twists of the
tilt identification, including the equal-characteristic untilt $F^\flat$ itself.
Their twistor construction supplies an archimedean component over the
signed spherical base $\mathbb F_{1^2}$. They extend the finite curve to
this base and glue the two components at their generic points
\cite[Definitions~1.1 and~7.1]{CCTwistor2026}.

For finite sets of places of $\Q$ containing $\infty$, Connes and Consani
study the union of the generic orbit and the strata with one zero coordinate
\cite[\S\S5.3--5.4]{CCKnots2025}. For every number field $K$, our atlas
glues $C_K\times_{K_v^\times}K_v$ along $C_K$, using the idele topology
away from $v$. Its map to $\A_K/K^\times$ sends each chart to classes
that are idelic away from that place (Proposition~\ref{prop:geo:atlas}).
Proposition~\ref{prop:ana:normal-quotient} computes the transition of
the valuation cover in these coordinates, and
Proposition~\ref{prop:ana:ball-module} sends the resulting shift to
dilation on Schwartz functions.
Connes and Consani's class-field functor sends finite abelian extensions
of $\Q$ contravariantly to finite, possibly ramified covers of an adelic
space, with arithmetic Frobenius monodromy on the unramified prime orbits
\cite[Theorems~3.9, 3.12 and~4.2]{CCKnots2025}.
Our $Q_{v,L}$ are valuation-period torsors. Under a field extension the
norm multiplies the valuation by the residue degree, and the ramification
index enters the norm formula on residue units
(Theorems~\ref{thm:geo:norm} and~\ref{thm:res:cover-norms}).
Proposition~\ref{prop:res:frobenius} realizes Frobenius on the selected
residue field by a power change of its cyclic marking.

There is also a sheaf-cohomological precedent for the Fourier--Poisson
action. On the three-point space $\mathbf P^1_{\mathbb F_1}$, Connes
and Consani define a sheaf $\Omega$ on two affine opens. On the
summands of tests with zero value at the origin and zero integral,
the two restrictions are orbital summation and its norm-weighted
inverse-coordinate form. They vanish on the polar summands.
Poisson summation identifies their images on the overlap, and
$H^1(\mathbf P^1_{\mathbb F_1},\Omega)$ realizes the zeros of Hecke
$L$-functions \cite[equations~(84)--(89), Lemma~5.3 and Theorem~5.5]{ConnesConsaniProjective}.
Grouping the finite and archimedean charts of $Y_K$ gives two opens
with intersection $C_K$. For the sheaf $\mathscr F_K$ of
Section~\ref{cmp:sheaf}, the finite winding forces the restriction
from the finite open to vanish. A constant section on $C_K$ extends
to the archimedean open precisely when its value belongs to
$\theta H_\cap$. The Mayer--Vietoris connecting map on these constant
sections embeds $H_-/\theta H_\cap$ in $H^1(Y_K,\mathscr F_K)$
(Proposition~\ref{cmp:mv}).

\subsection{Section invariants and analytic realizations}

Meyer's representation of the idele class group supplies the complete
trace and its spectral interpretation
\cite[Theorems~5.8 and~5.11]{Meyer}.
Our normal-coordinate maps realize winding operators in this
representation, and the fixed-weight quotient retains its finite zero
jets. On the section side, theta size and Euclidean duality come from
arithmetic Poisson summation
\cite[Proposition~1 and \S8]{GeerSchoof}
\cite[\S\S2.2, 5.5, 6.1]{Bost2020}.
Bost--Charles develop Hermitian section constructions and theta
invariants \cite[\S0.2.3.3--0.2.4.1]{BostCharlesTheta}.
Here these tools give explicit identities for the lattices imposed by
the integral jet condition. The Bombieri--Weyl normalization
\cite[\S2, equation~(10) and the following display]{AldazRender2024}
and orthogonal symmetrization yield the contractive product estimate
of Proposition~\ref{prop:pos-bombieri-product}.

The retained jet information also has cycle-theoretic realizations.
For the smooth projective curve and nonempty effective Cartier modulus
used here, the weight-one relative motivic complex identifies relative
units and the relative Picard group
\cite[Theorem~4.1 and Corollary~4.2]{BindaSaito}.
Section~\ref{sec:cycles} applies this to the augmentation-one
cyclotomic jet. The arithmetic Chern-character formula of
Gillet--Soul\'e supplies the finite and metric terms of the indicated
Koszul complex \cite[Theorem~1.1]{GilletSoule2014}.
The principal twist is trivial in the ordinary Picard group, but its
relative jet class and integral lattice indices retain further data.

Normed geometry provides another comparison.
Scholze's Berkovich spectra parametrize bounded multiplicative
seminorms and detect the stated arc covers
\cite[Definitions~2.13 and~3.1]{ScholzeBerkovich};
Aoki defines absolute semivalues and absolute values on ring stacks
\cite[Definition~4.31]{AokiNormed}.
Here the completions $K_v$ specify the charts, and the positive
section norms are defined on explicit coefficient modules.
Their product and division bounds, together with the comparison to
finite-cover metrics, are proved in Sections~\ref{sec:section-metrics}
and~\ref{sec:annular-sections}.

\subsection{The middle cohomology and the intersection problem}

For a smooth projective geometrically connected curve $X/\F_q$, let
$F$ be its $q$-power Frobenius endomorphism and put
$H^i=H^i_{\mathrm{\acute et}}(X_{\overline{\F}_q},\Q_\ell)$
for $\ell\nmid q$. The trace formula and the numerator of the zeta
function give, for $n\ge1$,
\[
 1+q^n-\#X(\F_{q^n})=\Tr((F^*)^n\mid H^1),\qquad
 P_1(T)=\det(1-TF^*\mid H^1).
\]
The perfect alternating cup product on $H^1$ pairs the eigenvalues
$\alpha$ and $q/\alpha$. Purity asserts that $|\alpha|=q^{1/2}$
under every complex embedding
\cite[\S1.4, (1.5.1), (1.5.4), Theorems~1.6 and~2.3, \S\S2.4--2.6]{DeligneWeilI}.
The intersection proof uses divisors on $X\times X$. Project the
diagonal and the graph of $F^n$ perpendicular to the two fibre
classes, and take minus the intersection product there. Their inner
product is $1+q^n-\#X(\F_{q^n})$ and their squared norms are $2g$
and $2gq^n$, where $g$ is the genus. Cauchy--Schwarz gives the
bound $2gq^{n/2}$ for the absolute trace
\cite[Definitions~1.1--1.2 and Lemma~1.3, p.~127]{HallouinPerret}.

On $Y_K$ the norm flow records a closed prime orbit with its period
$\log q_v$ after quotienting by $C_K^1$
(Proposition~\ref{prop:geo:norm-flow}). In the Fourier--Poisson
complex the trivial and norm characters account for the two poles of the completed Dedekind zeta function,
and the middle representation accounts for the nontrivial zeros in
Meyer's distributional trace \cite[Theorems~5.8 and~5.11]{Meyer}.
Section~\ref{sec:sheaf-comparison} places that middle representation
in ordinary sheaf cohomology. Its image restricts to zero on both
opens of the finite--archimedean cover, so the cup product of any
two image classes vanishes, also after conjugating either factor.
Theorem~\ref{thm:st:comparison} recovers the complete arithmetic trace
from a normalized relative-determinant pairing of sections.
A geometric proof of the required sign on tests with both pole moments
zero remains open. Any comparison using the section
lattices has to account for the duality quotient of
Proposition~\ref{prop:prod-jet-duality-quotient} when it identifies
opposite jet conditions. On a normal integral projective variety of dimension $n\ge1$ over
$\overline{\Q}$, Yuan--Zhang prove
$\overline M^2\cdot\overline L_1\cdots\overline L_{n-1}\le0$
for adelic $\Q$-line bundles when $\overline M$ is integrable,
the $\overline L_i$ are nef, the underlying $L_i$ are big, and
$M\cdot L_1\cdots L_{n-1}=0$ \cite[Theorem~3.2]{YuanZhangHodge}.

\paragraph{Organization.}
After the notation, Sections~\ref{sec:source}--\ref{sec:cor:correspondence}
construct the winding coefficients, arithmetic input, curve, and
correspondence maps. Sections~\ref{sec:ana:normal}--\ref{sec:ana:determinants}
give the Schwartz realization, complete trace, and normal determinants.
Sections~\ref{sec:products}--\ref{sec:cycles} develop sections, jets,
derived contacts, and cycle classes. The metric and Poisson-quotient
constructions follow in Sections~\ref{sec:section-metrics}--\ref{sec:zero-comparison}.
Section~\ref{sec:assembly} assembles their compatibility statements. Section~\ref{sec:sheaf-comparison} constructs the coefficient sheaf, proves the comparison, and computes the cup product on its image.
Section~\ref{sec:section-trace-comparison} proves the section--trace
comparison and its contact-theta and Green-intersection formulas.
The appendix proves the dyadic estimate used in the quadratic lower bound
for the logarithmic section index.
\section{Notation and conventions}
\label{sec:notation}

For a place $v$, write $K_v$ for the completion and
$j_v:K_v^\times\to C_K$ for the map induced by the idele with all
other coordinates $1$. At finite $v$, let $\mathcal O_v$, $\mathfrak m_v$
and $k_v$ be the valuation ring, maximal ideal and residue field, and put
$q_v=|k_v|$. Choose a uniformizer $\varpi_v$ and normalize
\[
 |\varpi_v|_v=q_v^{-1},\qquad |x|_v=|x|\ (K_v=\R),\qquad
 |z|_v=z\bar z\ (K_v=\C).
\]
Then $\prod_v|a|_v=1$ for $a\in K^\times$.
Set $C_K^1=\ker(|\cdot|:C_K\to\R_{>0})$.
We use the trace character
$\psi_K=\psi_\Q\circ\Tr_{\A_K/\A_\Q}$ of $\A_K/K$,
whose real component for $K=\Q$ is $x\mapsto e^{2\pi ix}$,
and self-dual additive Haar measures.
A full lattice $\Lambda\subset K$ is a free $\Z$-module of rank
$[K:\Q]$ spanning $K$ over $\Q$.

All categories of spectra and modules are stable $\infty$-categories.
For a discrete group $G$, its sphere group algebra is
$\Sph[G]=\Sigma^\infty_+G$; spectral tensor products are derived.
The notation $HA$ denotes an Eilenberg--Mac Lane spectrum,
$\Perf_A$ the retracts of finite cell $A$-modules, and $h\mathcal D$
the homotopy category of $\mathcal D$.
For a space $X$, $\Shv(X;\Sp)$ uses its ordinary open-set site and
$\Delta_X$ denotes the constant-sheaf functor.
\section{Graph cospans and finite periods}
\label{sec:source}

We use spectra with homotopy-coherent group actions. For a discrete group
\(G\), write \(\Sph[G]=\Sigma^\infty_+G\); its module category is
\(\Fun(BG,\mathrm{Sp})\). Put
\[
 R=\Sph[\Z],\qquad R_L=\Sph[C_L],\qquad C_L=\Z/L\Z,
\]
and denote the marked generator by \(t\). Tensor products of spectra and
modules are derived tensor products.

\subsection{Equivariant graphs and sewing}

\begin{definition}[The graph-cospan category]
\label{def:src:cospans}
For \(r\geq0\), let \(F_r=\Z\times\{1,\ldots,r\}\), with translation
on the first factor and \(F_0=\varnothing\). The objects of
\(\mathcal C_{\tau,g}\) are the \(F_r\). A morphism \(F_r\to F_s\) is
an equivariant cospan
\[
 F_r\longrightarrow\widetilde Y\longleftarrow F_s,
\]
where \(\widetilde Y\) is a free \(\Z\)-CW graph with finite quotient.
The legs are replaced by their mapping cylinders, and cospans are identified
under equivariant homotopy equivalence relative to their boundary maps.
Composition is equivariant homotopy pushout along the common boundary.
Disjoint union is the symmetric monoidal product. The notation \(g\)
records the first Betti numbers of connected components of the quotient
graph.
\end{definition}

Let \(B=F_1\). The fold cospan \(B\sqcup B\to B\leftarrow B\), its
reverse, and the cospans with one empty boundary define maps
\(\mu,\Delta,\eta,\epsilon\). Translation defines an automorphism
\(\tau:B\to B\).

\begin{proposition}
\label{prop:src:frobenius}
The maps \((\mu,\eta,\Delta,\epsilon)\) make \(B\) a commutative
Frobenius object, and \(\tau\) preserves these maps. In particular,
\[
 (\mu\sqcup\mathrm{id})(\mathrm{id}\sqcup\Delta)
 =\Delta\mu
 =(\mathrm{id}\sqcup\mu)(\Delta\sqcup\mathrm{id}),
\qquad
 \tau\mu=\mu(\tau\sqcup\tau).
\]
\end{proposition}
\begin{proof}
After replacing the legs by mapping cylinders, the quotient graphs of
the two associativity composites are trees with the same labelled
boundary. Contracting their edges gives the same corolla. The identical
argument proves the unit and Frobenius identities; reversing the cospans
gives the coalgebra identities. Interchanging the two boundary orbits
proves commutativity and cocommutativity. Translating the middle orbit
intertwines each structure map with simultaneous boundary translation.
\end{proof}

For \(L\geq1\), set
\[
 H_L=\mu(\mathrm{id}\sqcup\tau^L)\Delta:B\longrightarrow B.
\]
Its quotient graph has two vertices and two edges. Lifting the edges
with displacements \(0\) and \(L\) shows that its components form
\(C_L\), with both boundary maps equal to \(k\mapsto k\bmod L\).

\begin{definition}[Spectral realization of a graph cospan]
\label{def:src:realization}
For a free \(\Z\)-graph \(\widetilde Y\), let
\(P(\widetilde Y)=\Sigma^\infty_+\widetilde Y\), regarded as an
\(R\)-module through translation. Applying \(P\) to the three entries
of a graph cospan gives a cospan of perfect \(R\)-modules.
\end{definition}

\begin{proposition}[The sewing triangle]
\label{prop:src:sewing}
Spectral realization defines a symmetric monoidal functor
\[
 \mathcal C_{\tau,g}\longrightarrow
 h\operatorname{Cospan}(\Perf_R).
\]
If \(\widetilde W=\widetilde Y\cup^h_{F_s}\widetilde Z\), there is a
natural cofiber sequence
\begin{equation}
 R^{\oplus s}\xrightarrow{(j_*,-i_*)}
 P(\widetilde Y)\oplus P(\widetilde Z)
 \longrightarrow P(\widetilde W).
 \label{eq:src:sewing}
\end{equation}
Moreover,
\begin{equation}
 P(\widetilde Y_{H_L})\simeq
 \cofib(t^L-1:R\longrightarrow R)\simeq R_L
 \quad\text{in }\Mod_R.
 \label{eq:src:period-cofiber}
\end{equation}
\end{proposition}
\begin{proof}
Choose one lift of every cell of the finite quotient graph. The cellular
construction of its suspension spectrum uses finitely many free
\(R\)-cells, so it is perfect. Suspension spectra preserve homotopy
pushouts. In a stable category the pushout of \(P(\widetilde Y)\) and
\(P(\widetilde Z)\) over \(R^{\oplus s}\) is the cofiber of the
difference of the two boundary maps, which proves
\eqref{eq:src:sewing} and compatibility with composition. Disjoint unions
give direct sums. In the graph of \(H_L\), contract the translates of the
edge with displacement zero. The remaining attaching map is \(t^L-1\).
Each component of the covering graph is contractible; its component
\(\Z\)-set is \(C_L\). This proves the second equivalence in
\eqref{eq:src:period-cofiber}.
\end{proof}

\subsection{The indexed finite-period system}

\begin{definition}[Finite periods and their transfers]
\label{def:src:periods}
Let \(O_L\) be the fibre over a marked point of an oriented circle of
the degree-\(L\) covering of that circle. It is a \(C_L\)-torsor, with
generator induced by positive winding. A marked sheet identifies it with
\(C_L\). For \(L\mid hL\), let
\(q_{hL,L}:O_{hL}\to O_L\) be the quotient map. It induces a map
\(q_*\) of permutation spectra. Define its additive stable transfer
\(q^!=D(q_*)\) by Spanier--Whitehead duality
\(D(X)=\operatorname{Map}(X,\Sph)\) and the canonical self-dualities
of finite permutation spectra in \(\Fun(BC_{hL},\mathrm{Sp})\), where
\(\Sph[O_L]\) carries the action through \(C_{hL}\to C_L\).
On the underlying spectra these maps are
\[
 q_*e_x=e_{q(x)},\qquad
 q^!e_y=\sum_{q(x)=y}e_x.
\]
After forgetting the group action, \(q_*q^!=h\,\mathrm{id}\).
The same identity holds equivariantly after applying
\(HA\otimes_{\Sph}-\) for an ordinary commutative ring \(A\).
Transfers of composites are composed in the reverse order.
\end{definition}

Finite-set duality makes the construction independent of an ordering
of the fibres and supplies its coherent composition and pullback
identities. If \(d\mid m,n\), then
\(O_m\times_{O_d}O_n\), with simultaneous translation, is a disjoint
union of \(\gcd(m,n)/d\) orbits of period \(\operatorname{lcm}(m,n)\):
its cardinality is \(mn/d\), and every translation orbit has the stated
length.

\begin{definition}[Primitive spectra]
\label{def:src:primitive}
Set \(U_1=\Sph\) with trivial action. If \(L=\ell^a\), where \(\ell\)
is prime and \(a\geq1\), put
\[
 U_{\ell^a}=
 \cofib\bigl(q^!: \Sph[C_{\ell^{a-1}}]
                 \longrightarrow\Sph[C_{\ell^a}]\bigr)
 \quad\text{in }\Fun(BC_{\ell^a},\mathrm{Sp}),
\]
where the source action is inflated along the quotient.
For general \(L\), take the external smash product over
\(\ell^a\Vert L\), and use the Chinese remainder isomorphism
\(C_L\simeq\prod_{\ell^a\Vert L}C_{\ell^a}\) that sends the marked
generator to the marked generator in every factor.
For a \(C_L\)-torsor \(P\), the same construction on its prime-power
quotients defines \(U(P)\).
\end{definition}

These primitive spectra are the Spanier--Whitehead duals of the
cyclotomic Moore spectra in Rezk's construction, as presented by Zhang
\cite[Construction~3.3.2 and the external smash product on pp.~25--26]{ZhangDirichlet}.
Indeed, that construction takes the fibre of the quotient map \(q_*\);
duality and the canonical self-dualities of finite permutation spectra
give the transfer cofiber above, with its cyclic action.

\begin{definition}[The level \(\mathbb F_{1^L}\)]
\label{def:src:f1level}
The symbol \(\mathbb F_{1^L}\) denotes the level-\(L\) entry
\[
 \bigl(O_L,\ R_L,\ \Mod_{R_L},\ U_L\bigr)
\]
in the system indexed by divisibility of positive integers. The transition
data are the quotient maps \(q_{hL,L}\), their additive transfers, the
induced ring-spectrum maps \(R_{hL}\to R_L\), and their extension and
restriction of scalars. The integer \(L\) is the covering degree and
the period of the component orbit of \(H_L\).
\end{definition}

\begin{lemma}[Refinement on cellular chains]
\label{lem:src:refinement}
Let \(A\) be a commutative ring, \(T=A[t,t^{-1}]\), and
\(D_L=[T\xrightarrow{t^L-1}T]\), in cohomological degrees \(-1,0\).
For \(N_{h,L}=\sum_{j=0}^{h-1}t^{jL}\), the quotient and transfer are
represented by
\[
 D_{hL}\xrightarrow{(N_{h,L},1)}D_L,
 \qquad
 D_L\xrightarrow{(1,N_{h,L})}D_{hL}.
\]
Their composite on \(D_L\) is chain homotopic to \(h\,\mathrm{id}\).
\end{lemma}
\begin{proof}
Both chain-map identities follow from
\((t^L-1)N_{h,L}=t^{hL}-1\). On zeroth homology the first map reduces
the exponent modulo \(L\), and the second sums its \(h\) lifts.
The composite is multiplication by \(N_{h,L}\) in both degrees.
The polynomial \(N_{h,L}-h\) is divisible by \(t^L-1\); multiplication
by the quotient, from degree zero to degree \(-1\), is the required
chain homotopy.
\end{proof}

\subsection{Primitive intersections and their endomorphisms}

\begin{proposition}
\label{prop:src:primitive-lattice}
The underlying spectrum of \(U_m\) is a wedge of \(\varphi(m)\)
spheres. Its integral homology is concentrated in degree zero and,
with the marked action,
\[
 H_0(U_m;\Z)\simeq\Z[t]/(\Phi_m(t)).
\]
There is an \(R\)-module equivalence
\begin{equation}
 U_m\simeq\cofib(\Phi_m(t):R\longrightarrow R).
 \label{eq:src:primitive-resolution}
\end{equation}
\end{proposition}
\begin{proof}
For \(m=\ell^a\), the transfer includes one sum of basis vectors for
each coset of \(C_\ell\). Within each coset that sum extends to an
integral basis. The quotient is free and its rotation module is
\[
 \Z[t]/(t^{\ell^a}-1,\,1+t^{\ell^{a-1}}+\cdots+
 t^{(\ell-1)\ell^{a-1}})=\Z[t]/(\Phi_{\ell^a}).
\]
Integral changes of basis lift to equivalences of finite wedges of
spheres. Thus the cofiber is a wedge of spheres of the indicated rank.
For relatively prime periods, the external tensor product of the
cyclotomic lattices is generated by the product of their distinguished
roots: each root can be recovered by a power of the product. Equality
of ranks identifies this tensor product with \(\Z[t]/(\Phi_m)\).

The class \(1\in\pi_0U_m\) defines an \(R\)-linear map \(R\to U_m\).
Its composite with \(\Phi_m(t)\) represents zero in \(\pi_0U_m\), so
choose a nullhomotopy and extend the map across its cofiber. The resulting
map induces an integral homology isomorphism. Its cofiber is connective
with zero integral homology, and is therefore contractible by the stable
Hurewicz theorem.
\end{proof}

\begin{definition}[Primitive contact]
\label{def:src:contact}
For positive integers \(m,n\), define
\[
 T_{m,n}=U_m\otimes_R U_n.
\]
The marked action on each factor specifies its \(R\)-module structure.
We give \(T_{m,n}\) the additional \(C_m\)-action induced by the first
factor, and use \(\End_{C_m}(T_{m,n})\) for the endomorphism ring spectrum
in \(\Fun(BC_m,\mathrm{Sp})\).
\end{definition}

\begin{proposition}[The underlying contact spectrum]
\label{prop:src:contact-spectrum}
Let \(m<n\). If \(n/m\) is not a positive power of a prime, then
\(T_{m,n}\simeq0\). If \(n=mp^a\), where \(p\) is prime and
\(a\geq1\), its underlying spectrum is a wedge of \(\varphi(m)\)
copies of \(\Sph/p\).
\end{proposition}
\begin{proof}
By \eqref{eq:src:primitive-resolution}, the contact is the cofiber of
\(\Phi_n(t)\) on \(U_m\). On its free integral homology this map is
multiplication by \(\Phi_n(\zeta_m)\) on \(\Z[\zeta_m]\). For \(m>1\),
the absolute determinant is the cyclotomic resultant
\cite[Theorem~9]{BHM}; it is \(1\) unless \(n/m=p^a\), and is
\(p^{\varphi(m)}\) in that case. For \(m=1\), the corresponding
formula is \(\Phi_n(1)=p\) for a prime power and \(1\) otherwise
\cite[Lemma~4]{BHM}.

In the prime-power case, write \(m=up^b\), with \(p\nmid u\).
Reduction modulo \(p\) gives
\(\overline\Phi_{up^r}=\overline\Phi_u^{\varphi(p^r)}\) for
\(r\geq1\), and hence \(\overline\Phi_m\mid\overline\Phi_n\).
It follows that \(\Phi_n(\zeta_m)=p\alpha\) with
\(\alpha\in\Z[\zeta_m]\). The resultant gives
\(\lvert N_{\Q(\zeta_m)/\Q}(\alpha)\rvert=1\), so \(\alpha\) is a
unit. In bases of the underlying wedges of spheres, the map is therefore
equivalent to \(p\) times an invertible integral matrix. Taking its
cofiber proves the assertion.
\end{proof}

\begin{lemma}[The scalar Moore endomorphism ring]
\label{lem:src:moore}
For \(M_p=\Sph/p\), composition gives
\[
 \pi_0\End(M_p)=
 \begin{cases}
  \F_p,&p\ne2,\\
  \Z/4,&p=2.
 \end{cases}
\]
\end{lemma}
\begin{proof}
If \(p\) is odd, the Moore cofiber sequence gives
\(\pi_1M_p=0\) and \(\pi_0M_p=\F_p\). Applying \([-,M_p]\) to it
identifies \([M_p,M_p]\) with \(\F_p\), generated by the identity.
For \(p=2\), it gives an extension of \(\Z/2\) by \(\Z/2\). With
\(\rho:\Sph\to M_2\), \(\beta:M_2\to\Sigma\Sph\), and the Hopf
map \(\eta:\Sigma\Sph\to\Sph\), the identity satisfies
\(2\,\mathrm{id}_{M_2}=\rho\eta\beta\ne0\)
\cite[Proposition~4.1]{StricklandMoore}. Thus the identity has additive
order four and generates the endomorphism ring.
\end{proof}

\begin{definition}[The prime-to-\(p\) primitive factor]
\label{def:src:primitive-factor}
For \(p\nmid m\), let \(e_m\) be the idempotent in
\(\Z_{(p)}[C_m]\) given by
\[
 e_m=\prod_{\ell\mid m}
 \left(1-\frac1\ell\sum_{h\in C_\ell}h\right),
\]
where \(C_\ell\) is the subgroup of order \(\ell\), and set
\(e_1=1\). Let \(A_m\) be the corresponding idempotent factor of the
commutative ring spectrum \(\Sph_{(p)}[C_m]\).
\end{definition}

Each factor in the product is idempotent, and these factors commute.
The algebra \(A_m\) is finite free over \(\Sph_{(p)}\), and
\(\pi_0A_m=\Z_{(p)}[t]/(\Phi_m)\). Indeed, after adjoining the
\(m\)-th roots of unity, \(e_m\) selects those characters that are
nontrivial on every subgroup \(C_\ell\), which are the primitive
characters. The character decomposition is valid over a finite
faithfully flat extension of \(\Z_{(p)}\), because \(m\) is a unit.
It descends to the displayed factor. The module \(\pi_0A_m\) is a
retract of a finite free \(\Z_{(p)}\)-module, hence is free of rank
\(\varphi(m)\). Lift a basis to a map
\(\Sph_{(p)}^{\vee\varphi(m)}\to A_m\). Since \(A_m\) is an
idempotent retract of a finite wedge of localized spheres, this map
induces an integral homology isomorphism. Its connective cofiber is
contractible by the stable Hurewicz theorem, proving the claimed
freeness over \(\Sph_{(p)}\).

\begin{theorem}[Indexed endomorphism coefficients]
\label{thm:src:endomorphisms}
Let \(p\nmid m\), \(a\geq1\), and \(E_p=\End(\Sph/p)\). There is
an equivalence of associative ring spectra
\begin{equation}
 \End_{C_m}(T_{m,mp^a})\simeq A_m\otimes_{\Sph_{(p)}}E_p.
 \label{eq:src:end-spectrum}
\end{equation}
Under this equivalence the marked rotation is the distinguished
cyclotomic element. In particular,
\begin{equation}
 E_m:=\pi_0\End_{C_m}(T_{m,mp})
 \simeq R_p[t]/(\Phi_m(t)),\qquad
 R_p=
 \begin{cases}\F_p,&p\ne2,\\ \Z/4,&p=2.\end{cases}
 \label{eq:src:end-ring}
\end{equation}
\end{theorem}
\begin{proof}
The primitive projector acts as the identity on \(U_m\) after
\(p\)-localization. The class \(1\) gives an equivariant map
\(A_m\to U_{m,(p)}\); its integral homology map is an isomorphism, so
it is an equivalence. The proof of
Proposition~\ref{prop:src:contact-spectrum} identifies
\(\Phi_{mp^a}(t)\in\pi_0A_m\) with \(p\) times a unit. Thus
\[
 T_{m,mp^a}\simeq A_m\otimes_{\Sph_{(p)}}M_p
\]
as a \(\Sph_{(p)}[C_m]\)-module supported on the factor \(A_m\).
Such modules have the same mapping spectra as \(A_m\)-modules.
The Moore spectrum is dualizable, and hence extension of scalars gives
\[
 \End_{A_m}(A_m\otimes M_p)
 \simeq A_m\otimes\End(M_p),
\]
compatibly with composition. The contact is \(p\)-local by
Proposition~\ref{prop:src:contact-spectrum}, so this also computes its
original \(C_m\)-equivariant mapping spectrum. Since \(A_m\) is finite
free over the localized sphere, taking \(\pi_0\) and applying
Lemma~\ref{lem:src:moore} gives \eqref{eq:src:end-ring}.
\end{proof}

\begin{proposition}[Linear maps and change of coefficients]
\label{prop:src:coefficients}
Fix \(p\nmid m\) and set \(T_m=T_{m,mp}\). The category whose objects
are retracts of finite sums of \(T_m\) and whose morphisms are
degree-zero \(C_m\)-equivariant maps is equivalent, through
\([T_m,-]_{C_m}\), to the category of finite projective \(E_m\)-modules.
If \(m\mid M\) and \(p\nmid M\), the subgroup inclusion
\(C_m\hookrightarrow C_M\), \(1\mapsto M/m\), induces
\begin{equation}
 E_m\longrightarrow E_M,\qquad t\longmapsto t^{M/m}.
 \label{eq:src:coefficient-map}
\end{equation}
Restriction corresponds to restriction of scalars. Induction followed by
the primitive projector for \(M\) corresponds to extension of scalars.
\end{proposition}
\begin{proof}
On finite sums of \(T_m\), the morphisms are matrices over the ring
\(E_m\). Splitting idempotent matrices on both sides gives the first
equivalence. A primitive \(M\)-th root raised to \(M/m\) is a primitive
\(m\)-th root, which proves \eqref{eq:src:coefficient-map}.
The induced map between primitive group-algebra factors is finite free
of rank \(\varphi(M)/\varphi(m)\): after a faithfully flat extension
splitting the prime-to-\(p\) roots, each primitive \(m\)-character has
that many primitive \(M\)-character lifts. Finite projectivity descends
from this splitting, and the base ring \(\pi_0A_m\) is semilocal,
so the resulting constant-rank projective module is free. Restriction of
\(A_M\otimes M_p\) therefore belongs to the additive idempotent
completion generated by \(A_m\otimes M_p\). The group-algebra
restriction and induction adjunction, followed by the primitive
projectors, gives the two scalar functors on matrices and their
idempotent summands.
\end{proof}

For \(m=p^d-1\), \(M=p^e-1\), and \(d\mid e\), compatible residue
factors in \eqref{eq:src:coefficient-map} give
\(\F_{p^d}\hookrightarrow\F_{p^e}\). The same functors give extension
and restriction of finite-dimensional vector spaces. At \(p=2\), first
reduce the corresponding matrices over the length-two Witt factors
modulo \(2\). On the resulting field the map \(t\mapsto t^p\) is
Frobenius. The trace and determinant of multiplication by an element,
viewed on the restricted vector space, are its field trace and norm.

\begin{proposition}[Completed primitive coefficients]
\label{prop:src:completed}
Let \(p\nmid m\), and let \(U_{m,p}\) be the \(p\)-completion of
\(U_m\). In \(p\)-complete spectra with \(C_m\)-action,
\[
 \pi_0\End_{C_m}(U_{m,p})
 =\Z_p[t]/(\Phi_m(t))
 \simeq\prod_{f\mid\overline\Phi_m}
 W\bigl(\F_p[t]/(f)\bigr),
\]
where \(f\) ranges over the monic irreducible factors and \(W(k)\)
is the unramified complete discrete valuation ring with residue field
\(k\) and uniformizer \(p\).
\end{proposition}
\begin{proof}
The primitive-factor argument above identifies \(U_{m,p}\) with the
primitive factor of \(\Sph_p^\wedge[C_m]\). Its endomorphism ring is
that factor, and its zeroth homotopy ring is the displayed polynomial
quotient. The reduction of \(\Phi_m\) is separable because \(p\nmid m\).
Hensel lifting of its coprime factors decomposes this finite etale
\(\Z_p\)-algebra into the indicated unramified factors.
\end{proof}
\section{Marked graph returns and arithmetic lattices}
\label{sec:arithmetic}

\begin{definition}[Marked integral return]
\label{arith:def:return}
A marked integral return is a triple \((Y,f,e)\), where \(Y\) is a finite connected
based graph, \(f:Y\to Y\) is a based cellular map, and \(e\in H_1(Y,\Z)\).
Write \(M=H_1(Y,\Z)\), \(d=\operatorname{rank}M>0\), and \(U=f_*\).
The return is cyclic and irreducible if
\(e,Ue,\ldots,U^{d-1}e\) is a \(\Z\)-basis of \(M\) and the characteristic
polynomial of \(U\) is irreducible over \(\Q\).
\end{definition}

For such a return, the map \(P(U)\mapsto P(U)e\) identifies
\(A_0=\Z[U]\) with \(M\). The rational algebra \(K=\Q[U]\) is a number field,
and \(M\otimes\Q\) is its regular module.

\begin{proposition}[Maximal order from a bounded family of covers]
\label{arith:prop:order}
Set
\[
 G_{ij}=\Tr_{K/\Q}(U^{i+j})\quad(0\leq i,j<d),
 \qquad D=|\det G|.
\]
For the \(U\)-stable lattices \(D M\subseteq L\subseteq M\), form the multiplier rings
\[
 A_L=\{a\in K:aL\subseteq L\}.
\]
This finite family has a largest member, equal to \(\mathcal O_K\).
Each lattice \(L\) is the image of the closed paths in a finite abelian cover of
\(Y\), and \(f\) lifts to that based cover.
\end{proposition}

\begin{proof}
The separability of \(K/\Q\) implies \(D\ne0\). An element of \(A_L\) acts on the
free integral lattice \(L\), so its characteristic polynomial is monic with integer
coefficients. Its action on \(K\) is multiplication by that element, and hence
\(A_L\subseteq\mathcal O_K\).

Let
\[
 A_0^\vee=\{a\in K:\Tr_{K/\Q}(aA_0)\subseteq\Z\}.
\]
The integrality of traces gives \(\mathcal O_K\subseteq A_0^\vee\). In the basis
\(1,U,\ldots,U^{d-1}\), the adjugate formula for \(G^{-1}\) gives
\(D A_0^\vee\subseteq A_0\). Consequently
\[
 D A_0\subseteq D\mathcal O_K\subseteq A_0.
\]
The lattice \(D\mathcal O_K\) is \(U\)-stable, and its multiplier ring is
\(\mathcal O_K\). It therefore realizes the largest member of the family.
There are finitely many intermediate lattices because \(M/D M\) is finite.

Let \(H_L\) be the inverse image of \(L\) under the abelianization
\(\pi_1(Y)\to M\). The associated based cover has finite deck group \(M/L\),
and the image of its fundamental group in \(M\) is \(L\).
The inclusion \(UL\subseteq L\) gives \(f_*(H_L)\subseteq H_L\), which is the based
lifting criterion.
\end{proof}

\begin{proposition}[Congruence completions]
\label{arith:prop:completion}
Put \(O=\mathcal O_K\). The lattices \(kDO\subseteq DO\), for positive integers
\(k\), and the marked multiplication action give the rings \(O/kO\) and their
transition maps. Their inverse limit, indexed by divisibility, is
\[
 \widehat O=\varprojlim_k O/kO
     \cong\prod_{v\nmid\infty}\mathcal O_v.
\]
Moreover
\[
 \A_{K,f}=\bigcup_{n\geq1}n^{-1}\widehat O,
 \qquad K_\infty=K\otimes_\Q\R,
 \qquad \A_K=\A_{K,f}\times K_\infty.
\]
The diagonal maps on the congruence quotients give the diagonal embedding
\(K\hookrightarrow\A_K\).
\end{proposition}

\begin{proof}
Multiplication by \(D\) identifies \(O/kO\) with \(DO/kDO\) as an \(O\)-module,
and transports its multiplication to this quotient. For \(k\mid\ell\), the inclusion
\(\ell DO\subseteq kDO\) gives reduction modulo \(k\).
The Chinese remainder theorem decomposes the inverse limit into its completions
at the nonzero prime ideals of \(O\).

An element of the finite adele ring is integral at all but finitely many places.
A positive integer can be chosen with sufficiently large valuation at the rational
primes below those places to make every component integral. This proves the
displayed union; the reverse inclusion follows from the definition of the
restricted product. The finite real algebra \(K\otimes_\Q\R\) is the product of
the real and complex completions. The compatibility of all diagonal maps follows
from their definition by evaluation of the same element of \(K\).
\end{proof}

\begin{definition}[Polynomial morphism of marked returns]
\label{arith:def:morphism}
Let \((M,U,e)\) and \((M',U',e')\) be the homology data of cyclic irreducible
returns. A polynomial morphism is a pair \((J,r)\), where
\(J:M\to M'\) is an injective integral homomorphism, \(r\in\Z[x]\),
\[
 J(e)=e',\qquad JU=r(U')J.
\]
Two pairs with the same \(J\) and the same value \(r(U')\) define the same morphism.
Composition is
\((J_2,r_2)\circ(J_1,r_1)=(J_2J_1,r_1\circ r_2)\).
\end{definition}

\begin{proposition}[Induced field maps]
\label{arith:prop:fieldmap}
A polynomial morphism induces the field embedding
\[
 \phi:K\hookrightarrow K',\qquad P(U)\longmapsto P(r(U')).
\]
It maps \(\mathcal O_K\) into \(\mathcal O_{K'}\), and induces compatible maps on
the congruence completions and adele rings.
\end{proposition}

\begin{proof}
Let \(p_U\) be the characteristic polynomial of \(U\). The intertwining identity
gives
\[
 p_U(r(U'))e'=Jp_U(U)e=0.
\]
Multiplication by an element of \(K'\) is
determined by its value on \(1=e'\), so \(p_U(r(U'))=0\) in \(K'\).
Evaluation therefore defines a unital map \(\Q[x]/(p_U)\to K'\), which is
injective because its source is a field. The same calculation gives
\(J(P(U)e)=\phi(P(U))e'\). Integral elements remain integral under a field
embedding. Reduction modulo positive integers and passage to inverse limits give
the finite completion map. Extending at the archimedean places and inverting
positive integers gives the adele map. Polynomial substitution proves
compatibility with composition.
\end{proof}

\begin{definition}[Lattice test and trace-dual lattice]
\label{arith:def:lattice}
For a full lattice \(\Lambda\subset K\), define
\[
 \widehat\Lambda=\Lambda\otimes_\Z\widehat\Z\subset\A_{K,f},
 \qquad
 \Lambda^\vee=\{a\in K:\Tr_{K/\Q}(a\Lambda)\subseteq\Z\}.
\]
The finite lattice test is the characteristic function
\(1_{\widehat\Lambda}\in\cS(\A_{K,f})\). Its all-place tests are
\(1_{\widehat\Lambda}\otimes\Phi_\infty\), with
\(\Phi_\infty\in\cS(K_\infty)\).
\end{definition}

\begin{proposition}[Fourier transform of a lattice test]
\label{arith:prop:latticefourier}
For the trace character and self-dual measure of Section~\ref{sec:notation},
\[
 \mathcal F_f(1_{\widehat\Lambda})
       =\operatorname{vol}(\widehat\Lambda)\,
             1_{\widehat{\Lambda^\vee}}.
\]
If \(\Lambda'\subseteq\Lambda\) has finite index \(a\), then
\(\operatorname{vol}(\widehat{\Lambda'})=a^{-1}
\operatorname{vol}(\widehat\Lambda)\).
\end{proposition}

\begin{proof}
Choose a \(\Z\)-basis of \(\Lambda\) and the dual basis for the nondegenerate
trace pairing. At every rational prime \(p\), the annihilator of
\(\Lambda\otimes\Z_p\) under the character induced by the trace is
\(\Lambda^\vee\otimes\Z_p\), because the standard character of \(\Q_p\) has
annihilator \(\Z_p\). Taking the product over \(p\) identifies the annihilator of
\(\widehat\Lambda\) with \(\widehat{\Lambda^\vee}\).
The integral of a character over a compact group is its volume for the trivial
character and zero otherwise. This proves the Fourier formula.
Finally, \(\widehat\Lambda/\widehat{\Lambda'}\) has \(a\) elements, and the
cosets have equal Haar measure.
\end{proof}
\section{The space and its winding covers}
\label{sec:geometry}

Write \(C=C_K\). For each place \(v\), the homomorphism
\(j_v:K_v^\times\to C\) sends a local element to the class of the idele
whose other components are \(1\). We use right actions in forming
associated quotients.

\begin{definition}[Normal charts and the associated space]
\label{def:geo:atlas}
The normal chart at \(v\) is the quotient space
\[
 N_v=(C\times K_v)/K_v^\times,
 \qquad (c,z)\cdot h=(cj_v(h),h^{-1}z).
\]
Its punctured part is identified with \(C\) by
\[
 \nu_v:N_v^\times\longrightarrow C,
 \qquad [c,z]\longmapsto cj_v(z).
\]
Let \(Y_K\) be the topological space obtained by gluing all \(N_v\)
along these identified open subsets. The closed orbit of the chart is
\[
 Z_v=\{[c,0]:c\in C\}\simeq C/j_v(K_v^\times).
\]
Thus the charts cover \(Y_K\), and \(N_v\cap N_w=C\) for \(v\ne w\).
\end{definition}

\begin{proposition}
\label{prop:geo:atlas}
Each chart \(N_v\) is locally compact and Hausdorff. The subset \(C\)
is open in \(Y_K\), and multiplication in the first coordinate defines
a continuous action of \(C\) on \(Y_K\). Each \(Z_v\) is a compact
\(C\)-orbit. There is a continuous \(C\)-equivariant map
\[
 \nu:Y_K\longrightarrow \A_K/K^\times.
\]
\end{proposition}

\begin{proof}
The map \(j_v\) is injective: if a single-place idele is diagonal, any
other component forces the diagonal element to be \(1\). It is proper.
Indeed, for a compact subset \(D\subset C\), the norm of an element of
\(j_v^{-1}(D)\) lies in a fixed compact interval in \(\R_{>0}\).
The corresponding local norm annulus in \(K_v^\times\) is compact,
and \(j_v^{-1}(D)\) is closed in it.

The action of \(K_v^\times\) on \(C\times K_v\) is free and proper.
For properness, the two \(C\)-coordinates of a compact subset of
\((C\times K_v)^2\) bound \(j_v(h)=c^{-1}c'\) in a compact subset
of \(C\); properness of \(j_v\) then bounds \(h\) in a compact subset
of \(K_v^\times\). The action map consequently has compact inverse
images of compact sets. Its orbit relation is closed and its quotient
map is open, so the quotient is Hausdorff. The image of a relatively
compact open neighbourhood upstairs has a compact neighbourhood
downstairs, proving local compactness.

The punctured identification has continuous inverse \(c\mapsto[c,1]\).
The zero-coordinate subset is closed and saturated. Left multiplication
by \(C\) commutes with the right action and agrees under all the
punctured identifications; it therefore descends to the glued space.

Let \(C^1=\ker(|\cdot|:C\to\R_{>0})\). This group is compact
\cite[\S5.8, p.~33]{Poonen}. At a finite place, every class modulo
\(j_v(K_v^\times)\) has a representative whose norm belongs to
\([1,q_v]\); this is a compact subset of \(C\), since the norm map
has a continuous section and compact kernel. At an infinite place,
every class has a representative in \(C^1\). Thus \(Z_v\) is compact.

For \(c\in C\), choose an idele \(a\) representing it and set
\[
 b_w=a_w\quad(w\ne v),\qquad b_v=a_vz.
\]
The class of \(b\) modulo \(K^\times\) is independent of the choice
of \(a\), and the right action defining \(N_v\) leaves it unchanged.
The continuous map on \(\A_K^\times\times K_v\) therefore descends
to \(N_v\). On the punctured part it is the usual map from idele
classes to adele classes. These maps agree on overlaps and define
\(\nu\).
\end{proof}

At a finite place put \(U_v=\mathcal O_v^\times\),
\(\ell_v=\log q_v\), and
\(H_{v,L}=\ord_v^{-1}(L\Z)\) for \(L\ge1\).

\begin{definition}[Valuation cover]
\label{def:geo:valuation-cover}
Over the chart \(N_v\) consider
\[
 P_v=(C\times K_v)/U_v.
\]
The right action of \(K_v^\times/U_v\simeq\Z\) is
\[
 [c,z]\cdot n=[cj_v(\varpi_v^n),\varpi_v^{-n}z].
\]
It is independent of the choice of uniformizer. On its punctured
part use the trivialization
\begin{equation}
 \label{eq:geo:cover-frame}
 [c,z]\longmapsto\bigl(cj_v(z),-\ord_v(z)\bigr)
       \in C\times\Z.
\end{equation}
The valuation cover \(Q_v\) on \(Y_K\) is obtained by gluing \(P_v\)
to \(N_w\times\Z\), for every \(w\ne v\), using
\eqref{eq:geo:cover-frame} on the common open subset.
Its period-\(L\) quotient is denoted by \(Q_{v,L}=Q_v/L\Z\).
\end{definition}

\begin{theorem}
\label{thm:geo:valuation-cover}
The map \(Q_v\to Y_K\) is a principal covering with discrete group
\(\Z\). For every \(L\ge1\), \(Q_{v,L}\to Y_K\) is a principal
\(C_L\)-cover, and
\[
 Q_{v,L}|_{N_v}\simeq(C\times K_v)/H_{v,L}.
\]
If \(L\mid M\), the quotient maps \(Q_{v,M}\to Q_{v,L}\)
compose in towers. Their suspension-spectrum transfer maps satisfy
the same composition law.
\end{theorem}

\begin{proof}
The function \(a([c,z])=\log|c|\) on \(P_v\) is well defined,
because \(|u|_v=1\) for \(u\in U_v\), and
\(a(x\cdot n)=a(x)-n\ell_v\). The inverse image of an open interval
of length less than \(\ell_v\) has disjoint nonidentity translates.
The orbit map is open, so its restriction to such an inverse image
is a homeomorphism onto an open subset of \(N_v\). The translates
give a covering trivialization. Its structure group is
\(K_v^\times/U_v\), identified with \(\Z\) by the valuation.

Formula \eqref{eq:geo:cover-frame} is invariant under \(U_v\) and
equivariant for translation of its second coordinate. An inverse is
\[
 (a,n)\longmapsto
 [aj_v(\varpi_v^n),\varpi_v^{-n}].
\]
The chart covers consequently glue as principal covers. Dividing by
\(L\Z\) gives the claimed local quotient. The quotient maps compose
by their definitions. Apply the canonical duality of finite
permutation spectra to each finite quotient of the cover. This
construction descends along the covering trivializations, including
their transition permutations and all composition homotopies.
Duality reverses composition, so the resulting transfers satisfy
the tower identity. The same finite-set construction gives
compatibility with pullback.
\end{proof}

\begin{definition}[Realization along a valuation cover]
\label{def:geo:realization}
For a space with a \(\Z\)-action, or a spectrum with a \(\Z\)-action,
its realization along \(Q_v\) is the associated locally constant
sheaf obtained by descending its constant sheaf on a trivializing
cover of \(Q_v\). The associated-bundle convention is
\[
 (p+n,x)\sim(p,n\cdot x).
\]
For a topological \(\Z\)-graph \(G\), we apply this construction to
the homotopy type \(\operatorname{Sing}(G)\), and then apply the
fibrewise suspension-spectrum functor.
\end{definition}

Associated descent and fibrewise suspension preserve homotopy colimits.
Consequently a homotopy pushout of \(\Z\)-graphs gives a homotopy
pushout of sheaves of spectra, and graph sewing gives its cofiber
triangle. For the permutation spectrum \(\Sph[C_L]\), the realization
is the fibrewise suspension spectrum of the finite cover \(Q_{v,L}\).

\begin{proposition}[Norm-flow charts]
\label{prop:geo:norm-flow}
The quotient \(C^1\backslash Y_K\) has generic open subset
\(\R_{>0}\). At a finite place its chart is
\[
 \bigl(\R_{>0}\times(\Z\cup\{+\infty\})\bigr)/\Z,
 \qquad (r,n)\cdot k=(rq_v^{-k},n-k).
\]
The punctured coordinate is \(rq_v^{-n}\), and its boundary is
\(\R_{>0}/q_v^\Z\). The boundary over \(Z_v\) in \(C^1\backslash Q_{v,L}\)
is \(\R_{>0}/q_v^{L\Z}\), of logarithmic length \(L\ell_v\).
At a real or complex place the quotient chart is \([0,\infty)\),
with coordinate \(r|z|_v\).
\end{proposition}

\begin{proof}
The norm induces \(C/C^1\simeq\R_{>0}\). After taking this quotient
in the first coordinate, quotient by the kernel of the local norm
in the second. At a finite place the result is
\(\Z\cup\{+\infty\}\), with the valuation topology and
\(\ord_v(0)=+\infty\). The remaining action is the displayed
integer action. On the boundary only the \(r\)-coordinate remains,
modulo multiplication by \(q_v^\Z\), or by \(q_v^{L\Z}\) at
period \(L\). At an infinite place the local norm gives the quotient
\([0,\infty)\); the remaining action of \(\R_{>0}\) identifies
\((r,s)\) with \(rs\).
\end{proof}

For a principal \(\Z\)-torsor \(P\), write \(fP\) for its pushout
along multiplication by an integer \(f\). A sum of such torsors
means their contracted product along addition of the structure groups.

\begin{theorem}[Norm maps]
\label{thm:geo:norm}
A finite extension \(L/K\) induces a continuous map
\(N_{L/K}:Y_L\to Y_K\), equivariant for the idele-class norm. On
the chart associated to \(w\mid v\) it is
\[
 [c,z]\longmapsto[N_{L/K}c,N_{L_w/K_v}z].
\]
These maps compose in towers of fields. For each finite place \(v\)
of \(K\), there is an isomorphism of principal \(\Z\)-covers
\begin{equation}
 \label{eq:geo:norm-cover}
 N_{L/K}^{*}Q_{K,v}\simeq
       \sum_{w\mid v}f(w/v)Q_{L,w}.
\end{equation}
The isomorphisms commute with reduction of the deck group modulo
\(r\), for every \(r\ge1\).
\end{theorem}

\begin{proof}
The identity \(N_{L/K}(j_w(h))=j_v(N_{L_w/K_v}h)\) proves that the
chart map respects the balancing relation. On each punctured chart
it is the idele-class norm; hence the chart maps glue. Local norms
and idele norms compose in towers.

Fix \(w\mid v\), and write \(f=f(w/v)\). The map
\[
 (C_L\times L_w)/\mathcal O_w^\times
 \longrightarrow (C_K\times K_v)/\mathcal O_v^\times,
 \qquad [c,z]\longmapsto[Nc,Nz]
\]
is well defined and is equivariant for \(\Z\xrightarrow{f}\Z\),
since \(\ord_v(Nz)=f\ord_w(z)\). It therefore identifies the
pushout of \(Q_{L,w}|_{N_w}\) along multiplication by \(f\) with
the pulled-back cover on that chart. In the right side of
\eqref{eq:geo:norm-cover}, all the other torsors are trivial on
\(N_w\). On a chart above a place different from \(v\), both sides
are trivial. On the common open subset the transition labels are
\(-\ord_v(Nz)=-f\ord_w(z)\), by
\eqref{eq:geo:cover-frame}. The chart isomorphisms thus agree on
all overlaps. Reduction modulo \(r\) is pushout of torsors, and
commutes with these constructions.
\end{proof}

\begin{corollary}
\label{cor:geo:norm-period}
At period \(r\), the generator of the \(w\)-valuation cover acts
on the target fibre \(C_r\) by addition by \(f(w/v)\). It has
\(\gcd(f(w/v),r)\) orbits, each of length
\(r/\gcd(f(w/v),r)\). For radial ball functions,
\[
 N^*1_{\{\ord_v\ge n\}}
   =1_{\{\ord_w\ge\lceil n/f(w/v)\rceil\}},
 \qquad N^*t_v^{f(w/v)}=t_wN^*,
\]
where \(t_v\) raises the ball index by \(1\).
\end{corollary}

\begin{proof}
The subgroup generated by \(f(w/v)\) in \(C_r\) has order
\(r/\gcd(f(w/v),r)\). Its cosets are the orbits, so their number
is \(\gcd(f(w/v),r)\). The ball identity follows from
\(\ord_v(Nz)=f(w/v)\ord_w(z)\); replacing \(n\) by
\(n+f(w/v)\) gives the shift identity.
\end{proof}
\section{Spherical normal coefficients}
\label{sec:normal}

All tensor products of spectra and of spectral modules are derived.
Write \(R=\Sph[\Z]=\Sph[t,t^{-1}]\), and let
\(\Sph_{\mathrm{aug}}\) be the \(R\)-module on which \(t\) acts as
the identity.

\begin{definition}[Spherical functions on the valuation line]
\label{def:nor:functions}
Give \(\overline\Z=\Z\cup\{-\infty,+\infty\}\) its order topology.
For a finite clopen partition \(\mathcal P\) of \(\overline\Z\),
form the commutative sphere algebra \(\Sph^{\mathcal P}\). Refinement
acts by diagonal maps on its factors. Define
\[
 \overline A=\mathop{\operatorname{colim}}_{\mathcal P}
                    \Sph^{\mathcal P},
 \qquad
 D=\fib(\mathrm{ev}_{-\infty}:\overline A\to\Sph),
 \qquad
 A_c=\fib(\overline A\to\Sph\times\Sph).
\]
The last map evaluates at the two endpoints. The translation of
\(\overline\Z\) induces an automorphism \(\sigma\), fixed at the
endpoints. For the characteristic section \(b_n\) of the clopen
tail \([n,+\infty]\), the convention is
\(\sigma(b_n)=b_{n+1}\). Put \(\delta_n=b_n-b_{n+1}\).
\end{definition}

The valuation identifies \(K_v/U_v\) with
\(\Z\cup\{+\infty\}\). In this interpretation \(b_n\) represents
the compact-open ball \(\varpi_v^n\mathcal O_v\), and \(\delta_n\)
represents its punctured shell. The missing endpoint \(-\infty\)
is the point at infinity of this locally compact valuation space.

\begin{theorem}[The puncture sequence]
\label{thm:nor:puncture}
The \(R\)-linear maps
\[
 R\longrightarrow A_c,\quad t^n\longmapsto\delta_n,
 \qquad
 R\longrightarrow D,\quad t^n\longmapsto b_n
\]
are equivalences of \(R\)-modules, where \(t\) acts by \(\sigma\).
Under these equivalences, extension of compact-open sections from
the punctured valuation line and evaluation at zero form the cofiber
sequence
\begin{equation}
 \label{eq:nor:puncture}
 R\xrightarrow{\,1-t\,}R\longrightarrow\Sph_{\mathrm{aug}}.
\end{equation}
\end{theorem}

\begin{proof}
Partitions into finitely many singleton integers and two tails form
a cofinal family of clopen partitions. Hence \(\pi_i\overline A\)
consists of functions on \(\overline\Z\) with values in
\(\pi_i\Sph\), constant sufficiently far toward each end. Endpoint
evaluation splits as a map of spectra, using the two tails.
Consequently \(\pi_i A_c\) consists of the finitely supported
functions on \(\Z\), and \(\pi_iD\) consists of the functions
vanishing sufficiently far to the left and constant sufficiently
far to the right. The first have the shell basis \(\delta_n\).
Taking first differences gives the tail basis \(b_n\) for the
second: if the values are \(a_n\), the coefficient of \(b_n\) is
\(a_n-a_{n-1}\). Only finitely many coefficients are nonzero.

These basis descriptions hold for every \(i\). The displayed maps
are obtained from the free \(R\)-module property and the sections
\(\delta_0,b_0\), so they are coherent \(R\)-module maps and are
equivalences. The inclusion sends \(\delta_n\) to
\(b_n-b_{n+1}\), while positive-end evaluation sends every \(b_n\)
to \(1\). This identifies the endpoint fibre sequence with
\eqref{eq:nor:puncture}.
\end{proof}

The multiplication inherited by \(D\) from \(\overline A\) satisfies
\(b_nb_m=b_{\max(n,m)}\). The equivalence with \(R\) in the theorem
uses the translation module structure.

\subsection{The crossed product and its boundaries}

\begin{definition}[Normal algebra]
\label{def:nor:algebra}
The normal algebra is the \(E_1\)-ring spectrum
\[
 B=\overline A\rtimes_\sigma\Z
   =\bigoplus_{n\in\Z}\overline A U^n,
 \qquad
 (fU^n)(gU^m)=f\,\sigma^n(g)U^{n+m}.
\]
The two endpoint evaluations define a homomorphism
\[
 q:B\longrightarrow Q=R_-\times R_+,
 \qquad R_\pm=\Sph[u_\pm,u_\pm^{-1}],
 \qquad q(U)=(u_-,u_+).
\]
Write \(\mathscr K=\fib(q)\) for its two-sided ideal.
\end{definition}

\begin{proposition}[The matrix ideal]
\label{prop:nor:matrix}
There are matrix units
\[
 E_{ij}=\delta_iU^{i-j}\quad(i,j\in\Z),
 \qquad E_{ij}E_{kl}=
 \begin{cases}E_{il},&j=k,\\0,&j\ne k.\end{cases}
\]
They identify \(\mathscr K\) with the sphere algebra of finitely
supported matrices indexed by \(\Z\). If \(e=E_{00}\), then
\[
 eBe\simeq\Sph,
 \qquad Be\otimes_{\Sph}eB\xrightarrow{\ \sim\ }\mathscr K
\]
by multiplication, as \(B\)-bimodules.
\end{proposition}

\begin{proof}
The underlying fibre of \(q\) is
\(\bigoplus_{n\in\Z}A_cU^n\), since coproducts of spectra preserve
fibre sequences. The shell basis of \(A_c\) therefore identifies
its underlying spectrum with \(\bigoplus_{i,j}\Sph E_{ij}\).
The identity \(\sigma^r(\delta_k)=\delta_{k+r}\) gives the displayed
multiplication. On every finite set of indices it is the usual
matrix algebra; these finite matrix algebras give the stated
nonunital algebra by their corner inclusions.

The spectra \(Be\) and \(eB\) have respective bases \(E_{i0}\) and
\(E_{0j}\), while \(eBe\) has the single basis element \(e\).
Multiplication maps \(E_{i0}\otimes E_{0j}\) to \(E_{ij}\), proving
the bimodule equivalence.
\end{proof}

\begin{definition}
A homomorphism \(B\to Q\) of \(E_1\)-rings is \emph{Tor-unital}
if multiplication \(Q\otimes_B Q\to Q\) is an equivalence.
\end{definition}

\begin{proposition}[Boundary restriction]
\label{prop:nor:boundary}
The maps \(B\to Q\) and \(B\to R_+\) are Tor-unital. Their
restriction functors on left-module categories are fully faithful.
The kernel of \(Q\otimes_B-\) is equivalent to \(\mathrm{Sp}\)
under \(V\mapsto Be\otimes_{\Sph}V\), with inverse \(M\mapsto eM\).
\end{proposition}

\begin{proof}
The idempotent \(e\) acts trivially on \(Q\). Proposition
\ref{prop:nor:matrix} gives
\[
 \mathscr K\otimes_B Q
 \simeq Be\otimes_{\Sph}eQ\simeq0.
\]
Tensoring \(\mathscr K\to B\to Q\) with \(Q\) proves
\(Q\otimes_B Q\simeq Q\). The projection \(Q\to R_+\) is an
idempotent factor projection, so the same identity gives
\(R_+\otimes_B R_+\simeq R_+\). For either map, the counit of
extension and restriction of scalars is an equivalence; restriction
is therefore fully faithful.

If \(Q\otimes_B M=0\), the ideal sequence identifies
\(M\simeq\mathscr K\otimes_B M\simeq Be\otimes_{\Sph}eM\).
Conversely \(Q\otimes_B Be=0\). Finally
\(e(Be\otimes_{\Sph}V)\simeq V\), since \(eBe\simeq\Sph\).
\end{proof}

\subsection{Two coefficient lifts and the normal functor}

Identify \(R_+\) with \(R\) by \(u_+=t^{-1}\), and denote the
positive quotient by \(q_+:B\to R\). The homomorphism
\[
 r:R\longrightarrow B,\qquad r(t)=U^{-1}
\]
satisfies \(q_+r=\mathrm{id}_R\).

\begin{definition}[Normal module and coefficient lifts]
\label{def:nor:lifts}
Give \(D\) the right \(B\)-module structure
\[
 d\cdot(fU^n)=\sigma^{-n}(df).
\]
For a left \(R\)-module \(E\), define left \(B\)-modules
\[
 J(E)=B\otimes_{R,r}E,\qquad I(E)=q_{+,*}E.
\]
The quotient \(q_+\) induces \(J(E)\to I(E)\). The normal functor is
\[
 \mathcal N:\Mod_B\longrightarrow\mathrm{Sp},\qquad
 \mathcal N(M)=D\otimes_B M.
\]
\end{definition}

\begin{theorem}[Normal realization of the two lifts]
\label{thm:nor:lifts}
The right \(R\)-module \(D\), obtained by restriction along \(r\),
is free on \(b_0\). There are natural equivalences of spectra
\[
 \mathcal N(J(E))\simeq E,
 \qquad
 \mathcal N(I(E))\simeq\Sph_{\mathrm{aug}}\otimes_R E.
\]
The map induced by \(J(E)\to I(E)\) is the augmentation map in
the triangle obtained by tensoring \eqref{eq:nor:puncture} with \(E\).
\end{theorem}

\begin{proof}
For the module law, the two parenthesizations of multiplication
by \(fU^n\) and \(gU^m\) both give
\[
 \sigma^{-n-m}(df)\,\sigma^{-m}(g).
\]
Restriction gives \(d\cdot r(t)=\sigma(d)\), so Theorem
\ref{thm:nor:puncture} proves freeness. Tensor associativity then
gives \(D\otimes_B J(E)\simeq D\otimes_R E\simeq E\).

For the other lift, multiplication by \(e\) identifies \(De\)
with \(\Sph\). The map
\[
 D\otimes_B\mathscr K
 \simeq De\otimes_{\Sph}eB\longrightarrow D
\]
sends the row basis element \(E_{0j}\) to \(\delta_j\).
It is the inclusion \(A_c\to D\). Consequently
\(D\otimes_B Q\) is \(\Sph\) supported at the positive boundary,
with trivial action of \(u_+\). Passing to the positive factor gives
\(D\otimes_B R\simeq\Sph_{\mathrm{aug}}\). Tensoring over \(R\)
with \(E\) proves the second equivalence. Under these descriptions,
the comparison is positive-end evaluation, which is the final map
of \eqref{eq:nor:puncture}.
\end{proof}

\subsection{Sheaves of normal coefficients}

\begin{definition}[Sheaf conventions and descended normal data]
\label{def:nor:sheaves}
For an open subset \(V\subset Y_K\), let
\(\operatorname{Shv}(V;\mathrm{Sp})\) be the stable category of
sheaves of spectra on its ordinary open-set site. The constant-sheaf
functor is denoted by \(\Delta_V\). A locally constant sheaf in this
section is a sheaf identified with a constant sheaf on some open
cover. Modules over a sheaf of rings are formed in this category.

The principal cover \(Q_v\) and the actions
\[
 n\longmapsto\operatorname{Ad}(U^n)\quad\hbox{on }B,
 \qquad n\longmapsto\sigma^n\quad\hbox{on }D
\]
define by associated descent a sheaf of \(E_1\)-algebras
\(\cB_v\) and a sheaf \(\cD_v\). Let
\(R_{Y_K}=\Delta_{Y_K}R\). The descended maps are written
\[
 R_{Y_K}\xrightarrow{r_v}\cB_v
       \xrightarrow{q_v}R_{Y_K}.
\]
The sheaf \(\cL_v=\Sph_{Y_K}[Q_v]\) is the fibrewise free
sphere module on the torsor, with its right \(R_{Y_K}\)-action
induced by deck translation. For a sheaf \(F\) of left
\(R_{Y_K}\)-modules, put
\[
 J_v(F)=\cB_v\otimes_{R_{Y_K}}F,
 \qquad I_v(F)=q_{v,*}F,
 \qquad \mathcal N_v(M)=\cD_v\otimes_{\cB_v}M.
\]
Here \(\mathcal N_v\) takes values in
\(\operatorname{Shv}(Y_K;\mathrm{Sp})\).
For a constant source module \(E\), write
\(\Loc_v(E)=\cL_v\otimes_{R_{Y_K}}\Delta_{Y_K}E\).
\end{definition}

\begin{theorem}[The descended puncture sequence]
\label{thm:nor:descent}
The sheaf \(\cD_v\) is a right \(\cB_v\)-module. Restriction along
\(r_v\) identifies it with \(\cL_v\). The functor \(I_v\) is fully
faithful, and \(q_vr_v=\mathrm{id}\). For every sheaf \(F\) of
\(R_{Y_K}\)-modules, the comparison \(J_v(F)\to I_v(F)\) gives,
after applying \(\mathcal N_v\), the cofiber sequence
\[
 \cL_v\otimes_{R_{Y_K}}F
   \xrightarrow{1-t}\cL_v\otimes_{R_{Y_K}}F
   \longrightarrow\Sph_{Y_K,\mathrm{aug}}
                         \otimes_{R_{Y_K}}F.
\]
In particular,
\begin{equation}
 \label{eq:nor:global-puncture}
 \Loc_v(E)\xrightarrow{1-t}\Loc_v(E)
 \longrightarrow
 \Delta_{Y_K}(\Sph_{\mathrm{aug}}\otimes_R E).
\end{equation}
For \(E=\Sph[C_L]\), the first term is the fibrewise suspension
spectrum of \(Q_{v,L}\).
\end{theorem}

\begin{proof}
The identity
\[
 \sigma^n(d\cdot b)
  =\sigma^n(d)\cdot\operatorname{Ad}(U^n)(b)
\]
proves compatibility of the module structure with the transition
maps. Conjugation by \(U^n\) fixes \(r(R)\), and the endpoint map
is invariant under conjugation. Both algebra maps therefore descend.
In the ball basis, the transition of \(D\) is multiplication by
\(t^n\), which is the transition of the free sphere module on
\(Q_v\). This proves \(\cD_v|_{R_{Y_K}}\simeq\cL_v\).

Tor-unitality is an equivalence of sheaves of bimodules and can be
checked on an open cover trivializing \(Q_v\). Proposition
\ref{prop:nor:boundary} therefore proves full faithfulness of \(I_v\).
The local comparison of Theorem \ref{thm:nor:lifts} descends, since
positive-end evaluation and the shell inclusion commute with the
transition maps. Its augmentation quotient has trivial transition,
as \(t^n\) acts trivially on \(\Sph_{\mathrm{aug}}\). Tensoring
the descended sequence with \(F\) proves the claim. The final
assertion follows by forming the associated permutation spectrum
of \(Q_v\) with fibre \(C_L\).
\end{proof}

\subsection{Several places}

For a finite set \(S\) of finite places, use separate deck variables
\[
 R_S=\Sph[\Z^S]=\Sph[t_v^{\pm1}:v\in S],
 \qquad
 \cB_{v,S}=\cB_v\otimes_{\Sph_{Y_K}}
       \Delta_{Y_K}\Sph[\Z^{S\setminus\{v\}}].
\]
The maps \(r_v,q_v\) extend by the identity on the other variables;
denote the resulting maps by \(r_{v,S},q_{v,S}\).

\begin{definition}[Coefficients at a finite set of places]
\label{def:nor:coefficients}
For \(V\subset Y_K\) open, write \(R_{S,V}=\Delta_V R_S\) and put
\[
 \mathcal R_S(V)=\Mod_{R_{S,V}}
                       (\operatorname{Shv}(V;\mathrm{Sp})),
 \qquad
 \mathcal M_{v,S}(V)=\Mod_{\cB_{v,S}|_V}
                       (\operatorname{Shv}(V;\mathrm{Sp})).
\]
Define the stable category
\begin{equation}
 \label{eq:nor:coefficients}
 \cE_S(V)=\mathcal R_S(V)
       \mathbin{\times}_{\prod_{v\in S}\mathcal R_S(V)}
                         \prod_{v\in S}\mathcal M_{v,S}(V),
\end{equation}
where the product is a homotopy fibre product. Its two arrows are
the diagonal and the boundary functors
\(M_v\mapsto R_{S,V}\otimes_{\cB_{v,S}|_V}M_v\).
An object is a tuple
\[
 \bigl(F,(M_v)_{v\in S},(\eta_v)_{v\in S}\bigr),
 \qquad
 \eta_v:R_{S,V}\otimes_{\cB_{v,S}|_V}M_v\xrightarrow{\sim}F.
\]
Let \(p_S\) be its projection to \(F\). The sections \(J_S,I_S\)
of \(p_S\) use, respectively,
\(M_v=\cB_{v,S}|_V\otimes_{R_{S,V}}F\) and
\(M_v=q_{v,S,*}F\), with their canonical boundary equivalences.
\end{definition}

\begin{proposition}
\label{prop:nor:coefficients}
The categories \(\cE_S(V)\) satisfy descent as \(V\) varies.
The section \(I_S\) is fully faithful. The quotient maps define a
natural transformation \(J_S\to I_S\), and its normal realization
at each \(v\) is the descended puncture sequence, with the other
variables retained.

For \(S\subset T\), extension of scalars defines coherent functors
\(\cE_S(V)\to\cE_T(V)\) as follows: extend the common boundary
from \(R_S\) to \(R_T\), extend the existing components along
\(\cB_{v,S}\to\cB_{v,T}\), and use \(J\) for each new component.
These functors preserve the \(J_S\)-sections. On the
\(I_S\)-sections they carry a natural transformation to \(I_T\),
given on the new components by \(J\to I\).
\end{proposition}

\begin{proof}
Sheaves of modules are obtained by descent on open covers.
Taking the homotopy limit in \eqref{eq:nor:coefficients} commutes
with these descent limits. Each boundary functor is exact, so the
resulting category is stable. Full faithfulness of \(I_S\) follows
from full faithfulness of boundary restriction on every component:
the mapping-space fibre product for two objects in its image
reduces to the mapping space of their common boundary modules.
The comparison and normal sequence follow componentwise from
Theorem \ref{thm:nor:descent}.

For \(S\subset T\), write \(F_T=R_{T,V}\otimes_{R_{S,V}}F\).
For an existing component take
\(\cB_{v,T}|_V\otimes_{\cB_{v,S}|_V}M_v\); its boundary is
\(F_T\) by tensor associativity and \(\eta_v\). For a new place
\(w\in T\setminus S\), take
\(\cB_{w,T}|_V\otimes_{R_{T,V}}F_T\), whose boundary is again
\(F_T\). The added variables are central tensor factors, so these
operations compose under \(S\subset T\subset W\). Tensor
associativity also identifies the extension of \(J_S(F)\) with
\(J_T(F_T)\). For \(I_S(F)\), the old components become the
corresponding restricted boundary modules over \(R_T\); on the
new components the quotient map supplies the stated natural
transformation. Its coherence follows from that of the same
quotient maps and tensor products.
\end{proof}
\section{Angular covers and residue coefficients}
\label{sec:residue}

Fix a finite place \(v\) of \(K\). Write
\[
 F=K_v,\quad \mathcal O=\mathcal O_v,\quad
 \mathfrak m=\mathfrak m_v,\quad
 U=\mathcal O^\times,\quad U^r=1+\mathfrak m^r\ (r\geq1),
\]
and let \(k=\mathcal O/\mathfrak m\), \(q=\lvert k\rvert=p^d\).
The space \(Y_K\), its normal chart \(N_v\), and the principal
\(\Z\)-cover \(Q_v\to Y_K\) are those of
Definitions~\ref{def:geo:atlas} and~\ref{def:geo:valuation-cover}.

\subsection{The unit-cover system}

\begin{definition}[Angular covers]
\label{def:res:angular}
For an open subgroup \(H\subset U\), set
\[
 G_{v,H}=F^\times/H,\qquad
 N_{v,H}=(C_K\times F)/H,
\]
using the action \((c,z)h=(c j_v(h),h^{-1}z)\).
The group \(G_{v,H}\) is discrete and fits into the sequence
\begin{equation}
 1\longrightarrow U/H\longrightarrow G_{v,H}
 \xrightarrow{\ord_v}\Z\longrightarrow0.
 \label{eq:res:group-extension}
\end{equation}
Define \(\mathcal Q_{v,H}\) over \(Y_K\) by taking \(N_{v,H}\) over
\(N_v\), taking \(N_w\times G_{v,H}\) over \(N_w\) for \(w\ne v\),
and gluing over the common open \(C_K\) by
\begin{equation}
 [c,z]_H\longmapsto\bigl(cj_v(z),[z]^{-1}\bigr),
 \qquad z\ne0.
 \label{eq:res:punctured-frame}
\end{equation}
\end{definition}

\begin{proposition}
\label{prop:res:angular-covers}
The projection \(\mathcal Q_{v,H}\to Y_K\) is a principal
\(G_{v,H}\)-cover, with right deck action
\[
 [c,z]_H\cdot[h]=[c j_v(h),h^{-1}z]_H.
\]
An inclusion \(H'\subset H\) induces a coherent quotient map of these
covers. Extension of structure group through \(\ord_v\) gives
\(Q_v\); equivalently, \(\mathcal Q_{v,H}\to Q_v\) is a finite
principal \(U/H\)-cover.
\end{proposition}
\begin{proof}
An open subgroup of the compact group \(U\) is compact and of finite
index. The free proper \(F^\times\)-action on \(C_K\times F\)
descends to a free proper action of \(F^\times/H\) on \(N_{v,H}\).
For properness, lift compact subsets through the compact-group quotient
\(C_K\times F\to N_{v,H}\) and use properness of the original action.
A free proper action of a discrete group on a locally compact Hausdorff
space is a covering action: start with a relatively compact neighborhood,
discard the finitely many translates that can meet it, and shrink it to
be disjoint from all its nonidentity translates.

The first coordinate in \eqref{eq:res:punctured-frame} is invariant
under the deck group, and the second is multiplied by \([h]\).
Thus it is an equivariant trivialization. It glues the local cover to
the trivial covers on the other charts. The quotient maps for different
\(H\) agree in these trivializations, which proves their coherence.
Taking the quotient by \(U/H\) gives \((C_K\times F)/U\) on \(N_v\),
and the second coordinate becomes \(-\ord_v(z)\). This is the
transition function of \(Q_v\).
\end{proof}

On the boundary \(Z_v=C_K/j_v(F^\times)\), these maps restrict to
\(C_K/j_v(H)\to C_K/j_v(U)\to C_K/j_v(F^\times)\).
Thus the finite unit quotient and the valuation quotient arise from the
same sequence \eqref{eq:res:group-extension}.

\begin{definition}[A marked tame--wild refinement]
\label{def:res:frame}
Set \(m=q-1\). Choose a generator \(u\) of \(k^\times\) and a nonzero
homomorphism
\[
 \chi:U^1/U^2\longrightarrow C_p.
\]
Let \(H_\chi\subset U^1\) be the inverse image of its kernel. The pair
\((u,\chi)\) is a marking of the refinement
\begin{equation}
 \mathcal Q_{v,H_\chi}\longrightarrow
 \mathcal Q_{v,U^1}\longrightarrow Q_v.
 \label{eq:res:cyclic-refinement}
\end{equation}
\end{definition}

\begin{proposition}
\label{prop:res:cyclic-refinement}
Over \(Q_v\), the covers \(\mathcal Q_{v,H_\chi}\) and
\(\mathcal Q_{v,U^1}\) have marked deck groups \(C_{mp}\) and
\(C_m\), respectively. The map between them is induced by
\(C_{mp}\to C_m\), with kernel \(C_p\). Every quotient by a
subgroup of either cyclic group is a finite cover over \(Q_v\),
and its permutation spectra carry the transfers of
Definition~\ref{def:src:periods}.
\end{proposition}
\begin{proof}
The map \(1+x\mapsto x\bmod\mathfrak m^2\) identifies
\(U^1/U^2\) with the additive group of \(\mathfrak m/\mathfrak m^2\).
It is a \(d\)-dimensional vector space over \(\F_p\), so \(\chi\) is
surjective. Hensel lifting of the simple roots of \(X^m-1\) gives
\(U=\mu_m\times U^1\). Hence
\[
 U/H_\chi=\mu_m\times C_p,
 \qquad U/U^1=k^\times.
\]
Let \([u]\in\mu_m\) be the lift of \(u\), and let \(w_\chi\) be
the class in \(U^1/H_\chi\) with \(\chi(w_\chi)=1\). The element
\([u]w_\chi\) has order \(mp\) and maps to \(u\); it supplies the
two marked generators. Apply Proposition~\ref{prop:res:angular-covers}
and the finite-set transfer construction.
\end{proof}

\subsection{A sheaf of primitive contacts}

A principal finite-group cover and an action on an object \(M\)
determine an associated sheaf by descent along a trivializing open cover,
using the constant-sheaf convention of Section~\ref{sec:notation}.

\begin{definition}[The angular contact sheaf]
\label{def:res:contact-sheaf}
Let \(P_{mp}\) and \(P_m\) denote the two finite torsors over \(Q_v\)
in \eqref{eq:res:cyclic-refinement}. Their primitive permutation-spectrum
sheaves define
\[
 \mathcal T_v=U(P_m)\otimes_{\Delta_{Q_v}R}U(P_{mp}).
\]
It is a locally constant sheaf with values in
\(\Fun(BC_m,\mathrm{Sp})\), where \(C_m\) acts through the first
factor. At a point \(x\in Q_v\), write \(\mathcal T_{v,x}\) for its
fibre in that category.
\end{definition}

\begin{lemma}
\label{lem:res:contact-transport}
A choice of origin in a fibre of \(P_{mp}\) identifies
\(\mathcal T_{v,x}\) with \(T_{m,mp}\). Changing that origin by \(r\)
acts on the contact by \(t^{2r}\). The induced transport on
\(\pi_0\End_{C_m}(\mathcal T_{v,x})\) is the identity.
\end{lemma}
\begin{proof}
An origin in \(P_{mp}\) gives a compatible origin in its quotient
\(P_m\). These identify both primitive fibres with those of
Definition~\ref{def:src:primitive}. Translation of the origin by \(r\)
acts on the two factors by \(t^r\otimes t^r\), which is \(t^{2r}\)
after balancing over \(R\). A \(C_m\)-equivariant morphism commutes
with this power of the marked rotation, so its conjugation action on
the degree-zero endomorphism algebra is trivial.
\end{proof}

\begin{theorem}[Residue coefficients of the contact]
\label{thm:res:residue}
Let \(u\) be the generator in Definition~\ref{def:res:frame}. The
homomorphism \(\F_p[t]/(\Phi_m)\to k\), \(t\mapsto u\), selects a
factor of the fibre endomorphism algebra. This factor is
\[
 k\quad\text{if }p\ne2,
 \qquad
 W_2(k)\quad\text{if }p=2,
\]
where \(W_2(k)=W(k)/p^2W(k)\). In the second case its reduction modulo
\(2\) is \(k\).
\end{theorem}
\begin{proof}
Theorem~\ref{thm:src:endomorphisms} gives the fibre algebra
\(E_m=R_p[t]/(\Phi_m)\) before selecting a factor. Because
\(p\nmid m\), the reduction of \(\Phi_m\) is separable. If \(m>1\),
the order of \(p\) modulo \(m=p^d-1\) is \(d\): the exponent \(d\)
works, whereas \(0<p^r-1<m\) for \(0<r<d\). The Frobenius orbit of
\(u\) therefore has size \(d\), and its irreducible polynomial gives
the quotient \(k\). For odd \(p\), Chinese remainders select this
field as a factor of \(E_m\).

For \(p=2\), the primitive idempotent lifts uniquely through the
nilpotent ideal \(2E_m\). Its factor is the finite etale \(\Z/4\)-algebra
with residue field \(k\), and is consequently \(W_2(k)\), with marked
root the Teichmuller lift of \(u\). One can construct it by the coprime
Hensel factorization of \(\Phi_m\); uniqueness also follows from the
finite etale equivalence for henselian local rings
\cite[Tag~04GK]{Stacks}. If \(m=1\), then \(q=2\) and the calculation
is \(E_1=\Z/4=W_2(\F_2)\).
\end{proof}

\begin{definition}[The selected contact sheaf]
\label{def:res:selected}
Let \(A=\Sph_p^\wedge[C_m]\). The primitive idempotent of
\(\F_p[C_m]\) indexed by the Frobenius orbit of \(u\) lifts to a
central idempotent \(\epsilon_u\in\Z_p[C_m]=\pi_0A\). Let
\(A_{\epsilon_u}\) be its idempotent factor, and put
\begin{equation}
 \mathcal T_{v,u}
 =\Delta_{Q_v}(A_{\epsilon_u})
       \otimes_{\Delta_{Q_v}A}\mathcal T_v.
 \label{eq:res:selected}
\end{equation}
Here \(\mathcal T_v\) is first formed as a locally constant sheaf of
\(A\)-modules in \(p\)-complete spectra.
\end{definition}

\begin{proposition}
\label{prop:res:selected}
The sheaf \(\mathcal T_{v,u}\) is a retract of \(\mathcal T_v\).
For every \(x\in Q_v\),
\begin{equation}
 \pi_0\End_{\Fun(BC_m,\mathrm{Sp})}(\mathcal T_{v,u,x})
 \simeq
 \begin{cases}
  k,&p\ne2,\\
  W_2(k),&p=2.
 \end{cases}
 \label{eq:res:fibre-end}
\end{equation}
The associated locally constant sheaf of degree-zero endomorphism
algebras has trivial transport and is constant with this value.
\end{proposition}
\begin{proof}
Every contact fibre is a finite wedge of Moore spectra by
Proposition~\ref{prop:src:contact-spectrum}, and is \(p\)-complete.
Its canonical \(\Sph_p^\wedge\)-action commutes with the retained
\(C_m\)-action and with the torsor transition maps. Trivialize the
finite torsors over an open cover, form the contact in
\(p\)-complete spectra on each open set, and descend the resulting
\(A\)-modules. This constructs the \(A\)-module sheaf used in
\eqref{eq:res:selected}.

Hensel lifting gives \(\epsilon_u\) as a central group-algebra
idempotent. Its action commutes coherently with transport, so extension
to the factor \(A_{\epsilon_u}\) splits the sheaf. On every fibre it
is the primitive projector of Theorem~\ref{thm:res:residue}. Mapping
spectra between \(p\)-complete objects agree with their mapping spectra
in spectra. The corner of the indexed endomorphism algebra is therefore
the selected factor in \eqref{eq:res:fibre-end}. The transport
calculation of Lemma~\ref{lem:res:contact-transport} descends to this
factor and is the identity.
\end{proof}

The completed primitive coefficient of
Proposition~\ref{prop:src:completed}, with the same choice of \(u\), is
\(W(k)\). The unique unramified subextension of \(F/\Q_p\) with
residue field \(k\) is \(\Q_p(\mu_{q-1})\)
\cite[Proposition~7.50 and Corollary~7.51]{MilneANT}.
In its ring of integers \(W(k)\subset\mathcal O\), the element \(p\) is a
uniformizer and \(\mu_{q-1}\) reduces onto \(k^\times\) by Hensel's lemma.
Thus Nakayama's lemma gives \(W(k)=\Z_p[\mu_{q-1}]\).
The ring \(\mathcal O\) is finite free and faithfully flat over this
unramified coefficient ring. For \(p=2\), reduction therefore gives an
injection \(W_2(k)\hookrightarrow\mathcal O/4\mathcal O\). Thus the finite
contact factor and the completed primitive factor have compatible
reductions in the normal coordinate ring.

\begin{proposition}[Change of marking and Frobenius]
\label{prop:res:frobenius}
For \(a\in(\Z/m\Z)^\times\), the residue comparisons satisfy
\[
 \operatorname{ev}_{u^a}
 =\operatorname{ev}_{u}\circ(t\mapsto t^a).
\]
Residue Frobenius is induced by the automorphism
\(t\mapsto t^p\) of the marked deck group \(C_m\).
This extends to an automorphism of the deck-group quotient
\(C_{mp}\to C_m\), with exponent
\(A\in(\Z/mp\Z)^\times\) on \(C_{mp}\) determined by
\[
 A\equiv p\pmod m,\qquad A\equiv1\pmod p.
\]
\end{proposition}
\begin{proof}
The first identity is evaluation of a polynomial. Modulo \(p\), raising
the marked root to \(p\) raises every element of its finite-field factor
to \(p\). Since \(p\nmid m\), the two congruences have a unique
solution modulo \(mp\); it is prime to both \(m\) and \(p\). On the
decomposition \(U/H_\chi=\mu_m\times C_p\), it acts by \(p\)-th
power on the tame factor and by the identity on the chosen wild quotient.
\end{proof}

\subsection{Norms, residue degrees and ramification}

\begin{proposition}
\label{prop:res:local-norms}
Let \(L_w/K_v\) be a finite extension, with residue degree \(f\) and
ramification index \(e\). For \(x\in L_w^\times\),
\(b\in\mathcal O_w^\times\), and \(a\in\mathcal O_w\),
\begin{align}
 \ord_v(N_{L_w/K_v}x)&=f\,\ord_w(x),
 \label{eq:res:valuation-norm}\\
 \overline{N_{L_w/K_v}b}&=N_{k_w/k_v}(\bar b)^e,
 \label{eq:res:unit-norm}\\
 \overline{\operatorname{Tr}_{L_w/K_v}a}
 &=e\,\operatorname{Tr}_{k_w/k_v}(\bar a).
 \label{eq:res:residue-trace}
\end{align}
\end{proposition}
\begin{proof}
For an integral nonzero \(x\), multiplication by \(x\) on the free
\(\mathcal O_v\)-module \(\mathcal O_w\) has determinant
\(N_{L_w/K_v}x\). The valuation of the determinant is
\(\operatorname{length}_{\mathcal O_v}(\mathcal O_w/x\mathcal O_w)\).
The latter quotient has \(\ord_w(x)\) successive factors isomorphic
to \(k_w\), each of \(\mathcal O_v\)-length \(f\). This proves
\eqref{eq:res:valuation-norm}; division extends it to \(L_w^\times\).

The \(k_v\)-vector space \(\mathcal O_w/\varpi_v\mathcal O_w\) has
the filtration by powers of \(\mathfrak m_w\), with \(e\) successive
quotients isomorphic to \(k_w\). Multiplication by \(a\) induces
multiplication by \(\bar a\) on each quotient. Traces add and
determinants multiply along this filtration. Reducing the integral
multiplication matrix modulo \(\varpi_v\) therefore gives
\eqref{eq:res:unit-norm} and~\eqref{eq:res:residue-trace}.
\end{proof}

\begin{theorem}[Norms of the covering system]
\label{thm:res:cover-norms}
Let \(L/K\) be a finite extension and \(N:Y_L\to Y_K\) the map of
Theorem~\ref{thm:geo:norm}. Fix a finite place \(v\) and an open
\(H_v\subset\mathcal O_v^\times\). For each \(w\mid v\), choose an
open \(H_w\subset\mathcal O_w^\times\) such that
\(N_{L_w/K_v}(H_w)\subset H_v\), and write
\[
 \nu_w:L_w^\times/H_w\longrightarrow K_v^\times/H_v
\]
for the norm homomorphism. Then there is an isomorphism of principal
\(K_v^\times/H_v\)-covers on \(Y_L\)
\begin{equation}
 N^*\mathcal Q_{K,v,H_v}
 \simeq\mathop{\bigotimes}_{w\mid v}
                 \nu_{w,*}\mathcal Q_{L,w,H_w}.
 \label{eq:res:norm-torsors}
\end{equation}
The right side is the contracted product of torsors for the abelian
target group. These isomorphisms commute with refinement of the open
subgroups and with towers of finite extensions.
\end{theorem}
\begin{proof}
On the chart at \(w\mid v\), the map
\[
 [c,z]_{H_w}\longmapsto[N_{L/K}c,N_{L_w/K_v}z]_{H_v}
\]
is well-defined and \(\nu_w\)-equivariant. It identifies the extension
of structure group of the source torsor with the pullback of the
target torsor. All factors associated to places other than \(w\)
are trivial on this chart. On a chart over a place not dividing \(v\),
all the torsors in \eqref{eq:res:norm-torsors} are trivial.
In the punctured trivializations, the transition label satisfies
\([z]^{-1}\mapsto[Nz]^{-1}\). Multiplication of these labels over
\(w\mid v\) is the target label, so the chartwise isomorphisms glue.
Quotienting subgroups and composing norms preserve these labels, which
proves both compatibilities.
\end{proof}

Such subgroups \(H_w\) exist by continuity of the norm. In particular,
one may refine any chosen \(H_w\) by intersection with
\(N_{L_w/K_v}^{-1}(H_v)\). Applying valuation to
\eqref{eq:res:norm-torsors} gives the coefficient \(f(w/v)\) in the
pullback of the valuation covers, while reduction on the unit groups
gives the ramification exponent \(e(w/v)\) in
\eqref{eq:res:unit-norm}.
\section{Correspondence coefficients and descent}
\label{sec:cor:correspondence}

Fix a Haar measure \(dc\) on \(C=C_K\) whose pushforward under
\(\log|\cdot|:C\to\R\) is Lebesgue measure. For \(a\in C\), let
\(\gamma_a:Y_K\to Y_K\) be the action defined by
\(\gamma_a[c,z]=[ac,z]\). Its graph is the span
\[
 Y_K\xleftarrow{\mathrm{id}}Y_K\xrightarrow{\gamma_a}Y_K.
\]
Composition of these spans gives the graph of \(\gamma_{ab}\).

\begin{definition}[The convolution algebra]
\label{def:cor:algebra}
Let \(\cD(C)\) be the Bruhat space of compactly supported smooth complex
functions on \(C\). Thus \(\cD(C)\) is the inductive limit of
\(C_c^\infty(C/U)\), pulled back to \(C\), over compact subgroups \(U\)
for which \(C/U\) is a Lie group. Each stage has its usual inductive
bornology. Put
\[
 B_C=\C[C^{\mathrm{disc}}],\qquad
 \mathscr D_K=\cD(C)\,dc,
 \qquad \mathscr A_K=B_C\oplus\mathscr D_K.
\]
An element of \(B_C\) is a finite linear combination of Dirac measures;
\(\mathscr D_K\) consists of the measures with a Bruhat test density.
Multiplication in \(\mathscr A_K\) is convolution of measures. The
bornology on \(\mathscr A_K\) is the direct sum of the fine bornology on
\(B_C\) and the Bruhat bornology on \(\mathscr D_K\).
\end{definition}

\begin{proposition}
\label{prop:cor:split}
The space \(\mathscr A_K\) is a unital commutative algebra of compactly
supported measures. Its unit is \(\delta_1\), \(\mathscr D_K\) is an
ideal, and the atomic projection
\(\pi:\mathscr A_K\to B_C\) is a split algebra homomorphism.
Multiplication is bounded for the stated bornology.
\end{proposition}
\begin{proof}
The norm quotient makes \(C\) nondiscrete, so Haar measure has no atoms.
Consequently an atomic measure that also has a smooth density is zero.
If \(\alpha,\beta\in B_C\) and \(f,g\in\cD(C)\), then
\[
 (\alpha+f\,dc)*(\beta+g\,dc)
 =\alpha*\beta+
     (\alpha*g+f*\beta+f*g)\,dc.
\]
Translation and convolution preserve \(\cD(C)\); hence the last three
terms are in \(\mathscr D_K\). This proves the algebra and ideal
statements. The atomic inclusion splits \(\pi\). A bounded subset of
\(B_C\) is contained in a bounded subset of a finite-dimensional span
of Dirac measures. Translation by the corresponding finite set is
bounded on the Bruhat space, and convolution is bounded on each common
smooth compact-support stage. These observations prove boundedness.
\end{proof}

The regular action on test functions is
\begin{equation}
 \rho(\mu)u(x)=\int_Cu(a^{-1}x)\,d\mu(a).
 \label{eq:cor:regular}
\end{equation}
Fubini's theorem gives \(\rho(\mu*\nu)=\rho(\mu)\rho(\nu)\).
This representation is faithful: convolution with a smooth approximate
identity converges to a compactly supported measure in the distribution
topology, so \(\mu*u=0\) for every \(u\in\cD(C)\) implies \(\mu=0\).

\begin{definition}[Spherical correspondence coefficients]
\label{def:cor:spherical}
Write \(\Sph_C=\Sph[C^{\mathrm{disc}}]\). Define the commutative ring
spectrum
\begin{equation}
 \mathscr A_K^{\mathrm{sph}}
  :=\Sph_C\mathbin{\times}_{HB_C}H\mathscr A_K,
 \label{eq:cor:spherical}
\end{equation}
where the first map is induced by \(\Sph\to H\C\) and the second by
\(\pi\). The pullback is taken in commutative ring spectra. Set
\(\mathrm{Corr}_K=\Perf(\mathscr A_K^{\mathrm{sph}})\).
\end{definition}

\begin{lemma}
\label{lem:cor:splitting}
The projection \(\mathscr A_K^{\mathrm{sph}}\to\Sph_C\) has a section
as a map of commutative ring spectra. As \(\Sph_C\)-modules,
\[
 \mathscr A_K^{\mathrm{sph}}
    \simeq \Sph_C\oplus H\mathscr D_K.
\]
The product on the second summand is convolution and its
\(\Sph_C\)-action factors through the translation action of \(B_C\).
\end{lemma}
\begin{proof}
The atomic inclusion \(B_C\to\mathscr A_K\), together with the identity
of \(\Sph_C\), supplies the section of the pullback projection. The
forgetful functor from commutative ring spectra to spectra preserves
limits. Its fiber is therefore
\(\fib(H\mathscr A_K\to HB_C)=H\mathscr D_K\). The section splits
this fiber sequence as modules. The multiplication is the one induced
from \(\mathscr A_K\).
\end{proof}

Choose one uniformizer \(\varpi_v\) at every finite place and write
\(h_v=j_v(\varpi_v)\). For a finite set \(S\) of finite places, put
\[
 R_S=\Sph[\Z^S],\qquad
 H_S=\langle h_v:v\in S\rangle\subset C.
\]
The map \(\Z^S\to C\) is injective. Indeed, if the idele with
coordinates \(\varpi_v^{n_v}\) on \(S\) and \(1\) elsewhere is a
diagonal element \(k\in K^\times\), comparison at a place outside
\(S\) gives \(k=1\), and then every \(n_v\) is zero.

\begin{theorem}[Descent for winding coefficients]
\label{thm:cor:descent}
The homomorphism
\[
 R_S\longrightarrow\mathscr A_K^{\mathrm{sph}},
 \qquad t_v\longmapsto\delta_{h_v},
\]
admits an \(R_S\)-linear retraction. It is descendable, and the
extension-of-scalars functor induces an equivalence
\begin{equation}
 \Mod_{R_S}\simeq
 \Tot_{[n]\in\Delta}
 \Mod_{(\mathscr A_K^{\mathrm{sph}})^{\otimes_{R_S}(n+1)}}.
 \label{eq:cor:descent}
\end{equation}
For a tuple of positive integers \(\mathbf m=(m_v)_{v\in S}\), put
\[
 R_{S,\mathbf m}=\Sph\bigl[\textstyle\prod_{v\in S}C_{m_v}\bigr],
 \qquad
 A_{S,\mathbf m}=R_{S,\mathbf m}\otimes_{R_S}\mathscr A_K^{\mathrm{sph}}.
\]
There is an analogous equivalence with \(R_{S,\mathbf m}\) and
\(A_{S,\mathbf m}\). These equivalences identify mapping spectra and
restrict to perfect objects.
\end{theorem}
\begin{proof}
By Lemma~\ref{lem:cor:splitting}, it suffices to retract
\(R_S\to\Sph_C\) as modules. The decomposition of \(C\) into
\(H_S\)-cosets gives a direct-sum decomposition of \(\Sph_C\) as
an \(R_S\)-module. Projection onto the coset \(H_S\) is
\(R_S\)-linear and its target is \(\Sph[H_S]\simeq R_S\).
Thus the unit of \(\Mod_{R_S}\) belongs to the thick tensor ideal
generated by \(\mathscr A_K^{\mathrm{sph}}\). This is descendability
in the sense of \cite[Definition 3.18]{Mathew}, and
\cite[Proposition 3.22]{Mathew} gives
\eqref{eq:cor:descent}. Tensoring the retraction with
\(R_{S,\mathbf m}\) proves the assertion after finite cyclic base
change. Mapping spectra are part of the equivalence of stable
\(\infty\)-categories. Perfectness is detected by this base change
by \cite[Proposition 3.28]{Mathew}.
\end{proof}

In particular, indexed contact endomorphism spectra are recovered with
their Amitsur descent data. The fixed uniformizers make these
constructions compatible with inclusions of finite place sets.

\begin{proposition}[Cospans and coefficient sheaves]
\label{prop:cor:cospans}
For each finite place \(v\), the source realization of graph cospans,
followed by scalar extension along \(t\mapsto\delta_{h_v}\), gives
\begin{equation}
 \mathcal C_{\tau,g}\longrightarrow
 h\operatorname{Cospan}(\Perf_R)\longrightarrow
 h\operatorname{Cospan}(\mathrm{Corr}_K).
 \label{eq:cor:source}
\end{equation}
It preserves homotopy-pushout composition, the covering maps and
canonical transfers in the source, and their tower and pullback
homotopies.
For every open \(U\subset Y_K\), the coefficient category of
Definition~\ref{def:nor:coefficients} has an exact functor
\begin{equation}
 \Xi_{S,U}:\cE_S(U)\longrightarrow
       \Mod_{\Delta_U\mathscr A_K^{\mathrm{sph}}}(U),
 \qquad
 (F,(M_v),(\eta_v))\longmapsto
 \Delta_U\mathscr A_K^{\mathrm{sph}}
       \otimes_{\Delta_U R_S}F.
 \label{eq:cor:boundary-functor}
\end{equation}
These functors commute with restriction to open subsets. On the source
section \(J_S\), they agree with sheafification of the scalar extension
in \eqref{eq:cor:source}.
\end{proposition}
\begin{proof}
Derived scalar extension preserves colimits and is exact. It therefore
preserves the pushout that composes two cospans, as well as the source
maps and their homotopies. In \eqref{eq:cor:boundary-functor} the first
step is the exact boundary projection
\(p_S:\cE_S(U)\to\Mod_{\Delta_U R_S}(U)\).
Tensoring with the constant sheaf of \(\mathscr A_K^{\mathrm{sph}}\)
is exact and commutes with inverse image along an open inclusion.
Finally \(p_SJ_S=\mathrm{id}\), which proves the last assertion.
\end{proof}
\section{Schwartz realization of the normal coordinates}
\label{sec:ana:normal}

Let \(\Gamma=K^\times\). Write
\(\widetilde j_v:K_v^\times\to\A_K^\times\) for the idele embedding
with all other coordinates \(1\); its composite with the class map is
\(j_v\). For a finite place \(v\), put
\[
 I^{(v)}=\prod_{w\ne v}'K_w^\times,\qquad
 \A_K^{(v)}=\prod_{w\ne v}'K_w.
\]
The first restricted product has the idele topology and the second the
additive adele topology. If \(b\in I^{(v)}\), let \(c(b)\in C\) be the
class of the idele with \(v\)-coordinate \(1\) and remaining coordinates
\(b\).

\begin{proposition}[Normal quotient and transition function]
\label{prop:ana:normal-quotient}
There is a \(C\)-equivariant homeomorphism of quotient presentations
\begin{equation}
 (K_v\times I^{(v)})/\Gamma\ \xrightarrow{\ \sim\ }N_v,
 \qquad [(z,b)]\longmapsto[c(b),z],
 \label{eq:ana:normal-quotient}
\end{equation}
where \(k\in\Gamma\) acts by
\((z,b)\mapsto(k_vz,k^{(v)}b)\).
The pullback of the valuation cover \(Q_v\to N_v\) has a section
\(s(z,b)=[c(b),z]_{U_v}\), with \(U_v=\mathcal O_v^\times\), satisfying
\begin{equation}
 s(k_vz,k^{(v)}b)=s(z,b)\cdot(-\ord_v k).
 \label{eq:ana:frame-cocycle}
\end{equation}
For the associated-bundle convention \((p+n,x)\sim(p,t_v^n x)\),
coefficients in this moving section transform by \(t_v^{\ord_v k}\).
\end{proposition}
\begin{proof}
An idele lift \((a_v,b)\) of \(c\) gives
\([c,z]=[c(b),a_vz]\) in the balanced product defining \(N_v\).
Changing the lift by \(k\in K^\times\) gives the indicated diagonal
relation. Both constructions are continuous quotient maps and are
inverse to one another. The identity
\(c(k^{(v)}b)=c(b)j_v(k_v^{-1})\) implies
\eqref{eq:ana:frame-cocycle}, because the action of \(k_v^{-1}\) on
\(Q_v\) has valuation \(-\ord_v k\). If the section changes by this
deck translation, an invariant element of the associated bundle has
its coefficient multiplied by the inverse translation. This gives
\(t_v^{\ord_v k}\).
\end{proof}

We use the standard Bruhat--Schwartz spaces \(\cS(K_v)\) and
\(\cS(\A_K)\). At a finite place \(\cS(K_v)\) is the space of
locally constant compactly supported functions with its inductive
bornology; at an infinite place it is the usual Schwartz space with
its Fr\'echet bornology. The space \(\cS(\A_K)\) is the restricted
completed tensor product with distinguished vectors
\(1_{\mathcal O_v}\) at the finite places. The dilation convention is
\begin{equation}
 (\lambda_a\Phi)(x)=\Phi(a^{-1}x).
 \label{eq:ana:dilation}
\end{equation}

\begin{proposition}[The normal ball module]
\label{prop:ana:ball-module}
Let \(R_{v,\C}=\C[t_v,t_v^{-1}]\), with its fine bornology.
There is a bornological isomorphism
\begin{equation}
 \iota_v:R_{v,\C}\widehat\otimes\cS(\A_K^{(v)})
 \xrightarrow{\ \sim\ }\cS(\A_K)^{U_v},
 \qquad
 t_v^n\otimes\Psi\longmapsto
 1_{\varpi_v^n\mathcal O_v}\otimes\Psi.
 \label{eq:ana:ball-module}
\end{equation}
It intertwines multiplication by \(t_v\) with
\(\lambda_{\widetilde j_v(\varpi_v)}\). The diagonal action of \(k\in\Gamma\)
on its source is
\begin{equation}
 k(t_v^n\otimes\Psi)
  =t_v^{n+\ord_v k}\otimes\lambda_{k^{(v)}}\Psi.
 \label{eq:ana:schwartz-cocycle}
\end{equation}
\end{proposition}
\begin{proof}
A \(U_v\)-invariant Bruhat function on \(K_v\) is constant on each
valuation shell. Its shell values are zero for sufficiently negative
valuations and are constant for sufficiently positive valuations.
Taking their successive differences expresses it uniquely as a
finite linear combination of the balls
\(1_{\varpi_v^n\mathcal O_v}\). On every fixed support and local-
constancy stage this is an isomorphism of finite-dimensional spaces.
Passing to the inductive bornology and tensoring with the other
places proves \eqref{eq:ana:ball-module}. Dilation carries the
\(n\)-th ball to the \((n+1)\)-st. Finally,
\(k_v\varpi_v^n\mathcal O_v
 =\varpi_v^{n+\ord_v k}\mathcal O_v\), which proves
\eqref{eq:ana:schwartz-cocycle}. This is the ball basis also used in
\cite[Lemma 4.18]{Meyer}.
\end{proof}

The normal and Schwartz transitions in
\eqref{eq:ana:frame-cocycle}--\eqref{eq:ana:schwartz-cocycle} agree;
for a finite set \(S\), their exponent vector is \((\ord_v k)_{v\in S}\).

For derived coinvariants, use the exact category of smooth complete
convex bornological representations in which a short exact sequence
is admissible when it splits by a bounded linear map after forgetting
the group action \cite[\S4.3, pp.~29--30]{Meyer}.
Let \(H_\Gamma\) denote total derived coinvariants in this exact
category. In particular,
\[
 H_\Gamma(V)=\C\widehat\otimes^{\mathbf L}_{\cD(\Gamma)}V,
 \qquad \cD(\Gamma)=\C[\Gamma],
\]
where the group algebra has the fine bornology. To regard these
complexes in the strict derived category of complete convex
bornological vector spaces, apply the natural exact localization
from the linearly split derived category. We write
\(H_+=\cS(\A_K)/\Gamma\) for its degree-zero coinvariant space.
Meyer's coinvariant theorem states that
\begin{equation}
 H_\Gamma\cS(\A_K)\simeq H_+
 \label{eq:ana:meyer-coinvariants}
\end{equation}
in the derived category; see \cite[Theorem 5.1 and the paragraph
following its proof]{Meyer}.

\begin{proposition}[Coefficient realization]
\label{prop:ana:projection}
For \(E\in\Perf_{\Sph[\Z^S]}\), set
\(E_\C=H\C\otimes_{\Sph}E\), regarded as a perfect complex over
\(R_{S,\C}=\C[\Z^S]\). Let \(U_S=\prod_{v\in S}U_v\). There is a
natural equivalence
\begin{equation}
 H_\Gamma\!\left(
 E_\C\widehat\otimes^{\mathbf L}_{R_{S,\C}}
             \cS(\A_K)^{U_S}\right)
 \simeq
 E_\C\widehat\otimes^{\mathbf L}_{R_{S,\C}}H_+^{U_S}.
 \label{eq:ana:projection}
\end{equation}
For a finite place \(v\), evaluation at zero in the \(v\)-coordinate
gives a strictly exact sequence
\begin{equation}
 0\longrightarrow
 \cD(K_v^\times)\widehat\otimes\cS(\A_K^{(v)})
 \longrightarrow\cS(\A_K)
 \xrightarrow{\mathrm{ev}_{0,v}}\cS(\A_K^{(v)})
 \longrightarrow0,
 \label{eq:ana:evaluation}
\end{equation}
where \(\cD(K_v^\times)\) denotes locally constant compactly
supported functions with its inductive bornology. This sequence
gives a distinguished triangle after applying \(H_\Gamma\).
\end{proposition}
\begin{proof}
Averaging over the compact group \(U_S\) is a bounded idempotent,
commutes with the \(\Gamma\)-action, and is exact. Applying this
idempotent to \eqref{eq:ana:meyer-coinvariants} proves
\eqref{eq:ana:projection} for \(E_\C=R_{S,\C}\).
Both sides are exact functors of \(E_\C\) and preserve finite sums
and retracts. The formula follows for all perfect complexes. The
naturality includes every Laurent matrix and every source map.
For \eqref{eq:ana:evaluation}, a locally constant function vanishing
at zero vanishes on a neighborhood of zero, so its restriction to
\(K_v^\times\) has compact support there. Extension by zero gives
the inverse identification of the kernel. The map
\(\Psi\mapsto1_{\mathcal O_v}\otimes\Psi\) is a bounded linear
section of evaluation. Thus the sequence is linearly split and
strictly exact; its derived coinvariants form the asserted triangle.
\end{proof}

Complex coefficient realization sends a primitive contact to zero:
for \(m\ge1\), a prime \(p\), and \(a\ge1\),
\[
 H\C\otimes_{\Sph}T_{m,mp^a}\simeq0,
\]
because its underlying spectrum is a wedge of \(\Sph/p\) and
\(p\) is invertible in \(\C\).
Its arithmetic endomorphism coefficients are retained in the spherical
correspondence construction through Theorem~\ref{thm:cor:descent}.
The analytic operators in \eqref{eq:ana:ball-module} use the Laurent
source action before this contact specialization.

\section{The Fourier--Poisson complex}
\label{sec:ana:poisson}

Fix the standard character \(\psi:\A_K/K\to\C^\times\) obtained
from the trace to \(\A_\Q/\Q\), and use the self-dual additive Haar
measure. Write
\[
 \widehat\Phi(y)=\int_{\A_K}\Phi(x)\psi(xy)\,dx.
\]
Fourier inversion gives \(\widehat{\widehat\Phi}(x)=\Phi(-x)\).
The involution \(x\mapsto-x\) belongs to the diagonal action of
\(\Gamma\), so the induced Fourier map on \(H_+\) is an involution.

\begin{definition}[Weighted Schwartz spaces]
\label{def:ana:weighted}
For an interval \(I\subset\R\), define
\[
 \cS(C)_I
   =\{f:C\to\C:|\cdot|^\sigma f\in\cS(C)
                         \text{ for every }\sigma\in I\}.
\]
The bornology is the initial bornology for the multiplication maps
\(f\mapsto|\cdot|^\sigma f\) to \(\cS(C)\). Here \(\cS(C)\) is
the Bruhat--Schwartz space, constructed from Schwartz functions on
Lie quotients of \(C\). Put
\[
 H_-=\cS(C)_\R,
 \qquad
 \mathcal W=\cS(C)_{(1,\infty)}\oplus
                   \cS(C)_{(-\infty,0)}.
\]
For a function on \(C\), define
\(Jf(c)=|c|^{-1}f(c^{-1})\).
\end{definition}

This is the weighted-space convention of \cite[Definition 4.1 and
\S5.3]{Meyer}. The summation map is
\[
 \Sigma\Phi(c)=\sum_{a\in K^\times}\Phi(ax),
 \qquad c=[x]\in C.
\]
It descends to a bounded map on \(H_+\).
The maps
\begin{equation}
 i_+:H_+\longrightarrow\mathcal W,
       \quad [\Phi]\longmapsto(\Sigma\Phi,J\Sigma\widehat\Phi),
 \qquad
 i_-:H_-\longrightarrow\mathcal W,
       \quad f\longmapsto(f,f)
 \label{eq:ana:embeddings}
\end{equation}
are bornological embeddings \cite[Lemmas 5.3--5.4 and the definition
of \(i_-\) in \S5.3, p.~41]{Meyer}.

\begin{lemma}[Fourier covariance]
\label{lem:ana:fourier-twist}
For \(a\in\A_K^\times\),
\begin{equation}
 \mathcal F\lambda_a=|a|\lambda_{a^{-1}}\mathcal F,
 \qquad J\lambda_a=|a|\lambda_{a^{-1}}J.
 \label{eq:ana:fourier-twist}
\end{equation}
Both maps in \eqref{eq:ana:embeddings} are \(C\)-equivariant.
For a finite set \(S\), let
\[
 \iota_S:R_{S,\C}\to R_{S,\C},
 \qquad \iota_S(t_v)=q_v^{-1}t_v^{-1}.
\]
Let \(M\) be either \(\cS(\A_K)^{U_S}\) or \(H_+^{U_S}\).
The raw Fourier map is \(R_{S,\C}\)-linear from \(M\) to
\(\iota_S^*M\), where \(\iota_S^*\) denotes restriction of scalars.
For every perfect \(R_{S,\C}\)-complex \(E_\C\), its target after
coefficient substitution is naturally
\begin{equation}
 E_\C\otimes^{\mathbf L}_{R_{S,\C}}\iota_S^*M
 \simeq
 \iota_S^*\!\left(
       \iota_S^*E_\C\otimes^{\mathbf L}_{R_{S,\C}}M\right).
 \label{eq:ana:outer-twist}
\end{equation}
\end{lemma}
\begin{proof}
Substitute \(x=ay\) in the Fourier integral. The Jacobian is
\(|a|\), and the Fourier argument becomes \(ay\), which proves the
first identity. For the second, direct evaluation gives
\[
 J\lambda_af(x)=|x|^{-1}f(a^{-1}x^{-1})
              =|a|\lambda_{a^{-1}}Jf(x).
\]
The product formula lets Fourier descend through \(\Gamma\)-coinvariants:
it carries the relation \((\lambda_k-1)\Phi\) to
\((\lambda_{k^{-1}}-1)\widehat\Phi\). Summation is equivariant.
The two factors \(|a|\) and \(|a^{-1}|\) cancel in
\(J\Sigma\mathcal F\), proving equivariance of \(i_+\).
Since \(|h_v|=q_v^{-1}\), the first identity induces \(\iota_S\).
The relation on the left of \eqref{eq:ana:outer-twist} is
\(re\otimes m=e\otimes\iota_S(r)m\).
On the right it is
\(\iota_S(s)e\otimes m=e\otimes sm\); set \(s=\iota_S(r)\).
The outer restriction makes this identification linear for the
original \(R_{S,\C}\)-action. Applying the same argument to a
finite-cell resolution and then to retracts gives the derived formula.
\end{proof}

Define the two equivariant maps
\[
 \epsilon_0:H_+\to\C_1,\quad [\Phi]\mapsto\Phi(0),
 \qquad
 \epsilon_1:H_+\to\C_{|\cdot|},\quad [\Phi]\mapsto\int_{\A_K}\Phi,
\]
where \(\C_\chi\) denotes the one-dimensional representation of
character \(\chi\). Both maps are well defined by the product formula.
Their joint kernel will be denoted by \(H_\cap\).

\begin{proposition}[The common Poisson subspace]
\label{prop:ana:intersection}
The map \((\epsilon_0,\epsilon_1)\) is surjective. The summation map
restricts to an equivariant embedding
\(\theta:H_\cap\to H_-\), and
\[
 i_+(H_\cap)=i_+(H_+)\cap i_-(H_-)
\]
inside \(\mathcal W\). Its range is closed in the induced
bornological quotient construction.
\end{proposition}
\begin{proof}
Choose a Bruhat--Schwartz function supported away from zero with
nonzero integral, and another function equal to \(1\) at zero.
Subtracting a multiple of the first function from the second gives
joint surjectivity. Poisson summation on \(\A_K/K\) gives
\begin{equation}
 \Sigma\Phi(c)-J\Sigma\widehat\Phi(c)
       =|c|^{-1}\epsilon_1[\Phi]-\epsilon_0[\Phi].
 \label{eq:ana:poisson}
\end{equation}
Thus the two coordinates agree if and only if the two moments vanish.
The intersection of the two weighted Schwartz spaces is \(H_-\),
which identifies the common image. Fix an open subgroup \(U\) of the
finite integral units, and choose \(S\) sufficiently large that
\(U^S=\prod_{v\notin S}\mathcal O_v^\times\subset U\).
On the corresponding invariant stage the groups have Lie quotients.
The description in \cite[\S4.2, equation~(24), and Theorem~5.1]{Meyer},
followed by taking \(U/U^S\)-invariants, identifies \(H_+^U\)
with a complete Fr\'echet space.
The bornological embedding of
\cite[Lemmas~4.20, 5.3 and 5.4]{Meyer} restricts to a topological
embedding \(i_+:H_+^U\to\mathcal W^U\).
Indeed, its inverse on its image is bounded for the induced
bornology, hence continuous on this metrizable locally convex
subspace. Its image is therefore complete for the induced topology
and hence closed in \(\mathcal W^U\).
Equivariance and injectivity give
\(i_+(H_+^U)=i_+(H_+)\cap\mathcal W^U\).
By \eqref{eq:ana:poisson}, the inverse image of this closed subspace
under the continuous map \(i_-:H_-^U\to\mathcal W^U\) is
\(\theta(H_\cap^U)\). Thus \(\theta(H_\cap^U)\) is closed in \(H_-^U\),
and \(H_-^U/\theta(H_\cap^U)\) is Fr\'echet.
The compact-subgroup averaging operators commute with both
embeddings and split the inclusions between the quotient stages.
Passing through the resulting inductive system gives the asserted
bornological quotient.
\end{proof}

\begin{definition}[Fourier--Poisson complex]
\label{def:ana:FP}
Define the complex of smooth \(C\)-representations
\begin{equation}
 C_K^{\mathrm{FP}}=
 \left[
  \C_1\oplus H_\cap
       \xrightarrow{\ (0,\theta)\ }H_-
       \xrightarrow{\ 0\ }\C_{|\cdot|}
 \right]
 \label{eq:ana:FP}
\end{equation}
in cohomological degrees \(0,1,2\), with all other terms zero.
Its cohomology is
\[
 H^0=\C_1,\qquad
 H^1=H_-/\theta H_\cap,\qquad
 H^2=\C_{|\cdot|},
\]
where the middle term carries the quotient bornology of
Proposition~\ref{prop:ana:intersection}.
\end{definition}

The action of \(C\) on every term commutes with the differential.
Integration of this action gives an action of \(\mathscr A_K\).
The map \(\mathscr A_K^{\mathrm{sph}}\to H\mathscr A_K\)
therefore defines an exact realization
\[
 \mathcal R_K:\Perf(\mathscr A_K^{\mathrm{sph}})
                    \longrightarrow\mathcal D_{\mathrm{bor}}(\C),
 \qquad \mathscr A_K^{\mathrm{sph}}\longmapsto C_K^{\mathrm{FP}}.
\]
Here \(\mathcal D_{\mathrm{bor}}(\C)\) is the derived category of
complete convex bornological complex vector spaces with the strict
exact structure. The functor is obtained on finite free modules by
finite direct sums and matrices of the displayed operators, and on
perfect modules by finite cones and idempotent splittings.
Combined with \eqref{eq:cor:source}, it realizes source cospans as
cospans of complexes. The assembled equivariant diagram
\eqref{eq:ana:embeddings} admits the same coefficient substitution
as \eqref{eq:ana:projection}.

\section{The trace and the normal fixed strata}
\label{sec:ana:trace}

\begin{definition}[Trace distributions at a place]
\label{def:ana:weil}
Give \(K_v^\times\) the Weil Haar measure \(d_W^\times a\), normalized
by
\[
 \int_{1\le |a|_v<|b|_v}d_W^\times a=\log|b|_v
             \qquad (|b|_v>1).
\]
At a finite place its unit mass is \(\log q_v\). For
\(\varphi\in\cD(K_v^\times)\), define
\begin{equation}
 W_v^0(\varphi)
  =\int_{K_v^\times}
   \bigl(\varphi(a)-\varphi(1)1_{\mathcal O_v^\times}(a)\bigr)
        \frac{|a|_v}{|1-a|_v}\,d_W^\times a
       \quad(v\nmid\infty).
 \label{eq:ana:finite-PV}
\end{equation}
Choose a reference additive character \(\psi_{v,0}\) of conductor
\(\mathcal O_v\) at a finite place. At the real and complex places,
respectively, take
\[
 \psi_{v,0}(x)=e^{2\pi i x},\qquad
 \psi_{v,0}(z)=e^{4\pi i\Re z}.
\]
The associated self-dual measures \(dx_{v,0}\) give volume \(1\)
to \(\mathcal O_v\) at a finite place and are \(dx\) and
\(2\,dx\,dy\) at the infinite places. With
\(\mathcal F_{v,0}\) the resulting Fourier transform, define the
distribution \(P_v\) by
\[
 \langle P_v,\eta\rangle
   =\int_{K_v}\log|x|_v\,
              (\mathcal F_{v,0}\eta)(x)\,dx_{v,0}
       \qquad(\eta\in\cS(K_v)).
\]
At an infinite place put
\begin{equation}
 W_v^0(\varphi)=-\langle P_v,x\mapsto\varphi(1-x)\rangle.
 \label{eq:ana:infinite-PV}
\end{equation}
In \eqref{eq:ana:infinite-PV}, \(\varphi\) is extended by zero at
\(0\). The distribution \(P_v\) restricts to \(-d_W^\times x\)
away from zero. These are the local Fourier trace distributions of
\cite[Lemmas 3.4--3.5 and equation (19)]{Meyer}, expressed using
the reference characters above.
Define
\begin{equation}
 W_K(f)=\sum_v W_v^0(f\circ j_v)-\log|D_K|\,f(1),
 \qquad f\in\cD(C).
 \label{eq:ana:Weil}
\end{equation}
\end{definition}

At a finite place the subtraction in
\eqref{eq:ana:finite-PV} makes the integrand vanish near \(1\),
and \(W_v^0(1_{\mathcal O_v^\times})=0\). The Fourier normalization
gives the same value since \(P_v\) takes the value \(-\log q_v/(q_v-1)\)
on both \(1_{\mathcal O_v}\) and \(1_{1+\varpi_v\mathcal O_v}\),
using \(\mathcal F_{v,0}1_{\mathcal O_v}=1_{\mathcal O_v}\) and the
restriction of \(P_v\) away from zero, respectively.
For \(f\in\cD(C)\), only finitely many local terms in
\eqref{eq:ana:Weil} are nonzero. Indeed, choose a finite set of places
outside which \(f\) is unit-invariant, and a compact interval
containing the norms of its support. Outside a larger finite set,
only the unit part can meet that interval; there \(f\circ j_v\) is
\(f(1)1_{\mathcal O_v^\times}\), whose local distribution is zero.

The discriminant in \eqref{eq:ana:Weil} is the change from these
reference characters to the global trace character. Write
\(\psi_v(x)=\psi_{v,0}(\delta_vx)\). The self-dual measure for \(\psi_v\)
is \(dx_{v,\psi}=|\delta_v|_v^{1/2}dx_{v,0}\), and
\[
 \mathcal F_{v,\psi}\eta(x)
    =|\delta_v|_v^{1/2}(\mathcal F_{v,0}\eta)(\delta_vx).
\]
Define \(P_{v,\psi}\) by the same distribution pairing with
\(\mathcal F_{v,\psi}\) and \(dx_{v,\psi}\). Substitution of
\(y=\delta_vx\), followed by Fourier inversion at zero, gives
\[
 P_{v,\psi}=P_v-\log|\delta_v|_v\,\delta_0.
\]
The local trace term therefore changes by
\(\log|\delta_v|_v\varphi(1)\). At the infinite places \(\delta_v=1\);
at a finite place its valuation is the different exponent.
Consequently \(\prod_v|\delta_v|_v^{-1}=|D_K|\), giving the
discriminant term in \eqref{eq:ana:Weil}. This is the
change-of-character term in \cite[equation (20), \S5.1,
and Definition 5.7]{Meyer}.

\begin{theorem}[Complete trace]
\label{thm:ana:complete-trace}
For every \(f\in\cD(C)\), the operators \(\rho(f)\) on the three
cohomology spaces of \(C_K^{\mathrm{FP}}\) are summable, and
\begin{equation}
 \operatorname{str}_{H^*(C_K^{\mathrm{FP}})}\rho(f)
 =\int_C f(c)\,dc+\int_Cf(c)|c|\,dc
      -\Tr\!\left(\rho(f)\mid H^1\right)
 =W_K(f).
 \label{eq:ana:complete-trace}
\end{equation}
For a finite direct sum of copies of \(C_K^{\mathrm{FP}}\), the
supertrace of a matrix \((\rho(f_{ij}))\), with
\(f_{ij}\in\cD(C)\), is \(\sum_iW_K(f_{ii})\).
\end{theorem}
\begin{proof}
Proposition~\ref{prop:ana:intersection} identifies
\(H_+/(H_+\cap H_-)\) with \(\C_1\oplus\C_{|\cdot|}\), and
\(H_-/(H_+\cap H_-)\) with \(H^1\). These are the positive and
negative representations in \cite[\S5.5]{Meyer}.
Their summability and character are
\cite[Theorem 5.8]{Meyer}, with the local character normalization
computed above. Since the two positive characters occupy
even degrees in \eqref{eq:ana:FP}, its character is precisely the
supertrace in \eqref{eq:ana:complete-trace}. The two one-dimensional
traces are the displayed moments. Finite additivity gives the
matrix formula.
\end{proof}

\begin{theorem}[Fixed strata and their normal contribution]
\label{thm:ana:fixed}
Let \(a\in C\setminus\{1\}\). If
\(a\in j_v(K_v^\times)\), then
\(\operatorname{Fix}(\gamma_a)=Z_v\); otherwise the fixed locus is
empty. At most one place \(v\) has this property. For
\(a=j_v(a_v)\), the normal action at \(Z_v\) is multiplication by
\(a_v\), and the inverse-pullback normal factor is
\begin{equation}
 \frac1{|1-a_v^{-1}|_v}=\frac{|a_v|_v}{|1-a_v|_v}.
 \label{eq:ana:normal-weight}
\end{equation}
If \(f\in\cD(C)\) is supported away from \(1\), then
\begin{equation}
 \sum_v\int_{K_v^\times}
       f(j_v(a))\frac{d_W^\times a}{|1-a^{-1}|_v}
  =\operatorname{str}_{H^*(C_K^{\mathrm{FP}})}\rho(f).
 \label{eq:ana:fixed-trace}
\end{equation}
\end{theorem}
\begin{proof}
The action on the common open \(C\) is free. On
\(Z_v=C/j_v(K_v^\times)\), an element fixes a point if and only if
it belongs to \(j_v(K_v^\times)\); since \(C\) is abelian, it then
fixes every point. If an element is represented by both a single-
place idele at \(v\) and one at \(w\ne v\), their quotient is a
diagonal element of \(K^\times\). At a third place this element has
coordinate \(1\), so it is \(1\) in \(K\). Hence the two local
subgroups intersect only at \(1\).
Balancing in \(N_v\) gives
\([cj_v(a_v),z]=[c,a_vz]\). The normal difference for inverse
pullback is multiplication by \(1-a_v^{-1}\), whose additive
Jacobian is \(|1-a_v^{-1}|_v\). This proves
\eqref{eq:ana:normal-weight}, including the squared absolute value
at a complex place.

The map \(j_v\) is proper, so the pullback of the support of \(f\)
is compact and avoids \(1\). All integrals in
\eqref{eq:ana:fixed-trace} are ordinary integrals. Moreover, for
all but finitely many \(v\), the subgroup \(j_v(U_v)\) lies in an
identity neighborhood disjoint from the support. Nonunit elements
have norms \(q_v^n\), \(n\ne0\), and cannot meet its fixed compact
norm interval when \(q_v\) is sufficiently large. The sum is thus
finite. Each \(Z_v\) has quotient Haar mass \(1\) for the chosen
Haar measure on \(C\) and the Weil Haar measure on its stabilizer:
at a finite place the norm interval has length \(\log q_v\), equal
to the local unit mass; at an infinite place both measures push
forward to \(d\log r\). Thus the displayed integrals are the
normal contributions of the fixed strata. Since \(f(1)=0\),
Definition~\ref{def:ana:weil} and Theorem~\ref{thm:ana:complete-trace}
give \eqref{eq:ana:fixed-trace}.
\end{proof}

\begin{definition}[Source operators and transpose]
\label{def:ana:source-operators}
For \(S\) finite, write \(h^{\mathbf n}=\prod_{v\in S}h_v^{n_v}\)
and \(T_v=\lambda_{h_v}\). If
\(P=\sum_{\mathbf n}a_{\mathbf n}t^{\mathbf n}\in R_{S,\C}\), set
\[
 f_P(c)=\sum_{\mathbf n}a_{\mathbf n}
                       f((h^{\mathbf n})^{-1}c).
\]
The conjugate-linear involution on \(\mathscr A_K\) is specified by
\begin{equation}
 \delta_a^\sharp=|a|\delta_{a^{-1}},\qquad
 f^\sharp(c)=|c|^{-1}\overline{f(c^{-1})}.
 \label{eq:ana:sharp}
\end{equation}
The induced polynomial involution is
\[
 P^\sharp(t)=\sum_{\mathbf n}\overline{a_{\mathbf n}}
       \left(\prod_{v\in S}q_v^{-n_v}\right)t^{-\mathbf n}.
\]
\end{definition}

\begin{proposition}[Composition with the trace]
\label{prop:ana:source-trace}
For \(P,Q\in R_{S,\C}\) and \(f,g\in\cD(C)\),
\begin{align}
 P(T)\rho(f)&=\rho(f_P),\label{eq:ana:source-action}\\
 \operatorname{str}\!\left(
       P(T)\rho(f)Q^\sharp(T)\rho(g^\sharp)\right)
       &=W_K\!\left(f_P*(g_Q)^\sharp\right).
       \label{eq:ana:source-pairing}
\end{align}
If \(\int_C f(c)\,dc=\int_C f(c)|c|\,dc=0\), and the same holds
for \(g\), then the right side of
\eqref{eq:ana:source-pairing} is
\(-\Tr(\rho(f_P*(g_Q)^\sharp)\mid H^1)\).
\end{proposition}
\begin{proof}
Translation of the integration variable gives
\(\lambda_h\rho(f)=\rho(c\mapsto f(h^{-1}c))\), hence
\eqref{eq:ana:source-action}. The rule \(\sharp\) is an involution
and reverses convolution. For example, on the test-function pairing
\(\int_Cu(c)\overline{v(c)}|c|\,dc\), change of variables gives
\(\lambda_a^*=|a|\lambda_{a^{-1}}\); integrating this identity gives
\eqref{eq:ana:sharp}. Consequently
\((g_Q)^\sharp=Q^\sharp(\delta_h)*g^\sharp\).
All finite translates and convolutions remain in \(\cD(C)\), so
Theorem~\ref{thm:ana:complete-trace} applies to their product and
proves \eqref{eq:ana:source-pairing}.
The functionals \(f\mapsto\int f|\cdot|^j\), \(j=0,1\), are
convolution characters. Their common kernel is an ideal preserved
by translation, and \(\sharp\) exchanges the two characters with
complex conjugation. The two even-degree traces therefore vanish
under the stated moment conditions.
\end{proof}

\section{Normal determinants and the infinite factors}
\label{sec:ana:determinants}

\begin{definition}[Normal determinant factor]
\label{def:ana:determinant}
For a finite place \(v\) and \(\Re(s)>0\), specialize the normal
puncture differential of \eqref{eq:nor:puncture} by
\(t_v\mapsto q_v^{-s}\). Use the shell basis in the source and
the ball basis in the target, and define
\[
 \tau_v(s)=\det\bigl(1-q_v^{-s}:\C\to\C\bigr)^{-1}
          =(1-q_v^{-s})^{-1}.
\]
This is the inverse determinant of the specialized two-term
acyclic complex.
\end{definition}

\begin{proposition}
\label{prop:ana:determinant-laws}
The factor \(\tau_v(s)\) is independent of a change of the normal
covering frame. For positive integers \(L,h\), the period
refinement is compatible with
\[
 1-t^{hL}=(1-t^L)\sum_{j=0}^{h-1}t^{jL}.
\]
After complex specialization, the determinant lines of finite source
complexes are multiplicative in their sewing triangles.
\end{proposition}
\begin{proof}
A frame change multiplies both bases by the same Laurent unit, which
cancels in the determinant. The polynomial refinement identity gives
multiplicativity wherever its two factors are invertible. Finally,
\(\det C^\bullet=\bigotimes_i(\det C^i)^{(-1)^i}\) sends the
mapping cone of a sewing map to the target determinant line tensored
with the inverse source line.
\end{proof}

The self-dual local additive measures for the standard trace character
are \(dx\) at a real place, \(2\,dx\,dy\) at a complex place, and
have \(\mathcal O_v\)-volume \(q_v^{-d_v/2}\) at a finite place,
where \(d_v\) is the different exponent
\cite[Proposition 4.8]{Poonen}.
For the following Mellin integrals, use multiplicative measure
\(dx/|x|\) at a real place,
\(2\,dx\,dy/(\pi|z|^2)\) at a complex place, and unit mass \(1\)
at a finite place.

\begin{proposition}[The archimedean normal integrals]
\label{prop:ana:gaussians}
For \(\Re(s)>0\),
\begin{align*}
 \int_{\R^\times}e^{-\pi x^2}|x|^s\frac{dx}{|x|}
     &=\pi^{-s/2}\Gamma(s/2)=:\Gamma_\R(s),\\
 \int_{\C^\times}e^{-2\pi|z|^2}|z|^{2s}
                  \frac{2\,dx\,dy}{\pi|z|^2}
     &=2(2\pi)^{-s}\Gamma(s)=:\Gamma_\C(s).
\end{align*}
Under inverse-pullback dilation, the even real Taylor coefficient
of order \(2n\) has weight \(|a|^{-2n}\), and the complex radial
coefficient of \((z\bar z)^n\) has weight \(|a|_v^{-n}\).
\end{proposition}
\begin{proof}
The real integral is twice the positive half-line integral. Substitute
\(u=\pi x^2\) to obtain the first formula. In polar coordinates
the complex integral is
\(4\int_0^\infty e^{-2\pi r^2}r^{2s-1}\,dr\).
The substitution \(u=2\pi r^2\) gives the second formula.
The Taylor statements follow by substituting \(a^{-1}x\) and
\(a^{-1}z\) in the indicated monomials, using
\(|a|_v=|a|_{\mathrm{ordinary}}^2\) at a complex place.
\end{proof}

\begin{proposition}[The product of normal factors]
\label{def:ana:completed-product}
Let \(r_1,r_2\) be the numbers of real and complex places of \(K\).
Put \(V_K=\operatorname{vol}_{\psi_f}(\widehat{\mathcal O}_K)\),
using the finite self-dual additive measure. Define
\[
 \mathcal Z_K(s)=V_K^{-s}
     \Gamma_\R(s)^{r_1}\Gamma_\C(s)^{r_2}
     \prod_{v\nmid\infty}\tau_v(s),
 \qquad \Re(s)>1.
\]
\label{prop:ana:completed-product}
The product converges locally uniformly on \(\Re(s)>1\), and
\[
 \mathcal Z_K(s)
   =|D_K|^{s/2}\Gamma_\R(s)^{r_1}\Gamma_\C(s)^{r_2}\zeta_K(s).
\]
\end{proposition}
\begin{proof}
The local self-dual volumes give
\(V_K=\prod_{v\nmid\infty}q_v^{-d_v/2}=|D_K|^{-1/2}\).
For \(\Re(s)\ge1+\varepsilon\), one has
\(\sum_vq_v^{-\Re(s)}\le[K:\Q]\sum_pp^{-1-\varepsilon}<\infty\).
This proves local uniform convergence of the product. Unique
factorization of nonzero integral ideals gives
\(\prod_v(1-q_v^{-s})^{-1}=\zeta_K(s)\).
Substitution proves the formula, with the gamma normalizations of
\cite[\S4.7]{Poonen}.
\end{proof}
\section{Multiplication and integral sections}
\label{sec:products}

The valuation torsor gives graded coefficient lines with a section product.
The marked Laurent coordinate then gives divisor section lattices with
integral jet conditions.

\subsection{Graded descent along a valuation torsor}

Let $Q_v\to Y_K$ be the principal $\mathbb Z$-cover associated with the
finite place $v$. Put $A=\mathbb Z[t,t^{-1}]$ and let $\underline A$ denote
the constant sheaf of rings on $Y_K$.

\begin{definition}[The winding coefficient line]
\label{def:prod-line}
The invertible $\underline A$-module $\mathscr L_v$ is the associated sheaf
\[
 \mathscr L_v=\mathbb Z[Q_v]\otimes_{\mathbb Z[\mathbb Z]} A,
 \qquad n\cdot f=t^n f.
\]
For $m\ge0$ put $\mathscr L_v^{[m]}=\mathscr L_v^{\otimes_{\underline A}m}$,
with $\mathscr L_v^{[0]}=\underline A$, and define
\[
 \mathscr S_v=\bigoplus_{m\ge0}\mathscr L_v^{[m]}.
\]
\end{definition}

\begin{proposition}[Multiplicative descent]
\label{prop:prod-descent}
The sheaf $\mathscr S_v$ is a graded commutative $\underline A$-algebra.
Its multiplication identifies
\[
 \mathscr L_v^{[m]}\otimes_{\underline A}\mathscr L_v^{[n]}
       \simeq\mathscr L_v^{[m+n]}.
\]
The degree-one sheaf, after applying the Eilenberg--Mac Lane functor, is
the integral realization of the free winding coefficient along $Q_v$.
These constructions commute with pullback of the torsor and with extension
of the coefficient ring $A$.
\end{proposition}

\begin{proof}
For torsor frames $q_j=q_i+n_{ij}$, the line frames satisfy
$u_j=t^{n_{ij}}u_i$. The algebra transitions $A[u_j]\to A[u_i]$
therefore act in degree $m$ by $t^{mn_{ij}}$. The torsor cocycle gives
descent, and polynomial multiplication gives the stated tensor product.
Degree one has the free winding-module transitions. Pullback and coefficient
extension preserve these equations.
\end{proof}

\begin{definition}[Sections at finite period]
\label{def:prod-period}
For $L\ge1$ put
\[
 A_L=A/(t^L-1)=\mathbb Z[C_L],\qquad
 \mathscr L_{v,L}^{[m]}=\mathscr L_v^{[m]}\otimes_A A_L,
 \qquad
 \mathscr S_{v,L}=\bigoplus_{m\ge0}\mathscr L_{v,L}^{[m]}.
\]
A real coefficient fibre has the cone of coefficient-positive vectors
$f=\sum_{j=0}^{L-1}a_jt^j$, defined by $a_j\ge0$.
\end{definition}

\begin{proposition}[Effectivity and a common norm at finite period]
\label{prop:prod-finite-positive}
The coefficient-positive subsets are independent of the chosen frames.
They are closed under addition and graded multiplication, and their
intersection with their negatives is zero. On a complex coefficient fibre
set
\[
 \|f\|_{L,m}^2=e^{-2m\eta}L\sum_{j=0}^{L-1}|a_j|^2,
 \qquad \eta\in\mathbb R.
\]
This gives a positive Hermitian metric on each degree, preserved by the
deck transitions. Multiplication satisfies
\[
 \|fg\|_{L,m+n}\le\|f\|_{L,m}\|g\|_{L,n}.
\]
A multiplier $f\in A_L$ acts injectively on $A_L$ if and only if
$f(\zeta)\ne0$ for every $L$th root of unity $\zeta$.
\end{proposition}

\begin{proof}
Frame changes permute coefficients, preserving the cone and metric.
The cone is proper and closed under addition and cyclic convolution.
Young's inequality and Cauchy--Schwarz give
\[
 \|fg\|_{\ell^2}\le\|f\|_{\ell^1}\|g\|_{\ell^2}
       \le\sqrt L\,\|f\|_{\ell^2}\|g\|_{\ell^2}.
\]
The degree factors give the asserted estimate. Since $A_L$ is
$\mathbb Z$-free, injectivity is equivalent to a nonzero determinant;
over $\mathbb C$ that determinant is $\prod_{\zeta^L=1}f(\zeta)$.
\end{proof}

\begin{proposition}[Period changes and the metric]
\label{prop:prod-period-metric}
Let $M=hL$ with $h\ge1$. The quotient homomorphism $\pi:A_M\to A_L$ is a contraction
for the norms in Proposition~\ref{prop:prod-finite-positive}. Define
\[
 \operatorname{tr}_{L}^{M}(f)
    =\left(\sum_{j=0}^{h-1}t^{jL}\right)\widetilde f,
 \qquad
 \iota_L^M(f)=h^{-1}\operatorname{tr}_{L}^{M}(f),
\]
where $\widetilde f$ is any Laurent lift. The integral transfer has norm
factor $h$. Over $\mathbb Q$ or $\mathbb C$, $\iota_L^M$ is an isometric
multiplicative embedding, with
\[
 \pi\iota_L^M=1,
 \qquad
 \iota_L^M(1)=h^{-1}\sum_{j=0}^{h-1}t^{jL}.
\]
All three maps preserve coefficient positivity and commute with the graded deck
transitions. Quotients, integral transfers and normalized transfers compose
in towers of periods.
\end{proposition}

\begin{proof}
Cauchy--Schwarz on each quotient fibre gives
\[
 L\sum_{r<L}\left|\sum_{j<h}a_{r+jL}\right|^2
       \le Lh\sum_{k<M}|a_k|^2=M\sum_{k<M}|a_k|^2.
\]
The transfer repeats coefficients $h$ times, giving squared norm
$h^2L\sum_{r<L}|a_r|^2$. For $N=\sum_{j<h}t^{jL}$ one has
$N^2=hN$ and $N(t^L-1)=0$ in $A_M$. These identities give
independence of lifts, multiplicativity of $N/h$, and $\pi\iota_L^M=1$.
All maps commute with powers of $t$; grouping quotient fibres proves
the tower identities.
\end{proof}

\subsection{Cohomology of a boundary circle}

Put $B_v=\mathbb R_{>0}/q_v^{\mathbb Z}$. This circle is the image of
$Z_v$ in $C_K^1\backslash Y_K$. Quotienting the valuation cover over
$Z_v$ by $C_K^1$ gives
\[
 \overline Q_v=\mathbb R_{>0}\longrightarrow B_v,
 \qquad r\cdot n=rq_v^{-n}.
\]
The two spaces before this quotient are
$C_K/j_v(\mathcal O_v^\times)$ and $C_K/j_v(K_v^\times)$.
Use $\overline Q_v$ in Definition~\ref{def:prod-line} to form
coefficient lines on $B_v$, orienting its generator with deck label one.

\begin{proposition}[Boundary section ranks]
\label{prop:prod-circle}
In degree $m\ge0$ and at finite period $L$, the cohomology is computed by
\[
 [A_L\xrightarrow{t^m-1}A_L]
\]
in degrees zero and one. Both cohomology groups are free abelian of rank
$\gcd(L,m)$, with $\gcd(L,0)=L$.
For the coefficient ring $A$ and $m\ge1$, the same calculation gives
\[
 H^0=0,\qquad H^1=A/(t^m-1).
\]
For $A/(\Phi_u)$, if $u\nmid m$ then $H^0=0$ and $H^1$ is finite of
order $|\operatorname{Res}(\Phi_u,t^m-1)|$. If $u\mid m$, both groups
are $A/(\Phi_u)$.
\end{proposition}

\begin{proof}
The cellular differential is monodromy minus identity, here $t^m-1$.
The permutation $j\mapsto j+m$ on $C_L$ has $\gcd(L,m)$ orbits;
its invariants and coinvariants are free on those orbits. On $A$ the
differential is injective. On $A/(\Phi_u)$ it is zero if $u\mid m$,
and otherwise is injective with determinant
$\operatorname{Res}(\Phi_u,t^m-1)$.
\end{proof}

\subsection{Divisors on the marked projective coordinate}

Use the projective coefficient realization $X=\mathbb P^1_{\mathbb Z}$
with marked coordinate $t$. Let $P_i$ be finitely many distinct monic
cyclotomic polynomials with $P_i(1)\ne0$, and let $Z(P_i)$ be their
zero divisors in $X$. The finite family may enlarge under divisor addition
and factorization of $P_i(t^h)$; the new factors also have nonzero value at $1$.

\begin{definition}[Marked divisors and their section lattices]
\label{def:prod-divisor}
For integers $a,b,r_i$ put
\[
 D=a[0]+b[\infty]+\sum_i r_iZ(P_i),\qquad
 F_D(t)=\prod_iP_i(t)^{r_i},\qquad
 n(D)=a+b+\sum_i r_i\deg P_i.
\]
Define
\[
 S(D)=H^0(X,\mathcal O_X(D)).
\]
Effectivity of a nonzero rational section means
$\operatorname{div}(f)+D\ge0$ as a Cartier divisor.
\end{definition}

\begin{proposition}[Explicit sections and their product]
\label{prop:prod-divisor-sections}
If $n=n(D)\ge0$, then
\[
 S(D)=t^{-a}F_D(t)^{-1}\mathbb Z[t]_{\le n}.
\]
If $n(D)<0$, this group is zero. Multiplication of rational functions
gives
\[
 S(D)\otimes_{\mathbb Z}S(E)\longrightarrow S(D+E),
\]
and multiplication by any nonzero section is injective. The numerator
identifications
\[
 \nu_D:S(D)\longrightarrow\mathbb Z[t]_{\le n(D)},
 \qquad \nu_D(f)=t^aF_D f
\]
satisfy $\nu_{D+E}(fg)=\nu_D(f)\nu_E(g)$.
\end{proposition}

\begin{proof}
The divisor of $t^aF_D$ is
\[
 a[0]+\sum_i r_iZ(P_i)
       -\left(a+\sum_i r_i\deg P_i\right)[\infty].
\]
Thus $t^aF_D$ identifies $\mathcal O_X(D)$ with
$\mathcal O_X(n[\infty])$, whose sections are integral polynomials of
degree at most $n$ and vanish for $n<0$. Multiplication adds divisors
and is injective for a nonzero factor in $\mathbb Q(t)$. The numerator
identity follows from $F_{D+E}=F_DF_E$.
\end{proof}

For a change of winding frame by $k\in\mathbb Z$, put
\[
 D'=D-k\operatorname{div}(t),\qquad f'=t^kf.
\]
Then $n(D')=n(D)$ and $\nu_{D'}(f')=\nu_D(f)$.
For a tensor-degree-$m$ transition, take $k=mn_{ij}$.

\subsection{An integral condition on the boundary jets}

For $M\ge1$ put
\[
 J_M=\mathbb Z[x]/(x^M),\qquad
 J_{M,\mathbb Q}=\mathbb Q[x]/(x^M).
\]
The jet $j_M(f)=f(1+x)\bmod x^M$ is defined on $S(D)$ because
$P_i(1)\ne0$.

\begin{definition}[Sections with integral jets]
\label{def:prod-jet-sections}
Set
\[
 S_M(D)=\{f\in S(D):j_M(f)\in J_M\},\qquad
 \iota_M(D)=\log[S(D):S_M(D)]\quad(n(D)\ge0).
\]
\end{definition}

\begin{theorem}[Multiplication with the jet condition]
\label{thm:prod-integral-jets}
For $n(D)\ge0$, $S_M(D)$ is a full-rank finite-index sublattice of
$S(D)$. Multiplication restricts to
\[
 S_M(D)\otimes S_M(E)\longrightarrow S_M(D+E),
\]
and multiplication by every nonzero section remains injective.
Jet truncation gives $S_{M+1}(D)\subseteq S_M(D)$, and the frame
change $(f,D)\mapsto(t^kf,D-k\operatorname{div}(t))$ preserves these
lattices. For integers $h\ge1$, power pullback $\phi_h(t)=t^h$ sends $S_M(D)$ to
$S_M(\phi_h^*D)$, compatibly with products and composition.
\end{theorem}

\begin{proof}
Clearing denominators in the rational jet matrix gives
$cS(D)\subseteq S_M(D)$ for some $c>0$. Since $j_M$ is a ring map,
integral jets are closed under products; injectivity follows in
$\mathbb Q(t)$. Truncation is integral, and $(1+x)^k$ is an integral
jet unit for every $k\in\mathbb Z$. Power pullback induces the ring map
$J_M\to J_M$, $x\mapsto(1+x)^h-1$, which is well defined because its
image lies in $(x)$. These substitutions commute with products and compose.
\end{proof}

\begin{theorem}[The section index as a lattice expansion]
\label{thm:prod-index-expansion}
Assume $n(D)=M-1$ and put
\[
 u_D=j_M(t^{-a}F_D^{-1})\in J_{M,\mathbb Q}^{\times}.
\]
For commensurable lattices write
$e_+(L\mid B)=\log\#(L/(L\cap B))$. Then
\[
 \boxed{\quad \iota_M(D)=e_+(u_DJ_M\mid J_M).\quad}
\]
In numerator jet coordinates the section lattice is
\[
 \nu_D(S_M(D))=J_M\cap u_D^{-1}J_M.
\]
The equality is understood through the integral identification
$p(t)\mapsto p(1+x)\bmod x^M$.
\end{theorem}

\begin{proof}
The binomial change of basis $p(t)\mapsto p(1+x)$ is integral
unimodular. The section jet is $u_Dj_M(p)$, so integrality gives the
stated intersection. Multiplication by $u_D$ identifies
\[
 J_M/(J_M\cap u_D^{-1}J_M)
       \simeq u_DJ_M/(u_DJ_M\cap J_M),
\]
which proves the formula.
\end{proof}

\begin{corollary}[Signed cyclotomic section indices]
\label{cor:prod-signed-index}
Put
\[
 D_0=2Z(\Phi_2)-3Z(\Phi_4)+Z(\Phi_8),\qquad
 R(t)=\frac{\Phi_2(t)^2\Phi_8(t)}{\Phi_4(t)^3}.
\]
For $m\ge2$ define
\[
 D_m^\pm=(m-1)[\infty]\pm mD_0.
\]
Their section lattices have rank $m$, and
\[
 \iota_m(D_m^+)=\iota_m(D_m^-)
 \ge \left(\frac18m^2-O(m\log m)\right)\log2.
\]
\end{corollary}

\begin{proof}
Since $\deg D_0=2-6+4=0$, one has $n(D_m^\pm)=m-1$ and
$u_{D_m^\pm}=R(1+x)^{\mp m}$. Multiplication by $R(1+x)^m$ has
triangular matrix with diagonal one. For a determinant-one matrix $T$,
equal covolumes give $[TJ_m:I]=[J_m:I]$, where $I=TJ_m\cap J_m$;
applying $T^{-1}$ shows that $T$ and $T^{-1}$ have equal expansion
indices. The jets of $R^{\pm m}$ have only powers of $2$ as denominators,
so the lower bound of Theorem~\ref{thm:signed-cost} applies to their
rational-lattice index. Theorem~\ref{thm:prod-index-expansion} identifies
this index with the section index.
\end{proof}

\subsection{Dual sections and their boundary contribution}

Assume $n(D)\ge0$. Let $V=\mathbb Z[t]_{\le n(D)}$ and use $\nu_D$ to regard $S_M(D)$ as a
lattice $L\subseteq V$. Let
\[
 C_D:V_{\mathbb Q}\longrightarrow J_{M,\mathbb Q},
 \qquad C_D(p)=j_M(t^{-a}F_D^{-1}p).
\]
For a lattice $B$, write $B^\vee=\operatorname{Hom}_{\mathbb Z}(B,\mathbb Z)$
inside its rational dual. The map $C_D^*$ below is the algebraic transpose.

\begin{proposition}[The dual lattice]
\label{prop:prod-dual}
There is an equality of lattices in $V_{\mathbb Q}^*$,
\[
 \boxed{\quad L^\vee=V^\vee+C_D^*(J_M^\vee).\quad}
\]
This equality is preserved by winding frame changes.
\end{proposition}

\begin{proof}
Every functional on the right is integral on $L$, proving one inclusion.
Choose an integer $d>0$ with $dC_D(V)\subseteq J_M$. The image of $C_D$
modulo $J_M$ is a subgroup $H$ of $(d^{-1}J_M)/J_M$, and $V/L\simeq H$.
An element of $L^\vee/V^\vee$ is a character of $V/L$ with values in
$\mathbb Q/\mathbb Z$. Each character of $H$ extends to the finite group
$(d^{-1}J_M)/J_M$: this follows directly by reducing a subgroup of
$(\mathbb Z/d\mathbb Z)^M$ to Smith normal form. The characters of the
latter group are represented by $J_M^\vee$. Pullback along $C_D$ gives
the reverse inclusion. A frame change replaces $C_D$ by $UC_D$, where
$U$ is multiplication by $(1+x)^k$ and is an integral unimodular
automorphism of $J_M$. Thus $C_D^*U^*(J_M^\vee)=C_D^*(J_M^\vee)$.
\end{proof}

\begin{proposition}[The integral jet duality quotient]
\label{prop:prod-jet-duality-quotient}
In $J_{M,\mathbb Q}$ let $\tau$ take the coefficient of $x^{M-1}$,
and put $\langle f,g\rangle_\tau=\tau(fg)$. For a unit
$u\in J_{M,\mathbb Q}^\times$, set
\[
 L_-=J_M\cap u^{-1}J_M,\qquad L_+=J_M\cap uJ_M.
\]
Then
\[
 L_-^{\vee_\tau}=J_M+uJ_M,
\]
and there is a finite quotient
\[
 0\longrightarrow L_+\longrightarrow L_-^{\vee_\tau}
   \longrightarrow\mathcal B_M(u)\longrightarrow0
\]
whose order satisfies
\[
 \log\#\mathcal B_M(u)
   =e_+(uJ_M\mid J_M)+e_+(J_M\mid uJ_M).
\]
If multiplication by $u$ has determinant of absolute value one, this
quantity is $2e_+(uJ_M\mid J_M)$. In particular, for
$u=R(1+x)^{\pm m}$ and $M=m$, the numerator lattices of
$D_m^+$ and $D_m^-$ have duality quotient satisfying
\[
 \log\#\mathcal B_m(u)
 \ge\left(\frac14m^2-O(m\log m)\right)\log2.
\]
\end{proposition}

\begin{proof}
In the basis $1,x,\ldots,x^{M-1}$, the pairing matrix is the anti-diagonal
matrix with every anti-diagonal entry equal to one. Thus $J_M$ is
unimodular for this pairing. Multiplication by $u$ is self-adjoint for
$\langle\ ,\ \rangle_\tau$, since the jet ring is commutative.
Proposition~\ref{prop:prod-dual}, applied to multiplication by $u$ and
this integral identification with the dual, gives
$L_-^{\vee_\tau}=J_M+uJ_M$. Put $I=J_M\cap uJ_M$ and $S=J_M+uJ_M$.
The isomorphism $S/J_M\simeq uJ_M/I$ gives
\[
 [S:I]=[uJ_M:I][J_M:I],
\]
which proves the order formula. Equal covolumes give equality of the two
indices when the determinant has absolute value one. The final lower
bound follows from Corollary~\ref{cor:prod-signed-index}.
\end{proof}

The algebraic pairing $\tau(fg)$ and the positive metric in the
following theta identity are different pairings.

Equip $V_{\mathbb R}$ with a positive Euclidean metric and all sublattices
with the restricted metric. Write $B^\#$ for the Euclidean dual lattice
and
\[
 \Theta(B)=\sum_{x\in B}e^{-\pi\|x\|^2},\qquad
 h_\theta(B)=\log\Theta(B).
\]

\begin{theorem}[Theta duality for the integral-jet sections]
\label{thm:prod-jet-theta}
For the lattices above,
\[
 h_\theta(L^\#)-h_\theta(V^\#)
   =\iota_M(D)+h_\theta(L)-h_\theta(V).
\]
In particular,
\[
 \iota_M(D)-h_\theta(V)
 \le h_\theta(L^\#)-h_\theta(V^\#)
 \le\iota_M(D).
\]
If the coefficient metric has diagonal squared lengths $w_j^2$ with
\[
 T:=\sum_{j\ge0}\log\Theta(w_j\mathbb Z)<\infty,
\]
then the error in this identity belongs to $[-T,0]$, uniformly in
$D,M$ and $n(D)$.
\end{theorem}

\begin{proof}
Poisson summation gives
$h_\theta(B^\#)-h_\theta(B)=\log\operatorname{covol}(B)$.
Apply it to $L$ and $V$ and subtract. Their covolume ratio is $[V:L]$.
Since $L\subseteq V$, one has $0\le h_\theta(L)\le h_\theta(V)$.
For the diagonal metric, orthogonality gives
$h_\theta(V)=\sum_{j=0}^{n(D)}\log\Theta(w_j\mathbb Z)\le T$.
\end{proof}

The annular coefficient metrics with $w_j=e^{a|j|}(1+|j|)^s$, $a>0$,
satisfy the convergence condition, since
$\Theta(w_j\mathbb Z)-1=O(e^{-\pi w_j^2})$ as $j\to\infty$.
Thus Corollary~\ref{cor:prod-signed-index} also gives a quadratic
relative dual theta size with bounded error.

\subsection{Positive metrics and multiplication of sections}

Suppose the numerator spaces $V_n=\mathbb Z[t]_{\le n}$ carry positive
Hermitian norms satisfying
\[
 \|pq\|_{n+k}\le C\|p\|_n\|q\|_k
\]
with a fixed constant $C\ge1$. Transport these norms to $S_M(D)$ by
$\nu_D$. Then Theorem~\ref{thm:prod-integral-jets} gives
\[
 \|fg\|_{D+E}\le C\|f\|_D\|g\|_E.
\]
Rescaling all norms by $C$ makes multiplication contractive into
the common target metric of degree $n(D)+n(E)$.

\begin{proposition}[Small sections and theta monotonicity]
\label{prop:prod-small}
Use such a contractive family of norms and set
\[
 \widehat\Gamma_M(D)=\{f\in S_M(D):\|f\|_D\le1\}.
\]
These sets are finite. A nonzero $g\in\widehat\Gamma_M(E)$ induces an
injection
\[
 \widehat\Gamma_M(D)\longrightarrow\widehat\Gamma_M(D+E)
\]
and the inequality
\[
 h_\theta(S_M(D))\le h_\theta(S_M(D+E)).
\]
If the norm of the constant numerator $1$ is fixed, then the small-section
counts and theta sizes for all $D$ with $n(D)=0$ are uniformly bounded.
\end{proposition}

\begin{proof}
A bounded ball meets a finite-rank lattice in finitely many points.
Multiplication by $g$ is injective and norm-decreasing, so it gives the
asserted injection. Each summand of the source theta series is bounded
by the summand indexed by its image in the target series. For $n(D)=0$,
the numerator lattice of $S_M(D)$ is a subgroup of $\mathbb Z$ with
the fixed constant norm, which proves the final assertion.
\end{proof}

\subsection{Compatibility with primitive contact}

\begin{proposition}[Tensor products after cyclotomic specialization]
\label{prop:prod-contact}
Let $P,Q$ be distinct monic cyclotomic polynomials and $r,s\ge1$.
For any tensor degrees $m,n$, there is a natural equivalence
\[
 \left(\mathscr L_v^{[m]}\otimes_A A/(P^r)\right)
 \otimes_A^{\mathbb L}
 \left(\mathscr L_v^{[n]}\otimes_A A/(Q^s)\right)
 \simeq
 \mathscr L_v^{[m+n]}\otimes_A A/(P^r,Q^s).
\]
The final coefficient group is finite, of order
$|\operatorname{Res}(P,Q)|^{rs}$.
\end{proposition}

\begin{proof}
The winding coefficient lines are invertible, so their tensor product
can be taken outside the derived coefficient tensor product. The module
$A/(P^r)$ is torsion-free over $\mathbb Z$, and multiplication by $Q^s$
is injective because it is invertible after tensoring with $\mathbb Q$.
Its two-term resolution therefore has no negative homology after the
tensor product. Its cokernel is $A/(P^r,Q^s)$, and its determinant is
$\operatorname{Res}(P,Q)^{rs}$.
\end{proof}

\section{Primitive contacts and polynomial sections}
\label{sec:contact-sections}

Put $S=\Z[t,t^{-1}]$.

\begin{definition}
For a monic polynomial $P$ with constant term $\pm1$ and $r\ge1$ define
the integral thickening
\[
 M(P,r)=\cofib(P(t)^r:HS\longrightarrow HS)\simeq H(S/(P^r)).
\]
For distinct cyclotomic polynomials $P=\Phi_u$, $Q=\Phi_v$, write
$d=\deg P$ and $e=\deg Q$.
\end{definition}

\begin{theorem}[Quadratic contact index]\label{thm:contact}
For $r,s\ge1$,
\[
 M(P,r)\otimes^L_{HS}M(Q,s)\simeq H(S/(P^r,Q^s)),
\]
and the group on the right is finite, with
\[
 \boxed{\log\#S/(P^r,Q^s)=rs\log|\Res(P,Q)|.}
\]
There is a section multiplication presentation
\begin{equation}\label{eq:sylvester}
 0\longrightarrow\Z[t]_{<se}\oplus\Z[t]_{<rd}
 \xrightarrow{(A,B)\mapsto P^rA+Q^sB}
 \Z[t]_{<rd+se}\longrightarrow S/(P^r,Q^s)\longrightarrow0.
\end{equation}
\end{theorem}
\begin{proof}
The derived tensor product and index follow from
Proposition~\ref{prop:prod-contact} in tensor degree zero.
The map in \eqref{eq:sylvester} is the square Sylvester map. Its kernel
is zero by coprimality and the degree bounds. Division by the monic
$P^r$ shows that its cokernel is the same quotient. Since $t$ is
invertible modulo either polynomial, localization adds no relation.
\end{proof}

The spaces in \eqref{eq:sylvester} are ordinary torsion-free polynomial
section lattices. Homogenization realizes the map by sections on the
marked projective compactification $\mathbf P^1_\Z$. The coordinate $t$ identifies these sections with the marked Laurent
coefficient realization of the winding source.

\begin{corollary}[The centered prime coefficient]\label{cor:prime-weight}
For distinct positive integers $u,v$ define the symmetric quantity
\[
 J(u,v)=\frac{\log|\Res(\Phi_u,\Phi_v)|}
             {\sqrt{\varphi(u)\varphi(v)}}.
\]
If $u<v$, this is zero unless $v/u=p^a$ for a prime $p$.
If $v=up^a$ and $p\mid u$,
then $J(u,v)=(\log p)p^{-a/2}$.
\end{corollary}
\begin{proof}
The cyclotomic resultant formula is
$|\Res(\Phi_u,\Phi_v)|=p^{\varphi(u)}$ in the prime-power case and $1$
otherwise \cite[Theorem 9]{BHM}.
For $u=1$ this is the evaluation
$|\Res(\Phi_1,\Phi_v)|=\Phi_v(1)$, given in
\cite[Lemma 4]{BHM}.
For $p\mid u$, $\varphi(up^a)=p^a\varphi(u)$.
\end{proof}

\section{Dual contact lattices and a quadratic theta law}

Let $B_r=\Z[t]/P^r$. The elements
\[
 t^jP^k,\qquad 0\le j<d,\quad0\le k<r
\]
form an integral basis: successive monic division by $P$ proves both
generation and independence. Equip it with the Euclidean metric making
this basis orthonormal. Put $L_{r,s}=Q^sB_r$ and take its Euclidean dual
$L_{r,s}^*$ in $B_r\otimes\R$. For a Euclidean lattice $L$ write
\[
 \Theta(L)=\sum_{x\in L}e^{-\pi\|x\|^2},\qquad
 h_\theta(L)=\log\Theta(L),\qquad
 \vartheta_\Z=\sum_{n\in\Z}e^{-\pi n^2}.
\]

\begin{theorem}[Theta identity for contact lattices]\label{thm:thetaRR}
The contact lattices satisfy
\begin{align}
 h_\theta(L_{r,s}^*)-h_\theta(L_{r,s})
   &=rs\log|\Res(P,Q)|,\label{eq:thetaRR}\\
 rs\log|\Res(P,Q)|\le h_\theta(L_{r,s}^*)
   &\le rs\log|\Res(P,Q)|+rd\log\vartheta_\Z.\label{eq:theta-bound}
\end{align}
In particular, at $r=s=m$ the quadratic coefficient is the contact number, with an $O_P(m)$ remainder.
\end{theorem}
\begin{proof}
Theorem \ref{thm:contact} gives
$\covol(L_{r,s})=[B_r:L_{r,s}]=|\Res(P,Q)|^{rs}$.
Poisson summation gives
$\Theta(L^*)=\covol(L)\Theta(L)$, proving \eqref{eq:thetaRR}.
Since $L_{r,s}\subset B_r$ and the marked basis is orthonormal,
$1\le\Theta(L_{r,s})\le\vartheta_\Z^{rd}$. This proves
\eqref{eq:theta-bound}. The use of log theta as section size is the
usual lattice form of arithmetic Poisson duality
\cite[Proposition 1 and \S8]{GeerSchoof}; the quadratic coefficient here comes
from the thickened source index.
\end{proof}

\begin{proposition}[A perfect contact duality]\label{prop:linking}
Let $F=P^r$ have degree $n$, and let $\tau:B_r\to\Z$ take the
coefficient of $t^{n-1}$ in the monic remainder. Then
$(b,c)\mapsto\tau(bc)$ is a unimodular symmetric pairing. On
$A=B_r/Q^sB_r$,
\[
 \lambda(\bar b,\bar c)=\tau(bc/Q^s)\pmod\Z
\]
is a perfect symmetric pairing with values in $\Q/\Z$.
\end{proposition}
\begin{proof}
In the monomial basis the matrix entry $\tau(t^it^j)$ is zero when
$i+j<n-1$ and one when $i+j=n-1$, so its determinant is $\pm1$.
Multiplication by $Q^s$ is self-adjoint for this pairing and invertible
over $\Q$. Hence changing either lift changes $\lambda$ by an integer.
If $\lambda(\bar b,-)=0$, unimodularity implies $b/Q^s\in B_r$, proving
nondegeneracy. The two finite groups have equal order, so the induced
map to the character dual is an isomorphism.
\end{proof}

Under the simultaneous substitution $t\mapsto t^k$ with $k\ge1$, the marked basis
splits into $k$ residue blocks. Both the theta sum and the index are
products of the corresponding base quantities; their logarithms multiply
by $k$. This is a metric statement for the specified marked basis.
\section{Supported contacts and modulus data}
\label{sec:cycles}

Let $S=\Z[t,t^{-1}]$ and $X=\operatorname{Spec}S$. This is the
integral coefficient realization of the marked winding coordinate.
For distinct cyclotomic polynomials $P,Q$ and integers $r,s\ge1$, put
$T=V(P,Q)$ and $A_{r,s}=S/(P^r,Q^s)$. Write $K_0^T(X)$ for the Grothendieck group of
perfect complexes on $X$ whose cohomology is supported on $T$.

\begin{proposition}[The contact class with support]
\label{prop:cycles:support}
The Koszul complex of $P^r,Q^s$ represents $A_{r,s}$ in $K_0^T(X)$.
One has
\[
 [A_{r,s}]=rs[A_{1,1}],\qquad
 \sum_{x\in T}\operatorname{length}_{S_x}(A_{r,s,x})\log\#\kappa(x)
   =rs\log|\operatorname{Res}(P,Q)|.
\]
The image of this class in $K_0(X)$ is zero. Moreover $\operatorname{CH}^2(X)=0$.
\end{proposition}
\begin{proof}
The two polynomials are a regular sequence: $S/(P^r)$ is free over
$\Z$, and $Q^s$ is invertible on its rationalization, hence is a
nonzerodivisor. The Koszul resolution follows. Filtering successively
by the powers of $P$ and $Q$ gives $rs$ subquotients isomorphic to
$S/(P,Q)$. Each module has a finite free resolution, so additivity
applies with support. The arithmetic degree is the logarithm of the
order of the finite module, computed in Theorem~\ref{thm:contact}.
Forgetting support gives the alternating class $[S]-2[S]+[S]=0$.
Every closed point of $X$ has the form $(p,f(t))$, with $f\ne t$
irreducible over $\F_p$. On the vertical curve
$\operatorname{Spec}\F_p[t,t^{-1}]$, the divisor of the rational
function $f$ is that single point with multiplicity one. Hence every
codimension-two cycle is rationally equivalent to zero.
\end{proof}

\subsection{A nilpotent modulus in the marked coordinate}

Put $x=t-1$, $B_m=\Z_2[x]/(x^m)$ and
$A_m=\Q_2[x]/(x^m)$. For $m\ge3$ consider
\[
 R(t)=\frac{\Phi_2(t)^2\Phi_8(t)}{\Phi_4(t)^3}.
\]
The numerator and denominator have degree six, and $R(1)=1$.

\begin{proposition}[A relative unit with trivial transfer determinant]
\label{prop:cycles:jet}
In $A_m$, the augmentation-one unit $R(1+x)^m$ is nontrivial and
\[
 \log R(1+x)^m=\frac{3m}{4}x^2+O(x^3).
\]
Its determinant as a multiplication operator on the $\Q_2$-space $A_m$
is one. For the modulus $D_m=m[1]$ on $\mathbf P^1_{\Q_2}$,
\[
 \operatorname{Pic}^0(\mathbf P^1_{\Q_2},D_m)
 \simeq 1+x\Q_2[x]/(x^m).
\]
Thus the relative unit records information that the transfer determinant
does not retain.
\end{proposition}
\begin{proof}
Using $\Phi_2=t+1$, $\Phi_4=t^2+1$, $\Phi_8=t^4+1$, Taylor expansion
gives $R(1+x)=1+3x^2/4+O(x^3)$. The logarithm is a finite nilpotent
series. Its displayed coefficient is nonzero in $\Q_2$.
In the ordered monomial basis, multiplication by an augmentation-one
unit is triangular with all diagonal entries one.
The units and Picard sequence for the modulus gives
\[
 \operatorname{Pic}^0(\mathbf P^1_{\Q_2},D_m)
   =A_m^\times/\Q_2^\times.
\]
Every unit splits uniquely into its constant term and an
augmentation-one unit. This proves the identification. The kernel
of $K_1(A_m)\to K_1(\Q_2)$ is the same augmentation-one unit group,
since these are commutative local rings.
\end{proof}

Here $\operatorname{Pic}^0$ is the group of degree-zero line bundles
with a trivialization on $D_m$. This is the weight-one modulus
construction treated in \cite[Theorem 4.1 and Corollary 4.2]{BindaSaito}.
Its use here retains the chosen jet at $t=1$. At the integral prime
$2$, the cyclotomic divisors meet that modulus; the integral contact
is recorded by $B_m$ and its fractional lattices in $A_m$.

\subsection{Arithmetic intersection of the contact divisors}

Homogenizing $P^r,Q^s$ on $\mathbf P^1_\Z$ gives a finite intersection
$Z$. Put $a=r\deg P$ and $b=s\deg Q$. The vector-bundle complex
\[
 [\mathcal O(-a-b)\longrightarrow
   \mathcal O(-a)\oplus\mathcal O(-b)\longrightarrow\mathcal O]
\]
has sole cohomology $\mathcal O_Z$ in its last degree.
Twisting by $\mathcal O(a+b-1)$ recovers \eqref{eq:sylvester}, since
$H^i(\mathcal O(-1))=0$ for $i=0,1$.
Give its three vector-bundle terms smooth conjugation-invariant
Hermitian metrics.
The arithmetic Chern character identity
\cite[Theorem~1.1, p.~614]{GilletSoule2014} expresses the defect of this
generically exact complex as a finite-support class and a Bott--Chern
form. For this Koszul complex the finite contribution has degree
$\log\#\Gamma(Z,\mathcal O_Z)$
\cite[equation~(2.8) and Lemma~2.2, pp.~617--618]{GilletSoule2014}.
Theorem~\ref{thm:contact} computes this degree as
$rs\log|\operatorname{Res}(P,Q)|$. The form records the metrics on
these same bundles.
\section{Hermitian metrics on divisor sections}\label{sec:section-metrics}

The numerator identifications $\nu_D$ carry common positive Hermitian
metrics to the divisor section lattices $S_M(D)$.

\subsection{A tensor-compatible metric}

Let $V_n=\mathbb C[t]_{\le n}$ for $n\ge0$.  Homogenization identifies
$V_n$ with the homogeneous polynomials of degree $n$ in two variables.
Define
\begin{equation}\label{eq:pos-bombieri}
 \left\|\sum_{j=0}^{n}a_jt^j\right\|_{B,n}^{2}
   =\sum_{j=0}^{n}\binom nj^{-1}|a_j|^2.
\end{equation}
This is the Bombieri--Weyl normalization of the symmetric tensor metric;
the relation with the apolar normalization is
$\|P\|_{B,n}^2=\|P\|_{\mathrm{ap}}^2/n!$
\cite[\S2, equation (10) and the following display]{AldazRender2024}.

\begin{proposition}[Multiplication and the Fubini--Study metric]
\label{prop:pos-bombieri-product}
For $p\in V_n$ and $q\in V_m$,
\begin{equation}\label{eq:pos-contract}
 \|pq\|_{B,n+m}\le\|p\|_{B,n}\|q\|_{B,m}.
\end{equation}
All these norms are obtained from the tensor powers of one
Fubini--Study metric on $\mathcal O_{\mathbb P^1}(1)$:
\begin{equation}\label{eq:pos-fs}
 \|p\|_{B,n}^2
  =(n+1)\int_{\mathbb C}
       \frac{|p(z)|^2}{(1+|z|^2)^n}
       \frac{dA(z)}{\pi(1+|z|^2)^2}.
\end{equation}
The factor $n+1$ normalizes the section norm.
\end{proposition}
\begin{proof}
In $\operatorname{Sym}^{n}\mathbb C^2$, the monomial $X^{n-j}Y^j$
is the average of $\binom nj$ orthonormal elementary tensors, so its
squared norm is $\binom nj^{-1}$. Multiplication is orthogonal
symmetrization of the tensor product, proving \eqref{eq:pos-contract}.
Rotation makes monomials orthogonal in \eqref{eq:pos-fs}, and beta
integration gives
\[
 (n+1)\int_0^\infty \frac{u^j}{(1+u)^{n+2}}\,du
 =\frac{j!(n-j)!}{n!}.
\]
This proves the integral formula.
\end{proof}

For $n(D)\ge0$, give $f\in S_M(D)$ the norm
$\|f\|_D=\|\nu_D(f)\|_{B,n(D)}$. The numerator product identity gives
\begin{equation}\label{eq:pos-divisor-contract}
 \|fg\|_{D+E}\le\|f\|_D\|g\|_E.
\end{equation}
Thus Proposition~\ref{prop:prod-small} applies to these lattices with
one metric on each target. For a Laurent interval $I=[a,b]\cap\mathbb Z$,
set $\|t^ap\|_I=\|p\|_{B,b-a}$. Products use $I+J$, and $t^k$
is an isometry from the interval $I$ to $I+k$.

\subsection{The size of the section lattice}

Put $\Lambda_n=(\mathbb Z[t]_{\le n},\|\cdot\|_{B,n})$.
For lattice theta size and Poisson duality, see
\cite[\S0.2.4.1, equations (0.2.3)--(0.2.7)]{BostCharlesTheta}.

\begin{proposition}[Quadratic theta size]\label{prop:pos-theta-size}
With $\vartheta_{\mathbb Z}=\sum_{k\in\mathbb Z}e^{-\pi k^2}$,
\begin{equation}\label{eq:pos-theta-bound}
 0\le h_\theta(\Lambda_n)
       -\frac12\sum_{j=0}^{n}\log\binom nj
    \le(n+1)\log\vartheta_{\mathbb Z}.
\end{equation}
Consequently,
\[
 h_\theta(\Lambda_n)=\frac14n^2+O(n\log(n+2)).
\]
\end{proposition}
\begin{proof}
The monomials form an orthogonal integral basis.  Thus
\[
 \log\operatorname{covol}(\Lambda_n)
  =-\frac12\sum_{j=0}^{n}\log\binom nj.
\]
The dual lattice has orthogonal basis lengths $\binom nj^{1/2}\ge1$,
so its theta logarithm lies between zero and
$(n+1)\log\vartheta_{\mathbb Z}$.  Poisson summation gives
$h_\theta(\Lambda_n)=-\log\operatorname{covol}(\Lambda_n)
+h_\theta(\Lambda_n^*)$, proving \eqref{eq:pos-theta-bound}.
For the asymptotic, write
\[
 \sum_{j=0}^{n}\log\binom nj
  =(n+1)\log(n!)-2\sum_{j=0}^{n}\log(j!).
\]
Stirling's estimate and summation, or integration of
$-x\log x-(1-x)\log(1-x)$ on $[0,1]$, yield
$n^2/2+O(n\log(n+2))$.
\end{proof}

The construction extends to a number field $K$.  Give
$\mathcal O_K[t]_{\le n}$ the norm
\[
 \|p\|_{K,n}^2
   =\sum_{\sigma:K\hookrightarrow\mathbb C}
       \|\sigma(p)\|_{B,n}^2,
\]
with conjugate embeddings both included.  Multiplication remains
contractive since, for nonnegative $A_\sigma,B_\sigma$,
\[
 \sum_\sigma A_\sigma B_\sigma
 \le\Bigl(\sum_\sigma A_\sigma\Bigr)
       \Bigl(\sum_\sigma B_\sigma\Bigr).
\]
With $d=[K:\mathbb Q]$ and discriminant $\Delta_K$, the covolume is
\begin{equation}\label{eq:pos-number-field-covol}
 \log\operatorname{covol}(\mathcal O_K[t]_{\le n})
 =\frac{n+1}{2}\log|\Delta_K|
   -\frac d2\sum_{j=0}^{n}\log\binom nj.
\end{equation}
Each Minkowski coefficient lattice has covolume $|\Delta_K|^{1/2}$
and real dimension $d$; scaling by $\binom nj^{-1/2}$ proves
\eqref{eq:pos-number-field-covol}. Up to coordinatewise complex
conjugation, an isometry, its Euclidean dual is the orthogonal sum of
$\binom nj^{1/2}\mathcal O_K^\vee$, with $\mathcal O_K^\vee$ the
trace-dual lattice. Its theta sum is at most
$\Theta(\mathcal O_K^\vee)^{n+1}$, so Poisson summation gives
$h_\theta=dn^2/4+O_K(n\log(n+2))$. Together with the multiplicative
coefficient Gauss norms at finite places, these norms give an adelic
metric on each section lattice.

\subsection{Determinant norms of contact maps}

Let $P,Q\in\mathbb Z[t]$ be distinct monic cyclotomic polynomials,
let $r,s\ge1$, and put $a=r\deg P$, $b=s\deg Q$, $n=a+b$.
Equip the source and target of the Sylvester map
\[
 A_{P^r,Q^s}:\mathbb Z[t]_{<b}\oplus\mathbb Z[t]_{<a}
 \longrightarrow\mathbb Z[t]_{<n},\qquad
 (u,v)\longmapsto P^ru+Q^sv,
\]
with the Bombieri metrics in degrees $b-1,a-1,n-1$ respectively.
Define
\[
 G(N)=-\sum_{j=0}^{N}\log\binom Nj,\qquad
 E_B(a,b)=\frac12\{G(a+b-1)-G(a-1)-G(b-1)\}.
\]
The determinant formula for a map between Euclidean lattices gives
\begin{equation}\label{eq:pos-resultant-metric}
 \log\|\det A_{P^r,Q^s}\|
   =rs\log|\operatorname{Res}(P,Q)|+E_B(a,b).
\end{equation}
If $C_P=\|P^r\|_{B,a}$ and $C_Q=\|Q^s\|_{B,b}$, then
\begin{equation}\label{eq:pos-determinant-bound}
 rs\log|\operatorname{Res}(P,Q)|+E_B(a,b)
 \le\frac n2\log(C_P^2+C_Q^2).
\end{equation}
Indeed the two multiplication maps have norms at most $C_P,C_Q$,
so their direct sum has norm at most $(C_P^2+C_Q^2)^{1/2}$ and its
determinant norm is bounded by the $n$th power.

\section{A Hilbert space of Laurent sections}
\label{sec:annular-sections}

Here $t$ is the marked Laurent coordinate.

\begin{definition}\label{def:pos-annular}
For $a>0$ and $s>1/2$, set
\[
 w_{a,s}(n)=e^{a|n|}(1+|n|)^s,\qquad n\in\mathbb Z.
\]
Let $\mathcal H_{a,s}$ consist of Laurent series
$f(t)=\sum_{n\in\mathbb Z}c_nt^n$ with
\begin{equation}\label{eq:pos-annular-norm}
 \|f\|_{a,s}^2
  =\sum_{n\in\mathbb Z}w_{a,s}(n)^2|c_n|^2<\infty.
\end{equation}
The Hermitian inner product is linear in its first argument.
\end{definition}

Put
\[
 K_s=\left(\sum_{n\in\mathbb Z}(1+|n|)^{-2s}\right)^{1/2},
 \qquad C_s=2\max(1,2^{s-1})K_s.
\]

\begin{theorem}[Products and integral theta size]
\label{thm:pos-annular-algebra}
Every element of $\mathcal H_{a,s}$ is holomorphic on
$e^{-a}<|t|<e^a$ and continuous on its closure.  Laurent polynomials
are dense, and pointwise multiplication is continuous with
\begin{equation}\label{eq:pos-annular-product}
 \|fg\|_{a,s}\le C_s\|f\|_{a,s}\|g\|_{a,s}.
\end{equation}
The involution $f^\star(t)=\overline{f(1/\bar t)}$ is isometric.
The integral Laurent module $\mathbb Z[t,t^{-1}]$ has finite theta sum
\begin{equation}\label{eq:pos-annular-theta}
 \Theta_{a,s}
  =\prod_{n\in\mathbb Z}
       \left(\sum_{k\in\mathbb Z}e^{-\pi w_{a,s}(n)^2k^2}\right)
  <\infty.
\end{equation}
\end{theorem}
\begin{proof}
Cauchy--Schwarz gives $\sum_n|c_nt^n|\le K_s\|f\|_{a,s}$ and
uniform convergence on the closed annulus. Derivative series converge
on smaller annuli. Truncation proves density.

Write $A_n=e^{a|n|}|c_n|$ and $B_n=e^{a|n|}|d_n|$ for the two
coefficient sequences.  The exponential weight is submultiplicative,
and
\[
 (1+|k+\ell|)^s
 \le c_s\{(1+|k|)^s+(1+|\ell|)^s\},
 \qquad c_s=\max(1,2^{s-1}).
\]
Young's $\ell^2*\ell^1\to\ell^2$ inequality therefore bounds the
product norm by
\[
 c_s\bigl(\|(1+|n|)^sA_n\|_2\|B\|_1
       +\|A\|_1\|(1+|n|)^sB_n\|_2\bigr).
\]
Cauchy--Schwarz bounds each $\ell^1$ norm by $K_s$ times its weighted
$\ell^2$ norm, proving \eqref{eq:pos-annular-product}. Even weights
make coefficient conjugation and reversal isometric.

Integral finite-norm sequences have finite support since
$w_{a,s}(n)\to\infty$. Monotone convergence over finite intervals gives
\eqref{eq:pos-annular-theta}. Its finiteness follows by summing, at
$x=w_{a,s}(n)$, the bound
\[
 \sum_{k\in\mathbb Z}e^{-\pi x^2k^2}-1
  \le\frac{2e^{-\pi x^2}}{1-e^{-\pi x^2}}.
\]
\end{proof}

The integral coefficient module is the directed union of its finite
Laurent spans with their induced metrics, an ind-Euclidean lattice in
the terminology of \cite[Corollary~5.8.3; \S6.1, (6.1.1)--(6.1.3)]{Bost2020}.

The coefficient cone $c_n\ge0$ is closed, proper and preserved by
products and frame shifts. For each integer $h\ge1$, power pullback
sends it to the corresponding cone in $\mathcal H_{a/h,s}$
(Proposition~\ref{prop:pos-annular-covers}). Coefficient collection modulo
$L$ maps it onto the cone in $\mathbb R[C_L]$. These statements follow
coefficientwise and by continuity of multiplication.

On $\mathbb Z[t]_{\le N}$, the diagonal metric change gives
\begin{equation}\label{eq:pos-metric-comparison}
 \|p\|_{B,N}^2\le\|p\|_{a,s}^2
 \le e^{2aN}(N+1)^{2s}2^N\|p\|_{B,N}^2.
\end{equation}
The logarithmic covolume change from the Bombieri metric to the
annular coefficient metric is
\begin{equation}\label{eq:pos-metric-covol-change}
 \frac{aN(N+1)}2+s\log((N+1)!)
   +\frac12\sum_{j=0}^{N}\log\binom Nj.
\end{equation}
Both formulas follow from the squared-length ratio
$w_{a,s}(j)^2\binom Nj$. Every integral Laurent interval has annular
theta size at most $\log\Theta_{a,s}$.

\subsection{Bounded division and jet quotients}

\begin{lemma}[Division at a unit-circle point]
\label{lem:pos-division}
Let $|\zeta|=1$ and let $f\in\mathcal H_{a,s}$ satisfy $f(\zeta)=0$.
Then $q=f/(t-\zeta)$ belongs to $\mathcal H_{a,s}$ and
\begin{equation}\label{eq:pos-division-bound}
 \|q\|_{a,s}\le\frac{1}{1-e^{-a}}\|f\|_{a,s}.
\end{equation}
\end{lemma}
\begin{proof}
Write $f=\sum f_nt^n$.  The quotient coefficients are
\begin{equation}\label{eq:pos-division-tails}
 q_n=
 \begin{cases}
  \displaystyle\sum_{k>n}\zeta^{k-n-1}f_k,&n\ge0,\\[2mm]
  \displaystyle-\sum_{k\le n}\zeta^{k-n-1}f_k,&n<0.
 \end{cases}
\end{equation}
The sums converge absolutely, and $f(\zeta)=0$ gives
$f_n=q_{n-1}-\zeta q_n$ also at degree zero. For $n\ge0$,
$w(n)/w(n+\ell)\le e^{-a\ell}$ ($\ell\ge1$); for $n<0$,
$w(n)/w(n-\ell)\le e^{-a\ell}$ ($\ell\ge0$).
Young's inequality bounds the two disjoint coefficient parts by
$(e^a-1)^{-1}$ and $(1-e^{-a})^{-1}$ respectively. Adding their
squared bounds proves \eqref{eq:pos-division-bound}, and the Laurent
identity identifies the holomorphic quotient.
\end{proof}

\begin{theorem}[Closed cyclotomic ideals]\label{thm:pos-jet-quotient}
Let $P(t)$ be a monic polynomial of degree $d$ whose roots all lie on
the unit circle, and let $r\ge1$.  For each distinct root $\zeta$ of
$P$, let $m_\zeta$ be its multiplicity.  Then
\begin{equation}\label{eq:pos-jet-kernel}
 P^r\mathcal H_{a,s}
   =\{f\in\mathcal H_{a,s}:
       f^{(j)}(\zeta)=0\quad
       (0\le j<rm_\zeta,\ P(\zeta)=0)\}.
\end{equation}
This ideal is closed.  Multiplication by $P^r$ is injective, has a
bounded inverse on its image, and its inverse has norm at most
$(1-e^{-a})^{-rd}$.  Joint evaluation of the $rd$ jets induces a
linear isomorphism
\[
 \mathcal H_{a,s}/P^r\mathcal H_{a,s}
     \longrightarrow\mathbb C^{rd}.
\]
The quotient has a positive definite Hermitian inner product and a
continuous multiplication induced by pointwise multiplication.
\end{theorem}
\begin{proof}
Every derivative evaluation at $|\zeta|=1$ is continuous: its
coefficient bound follows from
\[
 \sum_{n\in\mathbb Z}
  |n(n-1)\cdots(n-j+1)|^2w_{a,s}(n)^{-2}<\infty.
\]
The inclusion from left to right in \eqref{eq:pos-jet-kernel} follows
from differentiation.  For the converse, factor $P^r$ into its $rd$
linear factors and apply Lemma \ref{lem:pos-division} successively.
The remaining zeros have the required multiplicities after each
division.  This also proves the inverse bound.  The right side of
\eqref{eq:pos-jet-kernel} is the intersection of finitely many closed
kernels.  Hermite interpolation by polynomials of degree below $rd$
proves surjectivity of the joint jet map.  The Hilbert quotient by
a closed subspace is positive definite; its multiplication is
bounded using minimum-norm lifts and
\eqref{eq:pos-annular-product}.
\end{proof}

The quotient inner product is explicit.  Index the jet coordinates by
$\nu=(\zeta,j)$ and put
\[
 b_{\nu,n}=n(n-1)\cdots(n-j+1)\zeta^{n-j},
 \qquad
 K_{\nu\mu}
  =\sum_{n\in\mathbb Z}
       b_{\nu,n}\overline{b_{\mu,n}}w_{a,s}(n)^{-2}.
\]
Surjectivity of the joint evaluation $J$ makes $K=JJ^*$ positive
definite. A jet vector $v$ has minimum-norm lift $J^*K^{-1}v$ and
squared quotient norm $v^*K^{-1}v$.

If the polynomial $P$ above has integral coefficients and constant term $\pm1$, the
lattice $\mathbb Z[t,t^{-1}]/(P^r)$ embeds in this jet quotient.

\subsection{Finite covers and deck-invariant metrics}

\begin{proposition}[Power pullback and averaging]
\label{prop:pos-annular-covers}
For an integer $h\ge1$, the map $F_hf(t)=f(t^h)$ satisfies
\[
 \|F_hf\|_{a/h,s}\le h^s\|f\|_{a,s}.
\]
The averaging map
\[
 A_h\left(\sum_{n\in\mathbb Z}c_nt^n\right)
     =\sum_{n\in\mathbb Z}c_{hn}t^n
\]
is a contraction from $\mathcal H_{a/h,s}$ to $\mathcal H_{a,s}$,
with $A_hF_h=1$.  The operator $F_hA_h$ is the orthogonal projection
onto exponents divisible by $h$ and is equal to
$h^{-1}\sum_{\zeta^h=1}f(\zeta t)$.  These maps compose under towers
of power covers.
\end{proposition}
\begin{proof}
For the pullback the exponential factor is unchanged, while
$1+h|n|\le h(1+|n|)$ bounds the Sobolev factor.  For $A_h$,
$w_{a,s}(n)\le w_{a/h,s}(hn)$, which proves contraction.
The remaining assertions follow coefficient by coefficient.
\end{proof}

For $b,k\in\Z$, with weights $w_{a,s}(n-b)$, multiplication by $t^k$ is isometric
from frame centre $b$ to $b+k$; products add the centres.

\begin{proposition}[The metric on a finite winding cover]
\label{prop:pos-finite-cover-metric}
For $L\ge1$, identify
$\mathcal H_{a,s}/(t^L-1)\mathcal H_{a,s}$ with
$\mathbb C[C_L]$ by collecting coefficients modulo $L$.
Set
\[
 \beta_{L,r}=\sum_{k\in\mathbb Z}w_{a,s}(r+kL)^{-2},
 \qquad 0\le r<L.
\]
In residue coefficients $b_0,\ldots,b_{L-1}$, the quotient norm is
\begin{equation}\label{eq:pos-cyclic-quotient}
 \|b\|_{\mathrm{quot}}^2
   =\sum_{r=0}^{L-1}\frac{|b_r|^2}{\beta_{L,r}}.
\end{equation}
Its average under the $L$ deck translations is
\begin{equation}\label{eq:pos-cyclic-average}
 \frac1L\sum_{j=0}^{L-1}\|t^jb\|_{\mathrm{quot}}^2
   =c_L\sum_{r=0}^{L-1}|b_r|^2,
 \qquad
 c_L=\frac1L\sum_{r=0}^{L-1}\beta_{L,r}^{-1}.
\end{equation}
Multiplying this averaged squared norm by $L/c_L$ gives the
deck-invariant squared norm $L\sum_r|b_r|^2$ on the finite cover.
\end{proposition}
\begin{proof}
Coefficient collection is continuous by weighted Cauchy--Schwarz.
Its kernel consists of functions vanishing at all $L$th roots of
unity, so Theorem \ref{thm:pos-jet-quotient} identifies its kernel
with $(t^L-1)\mathcal H_{a,s}$.  For fixed residue $r$, minimize
$\sum_k w_{a,s}(r+kL)^2|c_{r+kL}|^2$ subject to
$\sum_k c_{r+kL}=b_r$.  Cauchy--Schwarz gives minimum
$|b_r|^2/\beta_{L,r}$, attained at
\[
 c_{r+kL}=\frac{b_r w_{a,s}(r+kL)^{-2}}{\beta_{L,r}}.
\]
The residue classes are orthogonal, proving
\eqref{eq:pos-cyclic-quotient}.  Deck translation permutes the
residue coefficients cyclically; averaging therefore replaces each
coefficient weight by their arithmetic mean, giving
\eqref{eq:pos-cyclic-average}.
\end{proof}

\section{Weighted Poisson quotients and zero evaluation}
\label{sec:zero-comparison}

Let $d=[K:\Q]$, and let $\sigma_K(r)$ be the idele class with finite
coordinates $1$ and archimedean coordinates $r^{1/d}$. Since the
archimedean degrees sum to $d$, one has $|\sigma_K(r)|=r$, giving
$C_K\simeq C_K^1\times\R_{>0}$.
Fix a continuous unitary character $\chi$ of $C_K^1$, extended trivially
on $\sigma_K(\R_{>0})$, and normalize the Haar measure of $C_K^1$
to have mass one. We use the character component
$f(k\sigma_K(r))=\chi(k)^{-1}f(\sigma_K(r))$.
Write $t=\log|c|$ and center the radial function by
$g(t)=e^{t/2}f(\sigma_K(e^t))$. After centering, convolution in the character
component becomes ordinary convolution in $t$, and the involution
\eqref{eq:ana:sharp} becomes $g^\sharp(t)=\overline{g(-t)}$.

Let $E_\infty$ be the space of smooth functions decreasing, with all
derivatives, faster than every exponential in both directions, with the
Fr\'echet topology given by all positive integer exponential weights and Schwartz
seminorms. Let
$P\subset E_\infty$ be the centered image of $\theta H_\cap$ in the
chosen character component. Equivariance gives translation invariance,
and $P$ is a convolution ideal \cite[\S5.6, p.~45]{Meyer}.
For $f=\Sigma\Phi$ with both moments zero, \eqref{eq:ana:poisson}
gives $f^\sharp=\Sigma\overline{\widehat\Phi}$. Both moments of
$\overline{\widehat\Phi}$ vanish, and unitarity of $\chi$ keeps
$f^\sharp$ in the same component. Thus $P$ is $\sharp$-invariant. Put
\[
 F_g(z)=\int_\R e^{zt}g(t)\,dt.
\]
The Mellin identity for the theta image has the completed Hecke
$L$-function $L(s,\chi)$ as a factor on the open critical strip.
Thus, if $\rho$ is a nontrivial zero of $L(s,\chi)$ of order $m_\rho$, then
\begin{equation}
 F_p^{(j)}(\rho-\tfrac12)=0
 \quad(p\in P,\ 0\le j<m_\rho).
 \label{eq:zero:input}
\end{equation}
This is the divisibility used in the proof of the spectral computation
\cite[Theorem 5.11]{Meyer}.

The number field \(K\) is arbitrary. For the trivial character, these are the
nontrivial zeros of the Dedekind zeta function \(\zeta_K\). More generally,
every continuous unitary Hecke character of \(C_K\) has the form
\(\chi|\cdot|^{i\tau}\), with \(\tau\in\mathbb R\) and \(\chi\) normalized
as above. A pure norm twist translates the zeros vertically. Thus the
character components considered here cover the unitary Hecke family up
to this translation, including characters of infinite order
\cite[pp.~1--2 and Theorem 5.11]{Meyer}.

\begin{definition}[A fixed-weight Poisson quotient]
\label{def:zero:fixed}
For $\alpha>1/2$, set
\begin{align*}
 E_\alpha&=\{g\in C^\infty(\R):p_{N,j}(g)<\infty
                         \text{ for all }N,j\ge0\},\\
 p_{N,j}(g)&=\sup_{t\in\R}e^{\alpha|t|}(1+|t|)^N|g^{(j)}(t)|,\\
 P_\alpha&=\overline P^{\,E_\alpha},\qquad
 Q_\alpha=E_\alpha/P_\alpha.
\end{align*}
The topology on $E_\alpha$ is defined by these seminorms and that on
$Q_\alpha$ is the quotient topology.
\end{definition}

\begin{theorem}[Preservation of finite zero jets]
\label{thm:zero:jets}
The following statements hold.
\begin{enumerate}[label=\textup{(\roman*)}]
\item $E_\alpha$ and $Q_\alpha$ are nuclear Fr\'echet spaces.
Both $C_c^\infty(\R)$ and $E_\infty$ are dense in $E_\alpha$.
\item Convolution, $\sharp$, and smooth translations descend to
$Q_\alpha$.
\item For every actual nontrivial zero $\rho$ and $0\le j<m_\rho$,
\[
 \ell_{\rho,j}(g)=\int_\R t^j e^{(\rho-1/2)t}g(t)\,dt
\]
descends to a continuous functional on $Q_\alpha$. The joint evaluation
map for any finite set of distinct pairs $(\rho,j)$ is surjective.
\item The natural map $Q_\infty=E_\infty/P\to Q_\alpha$ is continuous
with dense image. Its transpose on continuous duals is injective and
sends each displayed zero jet to its original evaluation functional.
\end{enumerate}
\end{theorem}

\begin{proof}
Multiplication by $\cosh(\alpha t)$ is a Fr\'echet isomorphism
$E_\alpha\to\mathcal S(\R)$. Indeed its derivatives are bounded by
$C_j e^{\alpha|t|}$, and the derivatives of its reciprocal are bounded
by $C_j'e^{-\alpha|t|}$. Leibniz's formula proves the two families of
seminorm inequalities. Completeness and nuclearity follow, as does
nuclearity of the quotient by a closed subspace.

If $\eta$ is smooth, compactly supported and equal to one near zero,
then $\eta(t/R)g\to g$ in every $p_{N,j}$. The tail is bounded using
$p_{N+1,j}$ and a factor $O(R^{-1})$; derivatives of the cutoff satisfy
the same estimate with finitely many derivative seminorms. This proves
density. The inclusion $E_\infty\subset E_\alpha$ is strict, as shown
by $\exp(-\beta\sqrt{1+t^2})$ for $\beta>\alpha$.

The inequalities
\[
 e^{\alpha|t|}\le e^{\alpha|u|}e^{\alpha|t-u|},\qquad
 1+|t|\le(1+|u|)(1+|t-u|)
\]
give the continuous convolution estimate
\begin{equation}
 p_{N,j}(g*h)\le 2p_{N+2,0}(g)p_{N,j}(h).
 \label{eq:zero:convolution}
\end{equation}
Here $2=\int_\R(1+|u|)^{-2}du$. Choose compactly supported $g_n\to g$
and $p_n\in P$ tending to $p\in P_\alpha$. Then $g_n*p_n\in P$ and
converges to $g*p$, proving the ideal property. The involution preserves
all seminorms. For $T_a g(t)=g(t-a)$,
\[
 p_{N,j}(T_a g)\le e^{\alpha|a|}(1+|a|)^N p_{N,j}(g).
\]
Differentiation in $a$ is controlled by the next derivative seminorm,
so the translation action is smooth. Its stability and that of $\sharp$
extend from $P$ to its closure.

For $z=\rho-1/2$ one has $|\Re z|<1/2<\alpha$, and
\begin{equation}
 |\ell_{\rho,j}(g)|\le
 \frac{2j!}{(\alpha-|\Re z|)^{j+1}}p_{0,0}(g).
 \label{eq:zero:jet-bound}
\end{equation}
Equation~\eqref{eq:zero:input} and continuity show that each functional
annihilates $P_\alpha$. A linear relation among any finite set
of the functionals on $C_c^\infty$ would make
\[
 \sum_{\rho,j}a_{\rho,j}t^j e^{(\rho-1/2)t}
\]
the zero distribution. Independence of exponential polynomials gives
$a_{\rho,j}=0$ for all indices. The finite-dimensional target has no
nonzero functional annihilating the image, proving surjectivity. Finally,
density of $E_\infty$ gives density of the quotient map and injectivity
of its transpose.
\end{proof}

The weight of a zero jet can be read without a convention for dual group
actions. With $R_a=e^{a/2}T_a$, change of variable gives
\begin{equation}
 \ell_{\rho,k}(R_a g)=e^{\rho a}
 \sum_{j=0}^k\binom{k}{j}a^{k-j}\ell_{\rho,j}(g).
 \label{eq:zero:flow}
\end{equation}
This is the transpose action on functionals; the contragredient action
uses $R_{-a}$ instead. The evaluation formula describes the action on each finite zero-jet block.

\section{The arithmetic winding curve}
\label{sec:assembly}

\begin{definition}[Arithmetic winding curve]
\label{def:assembly:curve}
Fix a number field \(K\), a full lattice \(\Lambda\subset K\), the trace character
\(\psi_K\), and the uniformizers \(\varpi_v\) at the finite places. For each finite
place choose the residue frame \((u_v,\chi_v)\) of
Definition~\ref{def:res:frame}. The arithmetic winding curve associated to these
data, denoted by
\[
 \mathfrak W(K,\Lambda,\psi_K;
       (\varpi_v,u_v,\chi_v)_{v\nmid\infty}),
\]
is the following structured object:
\begin{enumerate}[label=\textup{(\roman*)}]
\item the \(C_K\)-space \(Y_K\), its open charts \(N_v\), their common open
      \(C_K\), and the boundary orbits \(Z_v\);
\item the principal covers \(Q_v\), their finite quotients, and the covers
      \(\mathcal Q_{v,H}\) for open subgroups \(H\subseteq\mathcal O_v^\times\), with
      the quotient maps and their finite stable transfers;
\item the sheaves \((\cB_v,\cD_v,\cL_v)\), their algebra and module maps,
      and the system of categories \(\cE_S\) of
      Definition~\ref{def:nor:coefficients};
\item the primitive contact sheaves and the selected summands
      \(\mathcal T_{v,u_v}\) of Definition~\ref{def:res:selected};
\item the correspondence ring spectrum \(\mathscr A_K^{\mathrm{sph}}\), the
      homomorphisms \(R_S\to\mathscr A_K^{\mathrm{sph}}\), the functors from
      the graph-cospan category, and the boundary functors
      of Section~\ref{sec:cor:correspondence};
\item the Schwartz space \(\cS(\A_K)\), its lattice tests from
      Definition~\ref{arith:def:lattice}, Fourier transform, theta maps,
      and the complex \(C_K^{\mathrm{FP}}\) of Definition~\ref{def:ana:FP},
      with the action of \(\mathscr A_K\);
\item the normal-coordinate realization maps
      \eqref{eq:ana:ball-module}--\eqref{eq:ana:projection}, the complete
      trace of Theorem~\ref{thm:ana:complete-trace}, and the normal
      determinant functions of Definition~\ref{def:ana:determinant};
\item the graded section algebras $\mathscr S_v$ and their finite-period
      specializations of Section~\ref{sec:products}, with coefficient
      cones, Hermitian metrics, quotient maps and additive transfers;
\item the marked projective coefficient realizations and their lattices
      $S_M(D)$, numerator identifications and jet maps, with the
      Bombieri and annular coefficient metrics of
      Sections~\ref{sec:section-metrics} and~\ref{sec:annular-sections};
\item the fixed-weight quotients $Q_\alpha$, for each unitary character
      component and $\alpha>1/2$, and the maps from the corresponding
      original Poisson quotient in Section~\ref{sec:zero-comparison}.
\end{enumerate}
The incidence maps and coefficient maps are those constructed in the preceding
sections. The finite-set structure maps use the same uniformizer at a place
whenever that place belongs to more than one set \(S\).
\end{definition}

The constituent operations fit together in the following form.

\begin{theorem}[Compatibility of the winding construction]
\label{thm:assembly:compatibility}
For the data of Definition~\ref{def:assembly:curve}, the following statements hold.
\begin{enumerate}[label=\textup{(\roman*)}]
\item The realization of \(\mathcal C_{\tau,g}\) in each valuation cover and
      in \(h\Cospan(\mathrm{Corr}_K)\) preserves sewing and the specified
      finite transfers. For the period \(L\), the associated permutation
      sheaf is \(\Loc_v(R_L)=\Sph_{Y_K}[Q_{v,L}]\).
\item The two normal coefficient lifts are related by the descended sequence
      \[
       \Loc_v(E)\xrightarrow{1-t_v}\Loc_v(E)
       \longrightarrow\Delta_{Y_K}
                 (\Sph_{\mathrm{aug}}\otimes_R E).
      \]
      The construction is compatible with the finite-set transition functors.
\item At the place \(v\), the chosen angular refinement realizes the
      primitive contact whose selected indexed fibre endomorphism ring is
      \(k_v\) in odd residue characteristic and \(W_2(k_v)\) in residue
      characteristic two.
\item For every finite set \(S\), extension to
      \(\mathscr A_K^{\mathrm{sph}}\) recovers \(R_S\)-modules and their
      mapping spectra through the Amitsur totalization
      \eqref{eq:cor:descent}. The same assertion holds at the indexed
      finite cyclic levels.
\item The transition functions of the normal cover are the valuation
      shifts in the Schwartz realization. The action of Laurent source
      operators and smoothing kernels has the complete pairing
      \eqref{eq:ana:source-pairing}, with its real and complex terms.
\item Tensor degree gives a descended section product on $\mathscr S_v$.
      At finite period the coefficient cone is proper and the displayed
      Hermitian norm makes multiplication contractive. The quotient
      maps are contractions; the normalized transfers are isometric
      multiplicative embeddings with the specified idempotent unit image.
\item On each marked projective coefficient realization, the integral-jet
      section lattices multiply into $S_M(D+E)$, and every nonzero
      multiplier is injective. For $n(D)\ge0$,
      Proposition~\ref{prop:prod-dual} and
      Theorem~\ref{thm:prod-jet-theta} give the dual lattice and relative
      theta identity. When $n(D)=M-1$,
      Theorem~\ref{thm:prod-index-expansion} gives the index.
\item In the annular coefficient completion, division by a cyclotomic
      polynomial on its vanishing ideal is bounded. The finite-period
      quotient metrics yield the deck-invariant coefficient metrics by
      averaging and normalization. For any finite set of $N$ distinct zero jets
      $\ell_{\rho,j}$ with $0\le j<m_\rho$, the evaluation map
      $Q_\alpha\to\C^N$ is continuous and surjective.
\end{enumerate}
\end{theorem}

\begin{proof}
The graph realization and its cofiber sequence are
Proposition~\ref{prop:src:sewing}. Associated descent on \(Q_v\) gives the
valuation realization of Definition~\ref{def:geo:realization}; the construction
on finite torsors commutes with the quotient and transfer maps of
Definition~\ref{def:src:periods}. Proposition~\ref{prop:cor:cospans} gives the
correspondence realization. The identification of the period sheaf and the
normal sequence follow from Theorem~\ref{thm:nor:descent}, and
Proposition~\ref{prop:nor:coefficients} gives the compatibility in the place set.

The cyclic refinement of Proposition~\ref{prop:res:cyclic-refinement} and the
selected sheaf of Proposition~\ref{prop:res:selected} give the indexed
endomorphism assertion. The Amitsur statement is
Theorem~\ref{thm:cor:descent}, including its finite-level base changes.
Propositions~\ref{prop:ana:normal-quotient} and~\ref{prop:ana:ball-module}
identify the normal and Schwartz transition functions. Finally,
Proposition~\ref{prop:ana:source-trace} computes the composed source action by
convolution, and Theorem~\ref{thm:ana:complete-trace} evaluates its trace.
Propositions~\ref{prop:prod-descent}, \ref{prop:prod-finite-positive},
and~\ref{prop:prod-period-metric} prove the graded finite-period assertions.
Theorem~\ref{thm:prod-integral-jets} gives the section product with the jet
condition, and its following index and duality statements give the stated
lattice identities. Theorems~\ref{thm:pos-jet-quotient}
and~\ref{thm:zero:jets} prove the two quotient assertions.
Proposition~\ref{prop:pos-finite-cover-metric} computes the finite-cover
metric.
\end{proof}

The section metrics and theta identities concern the specified coefficient
lattices. The normalized relative-determinant pairing of
Theorem~\ref{thm:st:comparison} recovers the complete arithmetic trace.
For tests on the norm quotient, \eqref{eq:st:combined} expresses it as a signed
sum of contact theta differences and regularized Green intersections.
Positivity on the zero representation, and a comparison with the positive
section metrics sufficient to prove it, remain open.
A duality comparison with oppositely twisted sections must incorporate
the finite quotient computed in Proposition~\ref{prop:prod-jet-duality-quotient}.
\section{Sheaf cohomology of the arithmetic winding curve}
\label{sec:sheaf-comparison}
For the space $Y_K$ and valuation covers of Section~\ref{sec:geometry},
we construct a sheaf $\mathscr F_K$ that compares the middle
Fourier--Poisson representation with the cohomology of the curve.
The finite charts use the valuation action on tests. At each
archimedean boundary the attaching map is the adelic Poisson
summation map of Proposition~\ref{prop:ana:intersection}.
The comparison of complex $C_K$-modules takes the form
\[
 H^1(C_K^{\mathrm{FP}})\ \xrightarrow{\ \iota\ }\
 H^1(Y_K,\mathscr F_K)\ \xrightarrow{\ r_{v,w}\ }\
 H^1(C_K^{\mathrm{FP}})^{C_K^1}.
\]
Their composite is the average of the coefficient action of $C_K^1$.
The first map uses all the charts. The second uses a finite place
$v$ and an archimedean place $w$.

\subsection{The middle Fourier--Poisson term in sheaf cohomology}
\label{cmp:sheaf}
Put $C=C_K$, $C^1=\ker|\cdot|$, $E=H_-$, $M=H_\cap$,
$P=\theta M$ and $Q=E/P$. Thus $[M\xrightarrow{\theta}E]$, in
degrees zero and one, is the middle summand of the Fourier--Poisson
complex of Definition~\ref{def:ana:FP}. In this comparison the coefficients on the ordinary
open-set site are complex vector spaces. Write $\Delta A$ for the
constant sheaf with fibre $A$ and $\operatorname{LC}(S,A)$ for
functions $S\to A$ constant on a neighbourhood of each point.

By Definition~\ref{def:geo:atlas} there is an open cover $Y_K=U_f\cup U_a$, where
$U_f=\bigcup_{v\nmid\infty}N_v$ and $U_a=\bigcup_{w\mid\infty}N_w$.
Their intersection is $C$. On a finite chart $N_v$ take the sheaf
of Definition~\ref{def:geo:realization} for the valuation cover with fibre $E$ and deck generator
$T_v=\lambda_{j_v(\varpi_v)}$. Its puncture frame is \eqref{eq:geo:cover-frame}.
On an infinite chart, for $j_w:C\hookrightarrow N_w$, use the
pushout of sheaves
\[
  \Delta M\ \amalg_{\,j_{w!}\Delta M}\ j_{w!}\Delta E .
\]
The arrows are the counit of restriction and extension by zero,
and the map induced by $\theta$, respectively. Each chart now
restricts to $\Delta E$ on $C$. Glue by the identity to obtain
$\mathscr F_K$ on $Y_K$. For $a\in C$ the comparison uses translation
by $a$ on the base together with its given action on the coefficient
spaces. This is well defined because that action commutes with each
$T_v$ and intertwines $\theta$.

\begin{proposition}\label{cmp:mv}
There is a morphism $[M\xrightarrow{\theta}E]\to
R\Gamma(Y_K,\mathscr F_K)$ in the derived category of complex vector
spaces. In degree one it induces a $C$-equivariant injection
\[
             \iota:Q\lhook\joinrel\longrightarrow H^1(Y_K,\mathscr F_K).
\]
Let $G=C^1/(C^1)^\circ$ and let $D_a\subset G$ be generated by the
images of $j_w(-1)$ at the real places. The kernel of restriction
to $H^1(U_f,\mathscr F_K)\oplus H^1(U_a,\mathscr F_K)$ identifies with
\[
 \operatorname{LC}(G,E)\big/\operatorname{LC}(G/D_a,P).
\]
In this quotient $\iota$ takes a class to its constant function.
Normalized Haar mean on $G$, followed by reduction modulo $P$,
is a left inverse there. Its constant classes form the invariant
subspace for translations by $C^1$ on the argument alone, with the
coefficient fibre fixed.
\end{proposition}

\begin{proof}
Fix a continuous splitting $\sigma:\mathbb R_{>0}\to C$ of the
norm. Every section of $\Delta E$ on $C$ factors through
$c\mapsto[c\sigma(|c|)^{-1}]\in G$ and has finite range, by
connectedness of the positive ray and compactness of $C^1$.
Suppose that such a section $f$ extends to a finite chart. Put
$h=j_v(\varpi_v)$. In the valuation cover the sequence
$[c,\varpi_v^n]_{\mathcal O_v^\times}$ tends to $[c,0]_{\mathcal O_v^\times}$.
In this lift a section has one constant value $m\in E$ for all
sufficiently large $n$. The puncture frame gives
$f(ch^n)=T_v^{-n}m$. Two of these values agree by finiteness,
so a positive power of $T_v$ fixes $m$. Iterating this power moves the norm to an end, where every
$H_-$ test tends to zero, so $m=0$. Thus $f(ch^n)=0$ for all large $n$.
The element $h\sigma(|h|)^{-1}$ of the compact group $C^1$
has arbitrarily large positive powers approaching the identity.
Along such powers $f(ch^n)=f(c)$ eventually. Hence the restriction
from $U_f$ is zero, and a finite-chart section vanishing on the
dense puncture is itself zero.

At $[c,0]\in N_w$ for $w\mid\infty$, the pushout stalk is $M$,
so a section agrees near that point with the image of a constant
$m\in M$. The paths $[c,r]$ and $[c,ru]$, with $r>0$ tending to
zero and $u$ any archimedean unit, give
$f(c)=f(cj_w(u))=\theta m$. The only angular components that
can act nontrivially on $G$ are the real signs. Conversely every
element of $\operatorname{LC}(G/D_a,P)$ extends to $U_a$: lift
its finitely many values to $M$, and on $N_w$ evaluate the lifted
function at the image of $c\sigma(|c|)^{-1}$ in $G/D_a$.
This is independent of the representative $[c,z]$, since the
image of $j_w(K_w^\times)$ is trivial in $G/D_a$. By Proposition~\ref{prop:ana:intersection}, $\theta$ is injective. Thus a section vanishing on the puncture has zero
boundary stalks as well.

Use finite restriction minus archimedean restriction for the
cover differential. Under the identifications just obtained the
two terms and their differential are
\[
 [\,\operatorname{LC}(G/D_a,P)\xrightarrow{-\,\mathrm{incl}}
                     \operatorname{LC}(G,E)\,].
\]
Send $m\in M$ to the archimedean section with constant value $-m$
(its puncture restriction is $-\theta m$), and send $e\in E$ to
the constant function with value $e$. These maps form a chain map
from $[M\xrightarrow{\theta}E]$. The maps from degree-zero sections
of each open and of their intersection to the corresponding derived
sections induce the map from the cover complex to $R\Gamma(Y_K,\mathscr F_K)$.
The Mayer--Vietoris triangle identifies its degree-one map with the
connecting injection of the displayed cokernel
\cite[Tags~01EB and 08BX]{Stacks}.

A function in $\operatorname{LC}(G,E)$ factors through a finite
quotient of the profinite group $G$, by compactness. Its Haar
mean is a finite sum and takes the denominator into $P$.
If the class of $f$ is fixed by base translations, every translate
differs from $f$ by an element of $\operatorname{LC}(G/D_a,P)$.
Averaging the finitely many translates shows that the constant
mean represents the same class. Conversely every constant class
is fixed. This also proves the claimed splitting.
\end{proof}

The finite winding action is needed for this injection. If all
$T_v$ in the finite-chart coefficient systems are replaced by the
identity while the same puncture frames are used, every constant
$e\in E$ extends to $U_f$. Its connecting class then vanishes.

\begin{proposition}\label{cmp:cup}
For $x,y\in Q$ the cup product of their images vanishes already in
$H^2(Y_K,\mathscr F_K\otimes_{\mathbb C}\mathscr F_K)$. The same
vanishing holds for $\iota(x)\smile\overline{\iota(y)}$ with
coefficients $\mathscr F_K\otimes_{\mathbb C}\overline{\mathscr F_K}$.
\end{proposition}

\begin{proof}
Both images restrict to zero on both opens by
Proposition~\ref{cmp:mv}. Set $Z_f=Y_K\setminus U_f$ and
$Z_a=Y_K\setminus U_a$. For sheaves $\mathcal A,\mathcal B$ the
support triangle lifts a class vanishing on $U_f$ to the
cohomology of the fibre of $\mathcal A\to Rj_{f*}(\mathcal A|_{U_f})$,
which is supported on $Z_f$. Here $j_f:U_f\hookrightarrow Y_K$ is the
open inclusion. A class vanishing on $U_a$ lifts in the same way to a
complex supported on $Z_a$ \cite[Tag~0G72]{Stacks}.
The derived tensor of the two supported complexes has zero
stalks because $Z_f\cap Z_a=\varnothing$. Tensoring sheaves of complex
vector spaces is exact on stalks. Naturality of cup
product sends the product of the two lifts to the product of the
original classes \cite[Tag~0FKU]{Stacks}. This proves the
first assertion. Conjugating coefficients commutes with
restriction, so the same argument gives the second.
\end{proof}

Consequently any pairing of these images obtained by a morphism
from the indicated tensor of sheaves and a subsequent trace is zero.
A dual class outside the restriction kernel need not have this
vanishing property.

\subsection{Periodization of the spherical Poisson image}
\label{cmp:periodization}
Put $P=\theta H_\cap$, $Q=H_-/P$, and
\[
 E_{\mathrm D}=H_-^{C_K^1},\qquad P_{\mathrm D}=P\cap E_{\mathrm D}=\theta(H_\cap^{C_K^1}),
 \qquad Q_{\mathrm D}=E_{\mathrm D}/P_{\mathrm D}\simeq Q^{C_K^1}.
\]
Here $C_K^1$ acts on $H_+$ after taking $K^\times$-coinvariants.
Its Haar projection commutes with $\theta$, and averaging a lift of
an invariant quotient class proves the last identification.
Proposition~\ref{prop:ana:intersection} gives a closed theta image on
the stage fixed by the finite integral units. Its intersection with
$E_{\mathrm D}$ is $P_{\mathrm D}$, which is therefore closed.
In the uncentered coordinate $t=\log|c|$, put $T_uf(t)=f(t+u)$ for
$u\in\R$. Thus $T_{\log q_v}$ is the action of $j_v(\varpi_v)$ on $E_{\mathrm D}$.
The space $P_{\mathrm D}$ is a convolution ideal. Indeed, integrals of compactly
supported smooth kernels against translates of an element of $P_{\mathrm D}$
belong to $P_{\mathrm D}$ by closedness. Cutting off a kernel in $E_{\mathrm D}$ gives
convergence of the convolutions in every exponentially weighted
Schwartz seminorm.

\begin{proposition}\label{cmp:nonresonance}
For every $h>0$, periodization $\Pi_h f(t)=\sum_{k\in\Z}f(t+kh)$
maps $P_{\mathrm D}$ onto $C^\infty(\R/h\Z)$. The induced operator
$1-T_h$ on $Q_{\mathrm D}$ is a continuous linear automorphism with continuous inverse.
\end{proposition}
\begin{proof}
First choose $\Phi\in\mathcal S(\mathbf A_K)$ with $\Phi(0)=1$ and
$\int\Phi=0$, as in Proposition~\ref{prop:ana:intersection}.
Let $F(t)$ be the norm-one average of $\Sigma\Phi$ at norm $e^t$,
using probability Haar measure. The weighted bounds in \eqref{eq:ana:embeddings} and the Poisson identity
\eqref{eq:ana:poisson} give $F(+\infty)=0$ and $F(-\infty)=-1$.
The difference $F-T_hF$ lies in $P_{\mathrm D}$ because in $H_+$ its source is
the $C_K^1$-average of $[\Phi]-\lambda_b[\Phi]$ for any $b\in C_K$
with $|b|=e^{-h}$, and both moments of this difference vanish.
Telescoping gives $\Pi_h(F-T_hF)=-1$.

For $a\in C_c^\infty(\R)$ put $M_a(s)=\int e^{st}a(t)\,dt$.
There is a radial Poisson test $p_a\in P_{\mathrm D}$ satisfying
\[
              M_{p_a}(s)=c_K(s-1)\zeta_K(s)M_a(s)
              \qquad(c_K>0).
\]
At the finite places take $\mathbf 1_{\prod_{v\nmid\infty}\mathcal O_v}$. On
$K_\infty^\times$ split log norm $t$ from its kernel and take
$(-a'(t)-a(t))b(y)$, where $b$ is a compactly supported smooth
function of integral one on the kernel. This has compact support
inside $K_\infty^\times$ and extends smoothly by zero. Its value at
zero is zero, and its additive integral vanishes since
$\int e^t(-a'-a)\,dt=0$. Apply $\Sigma$ and average over $C_K^1$.
Unfolding the integral for $\Re s>1$ gives the displayed identity.
With multiplicative measure giving $\mathcal O_v^\times$ mass one,
the integral over $\mathcal O_v\setminus\{0\}$ is $\sum_{n\ge0}q_v^{-ns}$.
The integral on $K_\infty^\times$ is $(s-1)M_a(s)$, and $c_K$
accounts for the quotient measure. The decay of $p_a$ makes the left
side entire, and meromorphic continuation of $\zeta_K$
\cite[\S5.11, Theorem~5.22(b)]{Poonen} extends the identity to $\C$.

Set $s_n=2\pi in/h$. For $n\ne0$ one has
$|\zeta_K(s_n)|^{-1}\ll_{K,h}(1+|n|)^B$ for some $B$.
For $|u|\ge2$, the convexity bound in the strip
$-1/2\le\Re s\le3/2$
\cite[Proposition~3.8]{HSWzeros}, absolute convergence to its
right, and Cauchy's formula give
$|\zeta_K(\sigma+iu)|+|\zeta'_K(\sigma+iu)|\ll_K(1+|u|)^b$
for $1\le\sigma\le2$.
Euler positivity from $3+4\cos x+\cos(2x)\ge0$ gives
\[
 \zeta_K(1+\delta)^3|\zeta_K(1+\delta+iu)|^4
                      |\zeta_K(1+\delta+2iu)|\ge1 .
\]
Thus $|\zeta_K(1+\delta+iu)|\gg_K\delta^{3/4}|u|^{-b/4}$.
Taking $\delta=|u|^{-5b-4}$ absorbs the displacement
$O_K(\delta |u|^b)$ to $\sigma=1$ for large $|u|$.
The same inequality as $\delta\downarrow0$ excludes a zero
at $1+iu$ for each fixed $u\ne0$. Cancelling the symmetric factor $s(s-1)$ in
\cite[\S2.1, (2.1)--(2.2), p.~279]{HSWzeros} gives
$\mathcal Z_K(s)=\mathcal Z_K(1-s)$ with the normalization of
Proposition~\ref{prop:ana:completed-product}. Stirling's formula gives,
for $|u|\ge2$,
\[
 |\zeta_K(iu)|\asymp_K |u|^{[K:\Q]/2}|\zeta_K(1-iu)|.
\]
This gives the asserted polynomial bound for large $|n|$. At the
remaining nonzero $s_n$, the functional equation transfers
nonvanishing since the gamma factors are finite and nonzero.

Use circle coefficients $\int_0^h e^{s_nt}b(t)\,dt$. The coefficients of $\Pi_h f$ are $M_f(s_n)$.
For prescribed mean-zero smooth circle data $(d_n)$ put $b_0=0$
and $b_n=d_n/(c_K(s_n-1)\zeta_K(s_n))$ for $n\ne0$. These
coefficients are rapidly decreasing. Choose the smooth period-$h$
function with unnormalized coefficients $b_n$ and multiply it
by a fixed compactly supported smooth function whose $h$-translates
sum to one. The resulting $a$ has $M_a(s_n)=b_n$ for every $n$,
so $\Pi_h p_a$ supplies the prescribed data, including zero mean.
The preceding telescoping test supplies every constant.

On $E_{\mathrm D}$, $\ker\Pi_h=(1-T_h)E_{\mathrm D}$. If $\Pi_h f=0$, the two sums
\[
 g(t)=\sum_{k\ge0}f(t+kh)=-\sum_{k<0}f(t+kh)
\]
give rapid exponential bounds for each derivative on the positive
and negative half-lines, respectively, and $(1-T_h)g=f$.
Surjectivity of $\Pi_h|_{P_{\mathrm D}}$ now gives surjectivity of
$1-T_h$ on $Q_{\mathrm D}$. If $T_hx=x$ there, the orbit $u\mapsto T_ux$
is $h$-periodic.
Every smooth circle kernel is a periodization of an element of
$P_{\mathrm D}$, and therefore annihilates this orbit by convolution.
Integration against $T_{-u}x$ gives the quotient of the ordinary
line convolution. A circle approximate identity converges to
$x$ by continuity, hence $x=0$. Finally $Q_{\mathrm D}$ is Fr\'echet, so the
open mapping theorem makes the inverse continuous.
\end{proof}

\subsection{The retraction on first cohomology}
\label{cmp:retraction}
Write $E=H_-$, $M=H_\cap$, $P=\theta M$ and $Q=E/P$.
On the overlap $C_K$, constant sections give the connecting map
$\iota:Q\to H^1(Y_K,\mathscr F_K)$ of Proposition~\ref{cmp:mv}. Take finite restriction minus archimedean
restriction in the Mayer--Vietoris complex, with positive overlap
coordinate for the connecting map. We take sheaf cohomology in
complex vector spaces.

\begin{proposition}\label{cmp:prop-retraction}
Fix a finite place $v$ and an infinite place $w$. Put
$V=\prod_{u\nmid\infty}\mathcal O_u^\times\subset C_K$,
$G_V=C_K/V$, and let $d$ be the order of the component of
$h=j_v(\varpi_v)\bmod V$ in $\pi_0(G_V)$. Put $T=T_v^d$.
There is a map commuting with the coefficient action,
\begin{equation*}\tag{R}\label{cmp:eq-pointed-map}
 B_{v,w}:H^1(Y_K,\mathscr F_K)\longrightarrow E/(T-1)P,
 \qquad B_{v,w}\iota([e])=[-(T-1)e].
\end{equation*}
Let $\operatorname{pr}:E/(T-1)P\to Q$ be the quotient map and
$\Pi_{\mathrm D}:Q\to Q_{\mathrm D}$ the probability Haar projector onto the
$C_K^1$-invariant \emph{coefficient} component of
Proposition~\ref{cmp:nonresonance}. One may choose $B_{v,w}$ so that
\[
 r_{v,w}=-(T-1)^{-1}\Pi_{\mathrm D}\operatorname{pr} B_{v,w}
       : H^1(Y_K,\mathscr F_K)\longrightarrow Q_{\mathrm D},
 \qquad r_{v,w}\iota=\Pi_{\mathrm D}.
\]
In particular $r_{v,w}$ retracts $\iota|_{Q_{\mathrm D}}$ and is equivariant
for simultaneous translation on the base and the coefficients.
\end{proposition}

\begin{proof}
Embed $V$ by taking archimedean components equal to one. Write
$q:Y_K\to Y'=V\backslash Y_K$. Since $j_v^{-1}(V)=\mathcal O_v^\times$,
the finite chart downstairs and its winding cover are
\[
  N'_v=(G_V\times D)/\mathbb Z,\qquad \widetilde N'_v=G_V\times D,
  \qquad (c,n)\cdot m=(ch^m,n-m).
\]
Here $D=K_v/\mathcal O_v^\times=\mathbb Z\cup\{+\infty\}$ by valuation,
with positive tails a basis at the fixed endpoint $+\infty$.
Logarithmic norm provides covering trivializations of width less
than $\log q_v$. The original cover maps equivariantly to its
pullback and induces a bijection on torsor fibres. On the other
charts both covers use the given puncture trivialization. Inverse
image preserves constant sheaves, extensions by zero across open
inclusions, and pushouts, so the finite and archimedean descriptions
give $\mathscr F_K=q^{-1}\mathscr F'_K$ for the corresponding sheaf on $Y'$.

For a sheaf $A$ on $Y'$, the unit $A\to Rq_*q^{-1}A$ has a natural
retraction. Each chart $N_u$ is Hausdorff and $V$-invariant. On that chart the
orbit projection onto $V\backslash N_u$ is closed with compact
fibres, hence proper. Its fibre at the orbit of $y$ is the profinite
space $V/\operatorname{Stab}_V(y)$. Proper base change
\cite[Tag~09V6]{Stacks} identifies the stalks of
$R^iq_*q^{-1}A$ with cohomology of this fibre with constant value
$A_{q(y)}$. Disjoint clopen refinement kills the groups for $i>0$.
The quotient charts cover $Y'$.
Average sections on each orbit with probability Haar measure on $V$,
taking values in the single stalk downstairs. To verify the sheaf
condition for this map, represent a section near each $gy$ by a
pullback from downstairs. Take a finite clopen partition
$V=\coprod V_i$ with each $V_i y$ in a neighbourhood having one
representative $a_i$. By compactness the containments persist on
a neighbourhood $W$ of $y$. On the open set $q(W)$ the mean is
$\sum_i\mu(V_i)a_i$. It fixes pullbacks and commutes with restriction,
so it carries connecting classes of constants from $Y_K$ to $Y'$.

The valuation map from $G_V$ to the finite ideal class group has open
kernel $K_\infty^\times/\mathcal O_K^\times$, where the units are
closed and discrete by Dirichlet's theorem. Hence $G_V$ is an abelian
Lie group with finitely many components. Its norm-one subgroup is compact, and
since $\log|h|=-\log q_v$, $G_V/\langle h\rangle$ is compact. Write
$O=(G_V/\langle h\rangle)^\circ$, a torus, and let $N_e$ be its
inverse image under the first-coordinate projection, with
$p:N_e\to O$. The winding coefficient on $N_e$ is $p^{-1}L$ for
the coefficient system $L$ associated to $G_V\to G_V/\langle h\rangle$ on
$O$. Its monodromy on $\Lambda=\pi_1(O,[1])$ is $T_v^{\nu(\cdot)}$,
where $\nu(\Lambda)=d\mathbb Z$. Indeed, a lift starting at $1$ can end at
$h^j$ if and only if $h^j\in G_V^\circ$. A lifted loop ending at
$h^j$ transports a vector $a$ to $T_v^ja$, since the associated
bundle identifies $(p+j,a)$ with $(p,T_v^ja)$ in
Definition~\ref{def:geo:realization}. This determines the sign of
$\nu$.

We claim $p^*:R\Gamma(O,L)\to R\Gamma(N_e,p^{-1}L)$ is an
isomorphism. Over a contractible coordinate box $U\subset O$ the
map is $U\times D\to U$ with constant coefficient. Split $D$ as
the compact clopen tail $D_+=\{0,1,\ldots,+\infty\}$ and the
discrete negative levels. Proper base change for $U\times D_+\to U$
kills the positive direct images. Its zeroth image is the constant
sheaf with fibre $\operatorname{LC}(D_+,E)$, as seen
on connected subboxes. Each negative level is a copy of $U$, and
cohomology on their disjoint union is a product, exact in vector
spaces. Thus $R^ip_*p^{-1}L=0$ for $i>0$ and the degree-zero direct image has fibre the space of sequences indexed by $\Z$ with values in $E$, eventually constant at the positive end. Endpoint evaluation splits off $L$ by constant
sequences. On its kernel, a primitive $\gamma\in\Lambda$ of
degree $d$ acts by an index shift of $d$ and an invertible
coefficient. For the forward shift the equation $a_n-U_0a_{n+d}=b_n$ has the
unique zero-tail solution $\sum_{j\ge0}U_0^jb_{n+jd}$ at each $n$.
Reversing the shift gives the same assertion. Complete $\gamma$ to
a lattice basis. The universal-cover cellular chains of the torus
are the tensor product of circle resolutions. On the kernel their
cochains are contractible since one commuting difference is
invertible. This proves the claim, compatibly with evaluation at
the point $1\in C_K/V$ over $[1]$. The group-cohomology model with these coefficients is described in \cite[\S3, following (3.4)]{BrownGroups}.

In the archimedean chart use the arc $[0,1]\to N'_w$ given by
$r\mapsto[1,r]$ modulo $V$. The endpoint $1$ lies over the same point. The pullback coefficient is
the pushout of $\Delta M$ and
$j_!\Delta E$ over $j_!\Delta M$ for $j:(0,1]\hookrightarrow[0,1]$.
The pushout is derived as well since $j_!\Delta M\to\Delta M$
is a monomorphism. For every vector space $B$,
$R\Gamma([0,1],j_!\Delta B)=0$ by the relative cohomology of
$([0,1],\{0\})$. Thus the arc has derived sections $M$ and endpoint
map $\theta$.

Restrict the two-open Mayer--Vietoris diagram on $Y'$ to $N_e$,
the arc at $w$, and the common point $1$. The evaluation squares
commute. The target homotopy fibre is therefore that of
$R\Gamma(O,L)\oplus M\to E$, with the difference of evaluations.
In normalized group cochains of the pointed torus, its degree-one
classes are pairs $(u,b)\in Z^1(\Lambda,E)\oplus E$ modulo
\[
 ((T_v^{\nu(\cdot)}-1)a,\ a-\theta m),\qquad a\in E,\quad m\in M.
\]
If $\nu(\alpha)=0$, the cocycle identity on $\alpha\gamma=\gamma\alpha$
gives $(T-1)u(\alpha)=0$. A $T$-periodic element of $H_-$ vanishes
by uniform decay as the norm tends to an end. Consequently $u$
factors through $d\mathbb Z$ and is uniquely determined by its
unrestricted value $u(\gamma)$. The first cohomology of the fibre
identifies with $E/(T-1)P$ by
$[(u,b)]\mapsto[u(\gamma)-(T-1)b]$, independent of the loop of
degree $d$. A connecting class is $(0,e)$, which gives
\eqref{cmp:eq-pointed-map}.

On $Q_{\mathrm D}$ the operator $T$ is $T_{d\log q_v}$, so
Proposition~\ref{cmp:nonresonance} supplies the inverse of $T-1$
appearing in $r_{v,w}$. Since $P$ is invariant under $T$,
$\operatorname{pr}$ is defined, and $\Pi_{\mathrm D}$ commutes with
$T$. Formula~\eqref{cmp:eq-pointed-map} now gives
$r_{v,w}\iota=\Pi_{\mathrm D}$.

All maps so far commute with coefficient translations. With coefficients
fixed, $G_V^\circ$ acts trivially on $H^1(Y',\mathscr F'_K)$.
Along a path of base translations the lifts to the winding
covers identify the two pullbacks of the coefficient systems on $[0,1]\times Y'$.
The archimedean pushouts have the same identification since the
action preserves the puncture. On each quotient chart the interval projection is proper
with contractible fibres. The base-change theorem above makes the unit on
pulled-back coefficients an isomorphism, so both endpoint restrictions invert
the same pullback. Average the pointed map on $Y'$ over the finite group
$G_V/G_V^\circ$ acting on the base alone, and compose with the $V$-average
above. Base translations fix connecting classes of constant functions, so
\eqref{cmp:eq-pointed-map} still holds. The resulting map is invariant under
base translations and commutes with coefficient translations.
\end{proof}

\section{A section formula for the complete arithmetic trace}
\label{sec:section-trace-comparison}

We retain the test space, measures, and conjugate-linear involution of
Sections~\ref{sec:notation} and~\ref{sec:ana:trace}. For
$f,g\in\cD(C)$ put
\[
 \Phi=f*g^\sharp,\qquad \varphi_v=\Phi\circ j_v.
\]
The construction below uses relative determinants of multiplication maps
on linear coefficient sections and an extension at the identity.
On the norm quotient it also has an expression in the contact theta
heights of Theorem~\ref{thm:thetaRR} and metric intersections on the marked
coefficient line. The local graph--diagonal interpretation of the Weil
kernel has a precedent in \cite[\S7, equation~(7.5)]{CCMWeil}.

\subsection{Relative determinants and the identity contribution}

The numerator lattice $S_M([\infty])$ is $\Z+\Z t$ for every $M$:
an integral linear polynomial has integral jets at $t=1$.
Homogenize it with basis $X,Y$, so $t=X/Y$, and extend scalars to $K_v$.
For $a\in K_v^\times\setminus\{1\}$ write $b=a^{-1}=u/w$, where $w\ne0$;
at a finite place choose $u,w\in\mathcal O_v$. Put
\[
 D_b=wX-uY,\qquad D_1=X-Y,\qquad D_\infty=Y,
 \qquad
 A_j(s,t)=sD_j+tD_b\quad(j=1,\infty).
\]
These are sums of two section multiplication maps from degree zero to
degree one. With the same source and target determinant frames define
\[
 \mathcal J_v(a)=\frac{|\det A_\infty|_v}{|\det A_1|_v}.
\]

\begin{lemma}[Relative section determinant]\label{lem:st:determinant}
The ratio is independent of the representation $a^{-1}=u/w$, and
\[
 \mathcal J_v(a)=|1-a^{-1}|_v^{-1}.
\]
At a finite place it is
$\#\coker(A_1)/\#\coker(A_\infty)$ for the displayed maps of
$\mathcal O_v$-lattices. At an infinite place it is a ratio of
determinant norms.
\end{lemma}
\begin{proof}
The two determinants in the basis $X,Y$ are $w-u$ and $-w$.
Their ratio gives the formula, and common rescaling cancels.
The lattice statement is the determinant formula for the index of an
injective square matrix over a discrete valuation ring.
At a complex place the normalized absolute value is the square of
the usual complex modulus, as in
Theorem~\ref{thm:ana:fixed}.
\end{proof}

For the reference additive characters of Definition~\ref{def:ana:weil},
let
\[
 G_\R(x)=e^{-\pi x^2},\qquad
 G_\C(z)=e^{-2\pi|z|_{\rm eucl}^2}.
\]
Both have integral one and are self-dual for the reference measures
$dx$ and $2\,dx\,dy$. Define
\[
 c_v=-\int_{K_v}\log|x|_v\,G_v(x)\,dx_{v,0}.
\]
Differentiating the Gaussian Mellin integrals of
Proposition~\ref{prop:ana:gaussians} gives
\[
 c_\R=\frac{\gamma+\log(4\pi)}2,\qquad
 c_\C=\gamma+\log(2\pi),
\]
where $\gamma$ is Euler's constant. Use the Weil multiplicative measures
of Definition~\ref{def:ana:weil} and set
\begin{align}
 I_v(\Phi)
 &=\int_{K_v^\times}
   \bigl(\varphi_v(a)-\Phi(1)1_{\mathcal O_v^\times}(a)\bigr)
      \mathcal J_v(a)\,d_W^\times a
       &&(v\nmid\infty),\label{eq:st:finite}\\
 I_v(\Phi)
 &=\int_{K_v^\times}
   \bigl(\varphi_v(a)-\Phi(1)G_v(1-a)\bigr)
      \mathcal J_v(a)\,d_W^\times a+c_v\Phi(1)
       &&(v\mid\infty).\label{eq:st:infinite}
\end{align}
Give $\mathcal O_K$ the Minkowski metric
$\sum_{\sigma:K\hookrightarrow\R}|\sigma(x)|^2+
2\sum_\tau|\tau(x)|^2$ on $K\otimes_\Q\R$, where $\tau$ runs over
one embedding from each conjugate complex pair, and write
$\mathcal O_K^\#$ for its Euclidean dual.
Define the sesquilinear form
\begin{equation}
 I_K(f,g)=\sum_v I_v(\Phi)
 -2\bigl(h_\theta(\mathcal O_K^\#)-h_\theta(\mathcal O_K)\bigr)\Phi(1).
 \label{eq:st:global}
\end{equation}

\begin{theorem}[Section--trace comparison]\label{thm:st:comparison}
For every number field $K$ and all $f,g\in\cD(C)$, the expressions above
are well defined and
\begin{equation}
 I_K(f,g)=W_K(f*g^\sharp)
 =\operatorname{str}_{H^*(C_K^{\mathrm{FP}})}\rho(f*g^\sharp).
 \label{eq:st:comparison}
\end{equation}
This includes arbitrary angular dependence and convolutions nonzero at
the identity. Replacing $f,g$ by $f_P,g_Q$ gives the corresponding formula
for Laurent source operators in Proposition~\ref{prop:ana:source-trace}.
If the two pole moments of each test vanish, then
$I_K(f,g)=-\Tr(\rho(f*g^\sharp)\mid H^1)$.
\end{theorem}
\begin{proof}
Lemma~\ref{lem:st:determinant} identifies the relative section
determinant with the inverse normal Jacobian in
Theorem~\ref{thm:ana:fixed}. Thus \eqref{eq:st:finite} is
\eqref{eq:ana:finite-PV}. Its integrand vanishes near $1$.
For all but finitely many finite places,
$\varphi_v=\Phi(1)1_{\mathcal O_v^\times}$, by the support and
unit-invariance argument following Definition~\ref{def:ana:weil}.

At infinity let $P_v=\mathcal F_{v,0}(\log|\cdot|_v)$, as in that
definition, and extend $\varphi_v$ by zero at $0$. Then
$\eta(x)=\varphi_v(1-x)$ is a Schwartz function and $\eta(0)=\Phi(1)$.
The distribution $P_v$ is $-d_W^\times x$ on tests vanishing at zero
\cite[Lemma~3.4]{Meyer}. Apply this to
$\eta-\eta(0)G_v$. Since $\widehat G_v=G_v$,
$-\langle P_v,G_v\rangle=c_v$, and hence
\[
 -\langle P_v,\eta\rangle
 =\int_{K_v^\times}\bigl(\eta(x)-\eta(0)G_v(x)\bigr)d_W^\times x
   +c_v\eta(0).
\]
Under $x=1-a$ the density $d_W^\times x$ becomes
$\mathcal J_v(a)d_W^\times a$. This proves
$I_v(\Phi)=W_v^0(\varphi_v)$. The subtracted function is
$O(|x|_{\rm eucl})$ at zero, so the integral is absolutely convergent
in both real dimensions one and two; its other ends are integrable
by compact support and Gaussian decay.

For a Euclidean lattice $L$, Poisson summation gives
$h_\theta(L^\#)-h_\theta(L)=\log\covol(L)$.
In the displayed Minkowski metric the discriminant determinant gives
$\covol(\mathcal O_K)=|D_K|^{1/2}$. Thus the last term of
\eqref{eq:st:global} is $-\log|D_K|\,\Phi(1)$.
Definition~\ref{def:ana:weil} and
Theorem~\ref{thm:ana:complete-trace} now give
\eqref{eq:st:comparison}. Finite source translates remain in $\cD(C)$,
and the moment-zero assertion is
Proposition~\ref{prop:ana:source-trace}.
\end{proof}

Adding $c\Phi(1)$ would preserve the off-identity kernel but change the form;
even on the moment-zero subspace,
\[
 (f*f^\sharp)(1)=\int_C|f(c)|^2|c|\,dc.
\]
The reference additive character and its self-dual Gaussian fix this
scalar. The normal factors give the same normalization.

\begin{proposition}[Normal-factor first variation]\label{prop:st:riesz}
Put $L_v(s)=\tau_v(s)$ at a finite place and
$L_v(s)=\Gamma_v(s)$ at an infinite place, in the normalization of
Section~\ref{sec:ana:determinants}. At finite places also put
$G_v=1_{\mathcal O_v}$ and use the reference additive measure of
$\mathcal O_v$-volume one. For $\varphi\in\cD(K_v^\times)$, extended
by zero at $0$, the expression
\[
 \mathcal R_{v,s}(\varphi)=
 \frac{L_v(1-s)}{L_v(s)}
 \int_{K_v}|x|_v^{s-1}\varphi(1-x)\,dx_{v,0}
 \qquad(\Re s>0,\ |s|\text{ small})
\]
extends holomorphically to a neighborhood of $s=0$ and satisfies
\[
 \mathcal R_{v,0}(\varphi)=\varphi(1),\qquad
 \left.\frac{d}{ds}\mathcal R_{v,s}(\varphi)\right|_{s=0}
       =W_v^0(\varphi).
\]
\end{proposition}
\begin{proof}
Write $d_W^\times x=\alpha_v\,dx_{v,0}/|x|_v$.
The normal factors give
\[
 \left.\frac{d}{ds}\frac{L_v(1-s)}{L_v(s)}\right|_{s=0}
 =\alpha_v=
 \begin{cases}
  (\log q_v)/(1-q_v^{-1}),&v\nmid\infty,\\
  1/2,&K_v=\R,\\
  1/(2\pi),&K_v=\C.
 \end{cases}
\]
Their constant term at $s=0$ is zero. Moreover
\[
 \int_{K_v}|x|_v^{s-1}G_v(x)\,dx_{v,0}
       =\frac{L_v(s)}{L_v(1)}.
\]
Subtract $\varphi(1)G_v(x)$ inside the integral. The difference
vanishes near zero at a finite place and vanishes to first order
in the Euclidean coordinate at infinity. Its Mellin integral is
holomorphic near zero. The removed term contributes
$\varphi(1)L_v(1-s)/L_v(1)$, proving the extension and
\[
 \mathcal R'_{v,0}(\varphi)
 =\int_{K_v^\times}
   \bigl(\varphi(1-x)-\varphi(1)G_v(x)\bigr)d_W^\times x
   -\frac{L_v'(1)}{L_v(1)}\varphi(1).
\]
At infinity $-L_v'(1)/L_v(1)=c_v$, so this is
\eqref{eq:st:infinite}. At a finite place the substitution $x=1-a$
changes $G_v(x)$ to $1_{\mathcal O_v}(a)$. Its excess over the unit
subtraction has mass
\[
 \int_{\varpi_v\mathcal O_v\setminus\{0\}}
       \mathcal J_v(a)d_W^\times a
 =\sum_{n\ge1}q_v^{-n}\log q_v
 =\frac{\log q_v}{q_v-1}
 =-\frac{L_v'(1)}{L_v(1)}.
\]
The last constant cancels this excess, giving
\eqref{eq:ana:finite-PV}.
\end{proof}

Set
$u_v(s)=\mathcal R_{v,s}(1_{\mathcal O_v^\times})$ at finite places
and $u_v(s)=1$ at infinite places. For each fixed $\Phi$ the germ
\[
 \mathcal A_s(\Phi)=|D_K|^{-s}\Phi(1)
 +\sum_v\bigl(\mathcal R_{v,s}(\varphi_v)-u_v(s)\Phi(1)\bigr)
\]
has only finitely many nonzero summands. The proposition gives
$u_v(0)=1$, $u_v'(0)=0$, and
$\mathcal A_0(\Phi)=\Phi(1)$,
$\mathcal A'_0(\Phi)=I_K(f,g)$.

\subsection{Finite places as contact theta differences}

Suppose now that $f,g$ are $C_K^1$-invariant. In the norm coordinate
write
\[
 u(t)=e^{t/2}f(e^t),\qquad w(t)=e^{t/2}g(e^t),\qquad
 \widetilde w(t)=\overline{w(-t)},\qquad H=u*\widetilde w.
\]
The notation $f(e^t)$ means the value on the norm fibre.
Then $u,w\in C_c^\infty(\R)$,
$\Phi(e^t)=e^{-t/2}H(t)$, and
$H(0)=\langle u,w\rangle_{L^2(\R)}$, with inner products linear
in the first variable.

For a rational prime $p$ and distinct positive integers $i,j$ let
\[
 L_{ij}=\Phi_{p^j}\bigl(\Z[t]/(\Phi_{p^i})\bigr),\qquad
 C_p(i,j)=
 \frac{h_\theta(L_{ij}^*)-h_\theta(L_{ij})}
      {\sqrt{\varphi(p^i)\varphi(p^j)}}.
\]
Use the orthonormal monomial basis on $\Z[t]/(\Phi_{p^i})$ and its
Euclidean dual. Theorem~\ref{thm:thetaRR} and
Corollary~\ref{cor:prime-weight} give
\begin{equation}
 C_p(i,j)=(\log p)p^{-|i-j|/2}.
 \label{eq:st:contact-coefficient}
\end{equation}

\begin{proposition}[Contact theta formula]\label{prop:st:finite-theta}
At a finite place over $p$ write $q_v=p^{d_v}$, where $d_v$ denotes
the residue degree in this subsection, and put $\ell_v=\log q_v$.
Choose an integer $A$ with $A+d_v n\ge1$ for every index $n$ meeting
the support of $u(s+n\ell_v)$ or $w(s+n\ell_v)$ for $0\le s<\ell_v$.
Then
\begin{equation}
 I_v(\Phi)=d_v\int_0^{\ell_v}
 \sum_{n\ne m} C_p(A+d_v n,A+d_v m)
          u(s+n\ell_v)\overline{w(s+m\ell_v)}\,ds .
 \label{eq:st:theta-pairing}
\end{equation}
The sum is restricted to indices meeting the respective two supports.
It is finite, its value is independent of $A$, and only finitely many
finite places contribute.
\end{proposition}
\begin{proof}
Set $r_v=q_v^{-1/2}$. The unit contribution in
\eqref{eq:st:finite} is zero. On the shells $|a|_v=q_v^{\pm k}$,
$k\ge1$, the centered weight is $r_v^k$, and each shell has Weil
mass $\ell_v$. Thus
\begin{equation}
 I_v(\Phi)=\ell_v\sum_{k\ge1}r_v^k
                     \bigl(H(k\ell_v)+H(-k\ell_v)\bigr).
 \label{eq:st:shell}
\end{equation}
Partition each correlation integral into translates of $[0,\ell_v)$
and use
$d_v C_p(A+d_v n,A+d_v m)=\ell_v r_v^{|n-m|}$.
This proves \eqref{eq:st:theta-pairing} and independence of $A$.
Compact support bounds the occurring indices. If $\ell_v$ exceeds
the support radius of $H$, every nonzero shell sample vanishes.
\end{proof}

There is also a formula on the valuation cover. With any fixed
additive Haar measure put
\[
 b_{v,n}=\frac{1_{\varpi_v^n\mathcal O_v}}
                  {\vol(\varpi_v^n\mathcal O_v)^{1/2}},\qquad
 (\mathcal B_vu)(s,x)=\sum_n u(s+n\ell_v)b_{v,n}(x).
\]
Nestedness of the balls gives
$\langle b_{v,n},b_{v,m}\rangle=r_v^{|n-m|}$, whence
\begin{equation}
 I_v(\Phi)=\ell_v\left(
 \langle\mathcal B_vu,\mathcal B_vw\rangle_{L^2([0,\ell_v)\times K_v)}
       -\langle u,w\rangle_{L^2(\R)}\right).
 \label{eq:st:ball}
\end{equation}

\subsection{Infinite places as regularized Green intersections}

Keep the centered norm profiles $u,w$ and $H=u*\widetilde w$.
On the marked coefficient line use the Fubini--Study Green kernel
\[
 G(z,z')=-\log|z-z'|+\tfrac12\log(1+|z|^2)
                         +\tfrac12\log(1+|z'|^2).
\]
In the coordinate $z=e^{t/2+i\theta}$, for $m\ge1$ define the smooth
compactly supported measure
\begin{align*}
 d\mu_{u,m}(z)&=2\sqrt m\,u(t)\cos(m\theta)\,dt\,\frac{d\theta}{2\pi},\\
 E_m(u,w)&=\iint G(z,z')\,d\mu_{u,m}(z)\,
                                  \overline{d\mu_{w,m}(z')},\\
 J_n(u,w)&=\frac n4 E_n(u,w).
\end{align*}
The measures have total mass zero. Their supports are compact annuli,
and the logarithmic diagonal singularity is integrable.
Set $\Psi=\Gamma'/\Gamma$ and
\[
 \begin{array}{c|c|c}
 K_v&\mathcal I_v&c_v^{\rm cen}\\ \hline
 \R&\{1,5,9,\ldots\}&\log\pi-\Psi(1/4)\\
 \C&\{1,3,5,\ldots\}&2\log(2\pi)-2\Psi(1/2).
 \end{array}
\]

\begin{theorem}[Green intersection formula]\label{thm:st:green}
For an infinite place and the radial test $\varphi_v(a)
=|a|_v^{-1/2}H(\log|a|_v)$,
\begin{equation}
 W_v^0(\varphi_v)=
 \lim_{N\to\infty}\left\{
  \sum_{\substack{m\in\mathcal I_v\\m\le N}}E_m(u,w)
  +\left(c_v^{\rm cen}
          -\sum_{\substack{m\in\mathcal I_v\\m\le N}}\frac4m\right)
       J_{N+1}(u,w)\right\}.
 \label{eq:st:green-limit}
\end{equation}
The error is
$O((1+\log(N+1))N^{-2}\|H''\|_\infty)$, with a constant depending
only on the real or complex type. For real $u,w$, each $E_m(u,w)$
is minus an arithmetic intersection of smooth metrized trivial
line bundles on $\mathbf P^1_\Z$. Complex profiles use sesquilinear
extension.
\end{theorem}
\begin{proof}
For $t\ne s$ the angle-dependent part of the Green kernel is
\[
 \sum_{k\ge1}\frac{e^{-k|t-s|/2}}k\cos(k(\theta-\phi)).
\]
Terms independent of the angles pair to zero. Angular integration
against $\cos(m\theta)\cos(m\phi)$ selects $k=m$ with factor $1/4$.
The factor $2\sqrt m$ in each measure therefore gives
\begin{equation}
 E_m(u,w)=\int_\R e^{-m|x|/2}H(x)\,dx.
 \label{eq:st:green-kernel}
\end{equation}
This calculation is justified first at unequal radii by the convergent
Fourier series and then by Fubini, since the original logarithmic
kernel is absolutely integrable against the smooth measures.

For the arithmetic interpretation, assume first that $u,w$ are real. Take
$dd^c=(i/2\pi)\partial\bar\partial$, so
$dd^c\log|z-z'|^2=\delta_{z'}$, and put
$U_{u,m}=2\int G(\,\cdot\,,z')\,d\mu_{u,m}(z')$.
The zero mass implies $dd^cU_{u,m}=-\mu_{u,m}$.
The potential is smooth by elliptic regularity; for real $u$ it
is invariant under complex conjugation. The Green star product
and arithmetic degree give
\[
 \widehat{\deg}\bigl((0,U_{u,m})\cdot(0,U_{w,m})\bigr)
 =\frac12\int_{\mathbf P^1(\C)}U_{u,m}\,dd^cU_{w,m}
 =-E_m(u,w)
\]
\cite[\S\S1.1.6(ii), 3.4.3, 4.2.11 and 4.3.1]{GilletSoule1990}.
Here $(0,U)$ is the trivial bundle with
$\|1\|^2=e^{-U}$ and its standard integral trivialization.
The complex case follows by sesquilinear extension.

It remains to identify the diagonal normalization. Averaging the
centered normal kernel over the two signs at a real place and over
the unit circle at a complex place gives the even densities
\[
 k_\R(x)=\frac{e^{-|x|/2}}{1-e^{-2|x|}},\qquad
 k_\C(x)=\frac{e^{-|x|/2}}{1-e^{-|x|}}
 \quad(x\ne0).
\]
For $x>0$ the real average is
$\frac12 e^{-x/2}((1-e^{-x})^{-1}+(1+e^{-x})^{-1})$.
For the complex average use
$\int_0^{2\pi}|1-re^{i\theta}|^{-2}d\theta/(2\pi)
=(1-r^2)^{-1}$ with $r=e^{-x/2}$.
Consequently
$k_v(x)=\sum_{m\in\mathcal I_v}e^{-m|x|/2}$.
The Fourier-normalized extension is
\begin{equation}
 W_v^0(\varphi_v)=c_v^{\rm cen}H(0)
       +\int_\R k_v(x)\bigl(H(x)-H(0)\bigr)\,dx.
 \label{eq:st:centered}
\end{equation}
To determine the scalar, use the homogeneous Fourier identity for the
reference character, valid for $0<\Re s<1$:
\[
 \mathcal F_{v,0}(|x|_v^{s-1})(y)
       =\gamma_v(s)|y|_v^{-s},\qquad
 \gamma_v(s)=\frac{\Gamma_v(s)}{\Gamma_v(1-s)}.
\]
Homogeneity and angular invariance determine the distribution up to
a scalar in this strip; there is no supported-at-zero term of that
degree. Evaluation on the self-dual Gaussian fixes the scalar by
its two Mellin integrals. For real $b$ with $\Re s-1<b<0$ this gives
the absolutely convergent Riesz convolution
\[
 [\mathcal F_{v,0}(|\cdot|_v^b)*|\cdot|_v^{s-1}](1)
       =\gamma_v(s)\gamma_v(1-s+b).
\]
One obtains the identity by Fourier transformation: the product is
$\gamma_v(s)|y|_v^{b-s}$, locally integrable at zero in the stated
range. Cutoffs give the distributional equality, and homogeneity
gives its value away from zero. Differentiation at $b=0$ is justified
by subtracting the test's value at zero near that singularity;
the other singularity and the tails, including their logarithmic
derivatives, remain integrable. Thus the value on the centered
mode $H(x)=e^{i\xi x}$, with $s=1/2+i\xi$, is
\[
 -\gamma_v(s)\gamma_v'(1-s)
 =-\frac{\Gamma_v'}{\Gamma_v}(s)
  -\frac{\Gamma_v'}{\Gamma_v}(1-s).
\]
This is real and even in $\xi$, so the centered distribution is even.
The first-order subtraction defining it gives order at most one at
zero. Its difference from the integrable first-difference extension
in \eqref{eq:st:centered} is therefore a multiple of $\delta_0$:
the only other supported term of that order is the odd distribution
$\delta_0'$. Taking $\xi=0$ gives
$-2\Gamma_v'(1/2)/\Gamma_v(1/2)=c_v^{\rm cen}$, proving
\eqref{eq:st:centered}. Constant and Fourier-mode tests are admissible
in this computation because the centered kernel decreases
exponentially at both ends.

Finally, evenness and Taylor's theorem give
\[
 \left|E_m(u,w)-\frac4mH(0)\right|
       \le\frac{16}{m^3}\|H''\|_\infty,\qquad
 |J_n(u,w)-H(0)|\le\frac4{n^2}\|H''\|_\infty.
\]
The linear Taylor term integrates to zero and
$\frac12\int_\R x^2e^{-m|x|/2}dx=16/m^3$.
The scalar-subtracted series therefore converges absolutely to the
integral in \eqref{eq:st:centered}. Its tail is $O(N^{-2}\|H''\|_\infty)$;
replacing $H(0)$ by $J_{N+1}$ in the finite scalar term costs
$O((1+\log(N+1))N^{-2}\|H''\|_\infty)$.
This proves \eqref{eq:st:green-limit}.
\end{proof}

\subsection{The combined formula}

For centered norm profiles, denote the finite expression
\eqref{eq:st:theta-pairing} by $C_v(u,w)$ and the finite expression
inside braces in \eqref{eq:st:green-limit} by $G_{v,N}(u,w)$.
The comparison theorem, the contact formula, and the Green formula give
\begin{equation}
 I_K(f,g)=W_K(\Phi)
 =\sum_{v\nmid\infty}C_v(u,w)
  +\lim_{N\to\infty}
     \left(\sum_{v\mid\infty}G_{v,N}(u,w)
                   -\log|D_K|\,J_{N+1}(u,w)\right).
 \label{eq:st:combined}
\end{equation}
The finite sum uses the existing cyclotomic contact lattices; the
infinite terms are signed metric intersections on the marked
coefficient line. The discriminant is the degree-zero theta defect
in \eqref{eq:st:global}.
On tests with both pole moments zero the middle trace is $-I_K$.

\appendix
\section{A dyadic expansion estimate}
\label{sec:dyadic-estimate}

Let
\[
 R(t)=\frac{\Phi_2(t)^2\Phi_8(t)}{\Phi_4(t)^3}.
\]
The period-one contact is measured on
$B_m=\Z_2[x]/x^m$, with $x=t-1$, inside $\Q_2[x]/x^m$.
Multiplication by $R(1+x)$ is invertible there, because $R(1)=1$.
Define the fractional lattice $L_m=R(1+x)^mB_m$.

\begin{definition}
For two commensurable $\Z_2$-lattices of the same rank, put
\[
 e_+(L\mid B)=\log\#(L/(L\cap B)).
\]
This is the expansion index relative to $B$.
\end{definition}

\begin{theorem}[Quadratic expansion of a cyclotomic jet]\label{thm:signed-cost}
The determinant of multiplication by $R(1+x)^m$ on $B_m\otimes\Q_2$
is one, but
\[
 \boxed{e_+(L_m\mid B_m)
 \ge\bigl(\tfrac18m^2-O(m\log m)\bigr)\log2.}
\]
\end{theorem}
\begin{proof}
Index the basis $1,x,\ldots,x^{m-1}$ by $0,\ldots,m-1$; below take
$m\ge2$, as suffices for the asymptotic statement.
Every such multiplication matrix is triangular with diagonal $R(1)^m=1$,
which proves the determinant statement. Factor, in $\Q_2[[x]]$,
\[
 R(1+x)=A(x)U(x),\quad A(x)=(1+x/2)^2,\quad
 U(x)=\frac{\Phi_8(1+x)/2}{(\Phi_4(1+x)/2)^3}.
\]
In the ring of integers of $\Q_2(\sqrt2)$, both $U(\sqrt2 y)$ and
its inverse have integral power series coefficients and constant term
one. Indeed
\[
 \Phi_4(1+\sqrt2y)/2=1+\sqrt2y+y^2,
\]
and
\[
 \Phi_8(1+\sqrt2y)/2
 =1+2\sqrt2y+6y^2+4\sqrt2y^3+2y^4.
\]
The same integrality holds for every positive or negative integral
power of $U$. Hence its coefficient of $x^j$ has valuation at least
$-j/2$. A $k$ by $k$ minor of its multiplication matrix on $B_m$
has valuation at least
\[
 -\tfrac12(\textstyle\sum\hbox{row indices}
                 -\sum\hbox{column indices})
 \ge-\tfrac12 k(m-k)\ge-\tfrac18m^2.
\]
Smith normal form says that $e_+(XB_m\mid B_m)/\log2$ is the largest
negative valuation of any minor of $X$, including the empty minor.
It follows that
\begin{equation}\label{eq:U-cost}
 e_+(U^{\pm m}B_m\mid B_m)\le\tfrac18m^2\log2.
\end{equation}

Set $k=\lfloor m/2\rfloor$. In the multiplication matrix of
$A^m=(1+x/2)^{2m}$, take rows $m-k,\ldots,m-1$ and columns
$0,\ldots,k-1$. Its determinant is
\[
 2^{-k(m-k)}
 \det\left[\binom{2m}{m-k+i-j}\right]_{0\le i,j<k}.
\]
The binomial determinant equals
\[
 \prod_{i=0}^{k-1}
 \frac{(2m+i)!\,i!}{(m-k+i)!\,(m+k+i)!}.
\]
For instance, the general formula with $2m$ replaced by $n$ and
$m-k$ by $d$ follows by induction from the adjacent-minor identity
\[
 D_k(n,d)D_{k-2}(n,d)
 =D_{k-1}(n,d)^2-D_{k-1}(n,d-1)D_{k-1}(n,d+1),
\]
with $D_0=1$ and $D_1=\binom nd$; the displayed factorial product
satisfies the same identity.
The formula $v_2(n!)=n-s_2(n)$ cancels the linear factorial terms.
The remaining binary digit sums have total absolute value
$O(m\log m)$. Therefore
\[
 e_+(A^mB_m\mid B_m)
 \ge\{k(m-k)-O(m\log m)\}\log2.
\]
Expansion indices satisfy the triangle inequality and are invariant
under simultaneous invertible linear maps. Since $A^m=R^mU^{-m}$,
\[
 e_+(A^mB_m\mid B_m)
 \le e_+(R^mB_m\mid B_m)+e_+(U^{-m}B_m\mid B_m).
\]
Using \eqref{eq:U-cost} and
$k(m-k)=\lfloor m^2/4\rfloor$ proves the result.
\end{proof}

\section*{Declaration}
ChatGPT was used in generating proofs.

\begingroup
\small
\bibliographystyle{amsplain}
\bibliography{Arithmetic_Winding_Curves}
\endgroup
\end{document}